\documentclass[11pt, a4paper]{article}

\usepackage[backref=page]{hyperref}

\usepackage[margin=2.5cm]{geometry}
\usepackage{mathtools}
\usepackage{amssymb}
\usepackage{amsthm}
\usepackage{amsmath}
\usepackage{empheq}

\usepackage{array,multirow,blkarray}

\usepackage[toc,page]{appendix}

\usepackage[nodayofweek]{datetime}

\usepackage{textcomp}

\usepackage{thmtools} 
\usepackage{thm-restate}

\usepackage{bbm}
\usepackage[]{dsfont}
\usepackage{authblk}
\usepackage{setspace}
\usepackage[shortlabels]{enumitem}
\usepackage[
 disable,
colorinlistoftodos, shadow, linecolor=white, backgroundcolor=white]{todonotes}

\usepackage{upgreek}

\newcommand{\Le}{\kappa}
\newcommand{\circuits}{\mathcal{C}}

\newcommand{\FeasibilityAlgorithm}{\textsc{Feasibility-Algorithm}}

\usepackage[linesnumbered,commentsnumbered,ruled,vlined]{algorithm2e}

\SetKw{run}{run}

\usepackage{bm}
\usepackage{xcolor}
\newcommand{\norm}[1]{\left\lVert#1\right\rVert}

\newcommand{\CPU}{\operatorname{\Psi}}
\newcommand{\nnz}{\operatorname{nnz}}

\usepackage{accents}
\newcommand{\ubar}[1]{\underaccent{\bar}{#1}}

\newenvironment{claimproof}[1][\proofname]
{
    \proof[#1 of Claim]
}
{
    \endproof
}

\newcounter{oraclecf}
\newcounter{algorithm saved}

\makeatletter
\newenvironment{oracle}[1][htb]{
    \RestyleAlgo{plainruled}
    \setlength{\algoheightrule}{1.2pt}

    \setcounter{algorithm saved}{\value{algocf}} \setcounter{algocf}{\value{oraclecf}}\begin{algorithm}[#1]}{\end{algorithm}
    \setcounter{oraclecf}{\value{algocf}}\setcounter{algocf}{\value{algorithm saved}}} 
\makeatother

\usepackage{cleveref}

\newcommand{\rk}{\operatorname{rk}}

\usepackage[noadjust]{cite}
\usepackage{booktabs}

\usepackage{algorithmic}

\newcommand{\im}{\operatorname{Im}}

\newcommand{\Diag}{\operatorname{diag}}

\hypersetup{
    colorlinks,
    linkcolor={red!50!black},
    citecolor={blue!50!black},
    urlcolor={blue!80!black}
}

\newtheorem{theorem}{Theorem}[section]

\newtheorem{lemma}[theorem]{Lemma}
\newtheorem{proposition}[theorem]{Proposition}
\newtheorem{remark}[theorem]{Remark}
\newtheorem{corollary}[theorem]{Corollary}

\newtheorem{claim}{Claim}[theorem]

\newtheoremstyle{case}{}{}{\itshape}{1em}{}{:}{ }{}
\theoremstyle{case}

\definecolor{darkred}{RGB}{180, 0, 0}
\definecolor{darkgreen}{rgb}{0, 0.6, 0}
\definecolor{darkblue}{RGB}{51,51,178}
\definecolor{lightgray}{RGB}{231,231,231}
\definecolor{lightblue}{RGB}{180,180,254}
\definecolor{lightred}{HTML}{FEB4B4}
\definecolor{darkcyan}{HTML}{7FBFBF}

\newcommand{\LPPrimalDual}{\textsl{Primal-Dual}}
\newcommand{\FeasLP}{\textsl{Feas-LP}}
\newcommand{\ProxFeas}{\textsl{Prox-Feas}}
\newcommand{\ProxFeasSolver}{\textsl{Prox-Feas-Solver}}

\newcommand{\ProxOpt}{\textsl{Prox-Opt-Solver}}

\newcommand{\InnerSysF}{\textsl{F-Primal}}

\newcommand{\InnerLoop}{\textsc{InnerLoop}}

\newcommand{\cl}{\operatorname{cl}}
\newcommand{\1}{\mathbbm{1}}

\newcommand{\supp}{\mathrm{supp}}

\DeclareMathOperator*{\argmin}{arg\,min}

\newcounter{myquestion}

\newcommand{\bl}[1]{{ #1}}

\newcounter{mycomment}

\renewcommand{\phi}{\varphi}
\renewcommand{\rho}{\varrho}

\newcommand{\R}{\mathbb{R}}
\newcommand{\Z}{\mathbb{Z}}

\newcommand{\eps}{\varepsilon}

\newcommand{\set}[1]{\left\{ #1 \right\}}
\newcommand{\pr}[2]{\langle #1, #2 \rangle}

\newcounter{Hequation}

\makeatletter
\g@addto@macro\equation{\stepcounter{Hequation}}
\makeatother

\author[1]{Daniel Dadush}
\author[2]{Bento Natura}
\author[2]{L{\'{a}}szl{\'{o}} A. V{\'{e}}gh}
\affil[1]{Centrum Wiskunde \& Informatica, The Netherlands}
\affil[2]{Department of Mathematics, London School of Economics and Political Science}

\title{Revisiting Tardos's Framework for Linear Programming: \\ Faster Exact Solutions
using Approximate Solvers
\thanks{This project has received funding from the European Research Council (ERC) under the European Union's Horizon 2020 research and innovation programme (grant agreements ScaleOpt--757481 and QIP--805241).
}
}
\date{\today}

\makeatletter
\newenvironment{tagequation}[1][Missing Tag]{
    \begin{equation}
    \tag{#1}
    \begin{aligned}
    \hphantom{\text{#1}}
    \end{aligned}
    \begin{aligned}
}{
    \end{aligned}
    \end{equation}
}
\makeatother

\usepackage{etoolbox}

\newtoggle{focs}
\togglefalse{focs}

\let\lIf\If
\let\lElse\Else
\let\lCase\uCase

\begin{document}
\maketitle

\begin{abstract}
    In breakthrough work, Tardos (Oper. Res. '86) gave a \emph{proximity} based
    framework for solving linear programming (LP) in time depending only on the
    constraint matrix in the bit complexity model. In Tardos's framework, one reduces
    solving the LP $\min \pr{c}{x}$, $Ax=b$, $x \geq 0$, $A \in \Z^{m \times n}$, to
    solving $O(nm)$ LPs in $A$ having small integer coefficient objectives and
    right-hand sides using \emph{any exact} LP algorithm. \bl{This gives rise to an LP algorithm in time poly$(n,m\log\Delta_A)$, where $\Delta_A$ is the largest subdeterminant of $A$. A significant extension to the real model of computation was given by Vavasis and Ye (Math. Prog. '96), giving a specialized interior point method that runs in time poly$(n,m,\log\bar\chi_A)$, depending on Stewart's $\bar{\chi}_A$, a well-studied condition number. 

    In this work, we extend Tardos's original framework to obtain such a running time dependence. In particular, we replace
    the exact LP solves with \emph{approximate} ones, enabling us to
    directly leverage the tremendous recent algorithmic progress for approximate
    linear programming. More precisely, we show that the fundamental ``accuracy''
    needed to \emph{exactly} solve any LP in $A$ is inverse polynomial in $n$ and
    $\log\bar{\chi}_A$. Plugging in the
    recent algorithm of van den Brand (SODA '20), our method computes an optimal
    primal and dual solution using ${O}(m n^{\omega+1} \log
    (n)\log(\bar{\chi}_A+n))$ arithmetic operations, outperforming the specialized
    interior point method of Vavasis and Ye  and its recent
    improvement by Dadush et al (STOC '20).   By applying the preprocessing algorithm
    of the latter paper, the dependence  can also be reduced from
    $\bar\chi_A$ to $\bar{\chi}_A^*$, the minimum value of $\bar\chi_{AD}$
    attainable via column rescalings. 
    Our framework is applicable to achieve the
    poly$(n,m,\log\bar\chi_A^*)$ bound using essentially any weakly
    polynomial LP algorithm, such as the ellipsoid method.}
    
    At a technical level, our framework combines together approximate LP solutions
    to compute exact ones, making use of constructive proximity theorems---which
    bound the distance between solutions of ``nearby'' LPs---to keep the required
    accuracy low.

    \end{abstract}
 
\newpage
\listoftodos
\newpage

\section{Introduction}
In this paper, we consider the task of computing exact primal and dual solutions
for linear programs (LP) in standard form:
\begin{equation}
\label{LP_primal_dual} \tag{LP}
\begin{aligned}
\min \; &\pr{c}{x} \quad \\
Ax& =b \\
x &\geq 0\, ,\\
\end{aligned}
\quad\quad\quad
\begin{aligned} 
\max \; & \pr{y}{b} \\
A^\top y + s &= c \\
s & \geq 0\, . \\
\end{aligned}
\end{equation}
Here, $A\in \R^{m\times n}$, ${\rm rank}(A) = m \leq n$, $b\in \R^m$, $c\in
\R^n$ are given in the input, and $x,s\in\R^n$, $y\in \R^m$ are the variables.
We consider the program in $x$ to be the primal problem and the program in $y,s$
to be the dual problem.

After the work of Khachiyan~\cite{Khachiyan79}, who gave the first polynomial
algorithm for LP using the ellipsoid method, Megiddo~\cite{Megiddo83} asked
whether there exists a ``genuinely polynomial'', now known as \emph{strongly
polynomial}, algorithm for LP. Informally, the goal is to find an algorithm that uses
${\rm poly}(n)$ \emph{basic arithmetic operations} (e.g.~addition,
multiplication, etc.), where each such operation must be performed on numbers of
size polynomial in the instance encoding length. While no such algorithm is
known, the search for a strongly polynomial LP algorithm has spurred tremendous
algorithmic advances for many classical combinatorial problems. 

Strongly polynomial algorithms have indeed been found for important
combinatorial classes of linear programs. Examples include
feasibility for two variable per
inequality systems~\cite{Megiddo83}, minimum-cost circulations~\cite{Goldberg89,Orlin93,Tardos85}, generalized flow maximization,
\cite{Vegh17,OV17}, and discounted Markov Decision Processes \cite{Ye2005,Ye2011}.

To generalize these results to larger problem classes, a natural
attempt is to seek abstract frameworks that capture known algorithms. In this
vein, a recurring principle in strongly polynomial algorithm design is that
``good enough'' approximate solutions can be used to glean combinatorial
information about exact optimal ones. Such information is used to
reduce the underlying instance in a way that preserves all optimal solutions. 

This was in fact the key idea in Tardos's seminal paper on minimum-cost circulations \cite{Tardos85}: solving a problem instance with a suitable rounded cost function reveals an arc that cannot be tight in any dual optimal solution; consequently, we can fix the flow value to 0.
 As another example, in
submodular function minimization any sufficiently small norm point in the base
polytope can be used to infer relations in a ring-family containing all
minimizers~\cite{IFF,DVZ18}. 

At a higher level, it can be useful to view
strongly polynomial algorithms as reductions from an exact optimization problem
to a suitable approximate version of itself. To achieve fast strongly polynomial
algorithms using these principles, important considerations are the complexity
of the individual approximate solves, e.g., the degree of accuracy required,
and the total required number of them.

\paragraph{Tardos's Framework for Linear Programming} 
Generalizing the above idea from minimum-cost flows to general linear programming,
Tardos~\cite{Tardos86} provided such a framework for
solving any standard form primal-dual LP with integer constraint matrix $A \in
\Z^{m \times n}$ using a number operations depending only on $n$ and the
logarithm of $\Delta_A$, the maximum absolute value of the determinant of any
square submatrix of $A$. This algorithm is strongly polynomial for
minimum-cost flow, \bl{noting that digraph incidence matrices are totally unimodular, and therefore  $\Delta_A = 1$.} At a high
level, Tardos's framework reduces getting exact LP solutions to getting exact
solutions for ``nearby LPs'' with simpler coefficient structure, heavily relying on
LP \emph{proximity} theorems (e.g., see~\cite{Hoffman52,Cook1986}). \bl{More
precisely, Tardos kmputing exact primal-dual solutions to $\max
~\pr{c}{x}, Ax=b, x \geq 0$ to computing exact primal-dual solutions to $O(nm)$
LPs in $A$ with ``rounded'' objectives $c'$ and right hand sides $b'$ having
integer coefficients of size $O(n^2 \Delta_A)$.}
In particular, after $O(n)$
such LP solves, one can determine a coefficient $x_i$ in the support of some
optimal solution, allowing to delete the $x_i \geq 0$ constraint. 
Due to their small coefficients, the LPs in the reduction can be solved using any
weakly polynomial algorithm. We note that the fundamental property enabling the
polynomial solvability of these rounded LPs is that the minimum non-zero slack
of their basic solutions, i.e., $\min \{x_i: x_i > 0\}$, is lower bounded by
$1/(n^{O(1)}\Delta)$ by Cramer's rule. 

\paragraph{Achieving $\bar{\chi}_A$ dependence} 
\bl{While Tardos's framework is powerful, it inherently relies on the 
condition measure $\Delta_A$. This is only defined for integer constraint matrices; one can obtain bounds for rational constraint matrices via multiplying by the least common denominator of the entries, but this leads to weak bounds that are highly volatile under small changes in the entries.}
\bl{
A significant strengthening of \cite{Tardos86} was given by 
Vavasis and Ye~\cite{Vavasis1996}. 
They gave an
interior point method (IPM) in the real model of computation based on \emph{layered least squares} (LLS) steps that
outputs exact primal-dual solutions in $O(n^{3.5} \log(\bar{\chi}_A+n))$
iterations.} Improved iteration bounds were later given for certain special
cases, in particular, $O(\sqrt{n}\log(\bar{\chi}_A+n))$ for homogeneous conic
feasibility~\cite{VavasisYeRealNumberData} and $O(n^{2.5}\log(\bar{\chi}_A+n))$
for LP feasibility~\cite{Ye:2006}. In a conceptual advance, Vavasis and Ye's
result showed that the polynomial solvability of LP does not require any minimum
non-zero slack assumption. 

The condition measure replacing $\Delta_A$ is Stewart's
$\bar{\chi}_A$~\cite{stewart}, which for integer matrices satisfies $\bar{\chi}_A \leq n \Delta_A$. \bl{In contrast with $\Delta_A$ that relies on the entry numerics, $\bar\chi_A$ is a geometric measure that depends only on the kernel of $A$;} 
Formally, letting $W
:= \ker(A)$ and $\pi_I(W) = \{x_I : x \in W\}$, one may define $\bar{\chi}_A :=
\bar{\chi}_W$ as the minimum number $M \geq 1$ such that for any $\emptyset \neq
I \subseteq [n]$ and $z \in \pi_I(W)$, there exists $y \in W$ with $y_I = z$ and
$\|y\| \leq M \|z\|$. In words, it represents the cost of lifting partial
fixings of coordinates into the subspace $W$. 

Very recently, the authors and Huiberts~\cite{DHNV20}, building on the work of
Monteiro and Tsuchiya~\cite{Monteiro2003, MonteiroT05}, gave an improved LLS
optimization algorithm and analysis requiring only $O(n^{2.5} \log n
\log(\bar{\chi}^*_A+n))$ iterations, where $\bar{\chi}^*_A$ is the minimum
$\bar{\chi}_{AD}$ over positive diagonal matrices $D > 0$. \cite{DHNV20} further
gave a nearly optimal rescaling algorithm which runs in $O(m^2n^2 + n^3)$ time
and computes $D > 0$ satisfying $\bar{\chi}_{AD} \leq n (\bar{\chi}^*_{A})^3$.
Thus, by suitable preprocessing, any algorithm achieving $\bar{\chi}_A$
dependence can be converted into one with $\bar{\chi}^*_A$ dependence.

\bl{A key tool in \cite{DHNV20} is to study the `circuit imbalance measure' $\kappa_A$. This closely approximates $\bar\chi_A$, with $\log(\bar{\chi}_A+n)=\Theta(\log(\kappa_A+n))$, and has very favourable combinatorial properties. Our approach also relies on $\kappa_A$ and $\kappa^*_A$, even though we state the results in terms of the better known
$\bar\chi_A$ and $\bar\chi^*_A$.

The condition number $\bar{\chi}^*_A$  can be smaller than $\bar\chi_A$ by an arbitrary factor, and in turn, $\bar\chi_A$ can be much smaller  $\Delta_A$ even for integer matrices $A$. Let $A\in \R^{n\times m}$ be the node-edge incidence matrix of an undirected graph on $n$ nodes and $m$ edges. If the graph has $k$ node-disjoint odd cycles, then $\Delta_A \ge 2^{k}$. However, it is easy to verify that for any graph, $\kappa_A\le 2$ (see the definition of $\kappa_A$ in Section~\ref{sec:chi-kappa}). Using Proposition~\ref{prop:kappa-chi}, we get the bound $\bar\chi_A\le 2m$.}

\bl{
\paragraph{Harnessing the progress in approximate solvers}
The complexity of fast approximate LP algorithms has seen substantial improvements in recent years~\cite{LS19,CLS19,vdb20,br2020solving, jiangFasterDynamicMatrix2020}.}
Taking the recent algorithm~\cite{vdb20}, given a
feasible LP $\min~\pr{c}{x}, Ax=b, x \geq 0$, having an optimal solution of
$\ell_2$ norm at most $R$, for $\eps > 0$ it computes a point $\tilde{x} \geq 0$
satisfying
\begin{equation}
\label{eq:apx-lp} \tag{APX-LP}
\iftoggle{focs}{
\begin{aligned}
	\pr{c}{\tilde{x}} &\leq \min_{Ax=b,x\geq 0} \pr{c}{x} + \eps \cdot \|c\|_2 \cdot R
\quad \\ \|A\tilde{x}-b\|_2 &\leq \eps \cdot(\|A\|_F\cdot R +
\|b\|_2),
\end{aligned}
}{
\pr{c}{\tilde{x}} \leq \min_{Ax=b,x\geq 0} \pr{c}{x} + \eps \cdot \|c\|_2 \cdot R
\quad \text{ and } \quad \|A\tilde{x}-b\|_2 \leq \eps \cdot(\|A\|_F\cdot R +
\|b\|_2),
}
\end{equation}
in deterministic time $O(n^\omega \log^2(n) \log(n/\eps))$, where $\omega < 2.38$ is the matrix multiplication exponent. 

Tardos's framework requires an exact black box solver for systems with the same matrix $A$ but replacing $b$ and $c$ by small integer vectors. It is possible to use the approximate solver \eqref{eq:apx-lp} to obtain exact optimal solution for integer matrices for  sufficiently small $\varepsilon$. \bl{ Assume $A\in \mathbb{Z}^{m\times n}$, $b\in \Z^m$, $c\in \Z^n$ and $\|b\|_\infty,\|c\|_\infty\le n^{O(1)}\Delta^t$, and let $\textrm{OPT}$ denote the optimum value of \eqref{LP_primal_dual}. 
We may call \eqref{LP_primal_dual} in a suitable extended system,  $\varepsilon=1/\left(n^{O(1)}\Delta_A^{O(t)}\right)$, and use a Carath\'eodory reduction to identify primal and dual optimal basic solutions. Integrality is used in multiple parts of such a reduction: e.g., for establishing a bound $R=n^{O(1)}\Delta_A^{O(t)}$ from Cramer's rule, and for showing that for any primal feasible solution $x$, $\pr{c}{x}<\textrm{OPT}$ implies $\pr{c}{x}<\textrm{OPT}-\varepsilon\|c\|_2 R$.}
 For a matrix $A\in \mathbb{R}^{m\times n}$, we cannot obtain an
exact solver by applying the approximate solver for high enough accuracy in terms of the condition numbers $\bar\chi_A$ or $\kappa_A$. This is the main reason why we cannot work with explicitly rounded systems, but require a more flexible approach.
\bl{Let us also note that recovering an exact solution from the approximate solver comes at a high arithmetic cost that we can save if using the approximate solution directly.}

\bl{
\paragraph{Fast algorithms with $\bar\chi_A$ dependence}
The layered least squares interior point methods discussed above 
represent substantial advances in the strongly polynomial
solvability of LP, yet it is highly non-obvious how to combine these techniques with
those of recent fast LP solvers.} For example, for the results
of~\cite{LS19,br2020solving}, one would have to develop analogues of LLS steps
for weighted versions of the logarithmic barrier. Furthermore, the proofs
of exact convergence are intricate and deeply tied to the properties of the
central path, and may leave one wondering whether the $\bar{\chi}_A$ solvability
of LP is due to ``IPM magic''. 
It  would therefore be desirable to have an
elementary proof of the $\bar{\chi}_A$ solvability of LP.

Partial progress on this question was given by Tun{\c{c}}el and
Ho~\cite{ho2002}, who generalized Tardos's framework in the real number model.
Firstly, they showed that one can still round instances to have minimum non-zero
slack $\tau_A > 0$, depending only on $A$. Second, they showed that applying the
Mizuno-Todd-Ye~\cite{MTY} predictor-corrector IPM on the homogeneous self-dual
formulation, these rounded instances can be solved ${\rm poly}(n,\log
\tau_A,\log (\Delta_A/\delta_A))$ time, where $\delta_A$ is the absolute value
of the minimum non-zero determinant of any square submatrix of $A$. Here, they
prove the relation $\bar{\chi}_A \leq n \Delta_A/\delta_A$ and note that
$\Delta_A/\delta_A$ can be arbitrarily larger than $\bar{\chi}_A$. Lastly, they
provide a different algorithm that removes the dependence on $\tau_A$, assuming
one has access to the Vavasis-Ye algorithm as a subroutine only on instances
with $b \in \set{\pm 1,0}^m, c \in \{0, \pm 1\}^n$.   

\subsection{Our Contributions}\label{sec:contributions}

As our main contribution, we provide a substantially improved Tardos style
framework for LP which achieves both $\bar{\chi}_A$ dependence and relies only
on approximate LP solves: \bl{we use the output  \eqref{eq:apx-lp}  of the approximate LP solvers in a black-box manner.}
 Our main result is summarized below. The more precise
technical statements are given as \Cref{thm:feas} for feasibility and
\Cref{thm:optimization-main} for optimization. \bl{The system \eqref{initialization_LP} is an extended system used for initialization, defined in Section~\ref{sec:initialization}.}

\begin{theorem}[Enhanced Tardos Framework for Feasibility]
\label{thm:tardos-v2-feas}
	 Assume we are given a feasibility LP $Ax=b$, $x\ge 0$ with data $A \in \R^{m \times n}$, ${\rm
rank}(A)=m$, and $b \in \R^m$.
\begin{enumerate}[label=(\roman*)]
	\item \label{it:primal-feas} If the primal program is feasible, then one can find a feasible solution $x$ using $O(m)$ approximate LP solves~\eqref{eq:apx-lp} with accuracy
$\eps = 1/(n\bar{\chi}_A)^{O(1)}$, on extended systems of the
form~\eqref{initialization_LP}, together with additional $O(mn^\omega
)$ arithmetic operations. This gives a total complexity $O(m n^{\omega}\log^2(n)\log(\bar{\chi}_A+n))$ using  the solver of van den Brand~\cite{vdb20}.
\item \label{it:primal-feas-cert} If the primal program is infeasible, then a Farkas certificate of infeasibility $y \in \R^m$, satisfying $A^\top y \geq
0$, $\langle b, y
\rangle < 0$ can be found using the amount of computation as in \ref{it:primal-feas}, and $O(nm^2 + n^\omega)\log\log(\bar\chi_A+n))$ additional arithmetic operations. 
\end{enumerate}
\end{theorem}
Next, we state our result for optimization:

\begin{theorem}[Enhanced Tardos Framework for Optimization]
\label{thm:tardos-v2-opt}
	 Assume we are given primal-dual~\eqref{LP_primal_dual} with data $A \in \R^{m \times n}$, ${\rm
rank}(A)=m$, $b \in \R^m$, $c \in \R^n$.
\begin{enumerate}[label=(\roman*)]
\item \label{it:opt-feas} If both primal and dual programs are feasible, then one can obtain
an optimal primal-dual pair $(x,y,s)$ of solutions, using at most $O(nm)$ approximate LP solves~\eqref{eq:apx-lp} as in \Cref{thm:tardos-v2-feas}\ref{it:primal-feas}, together with an additional $O(mn^{\omega + 1}
)$ arithmetic operations. This gives a total complexity $O(m n^{\omega+1}\log^2(n)\log(\bar{\chi}_A+n))$ using \cite{vdb20}.
\item \label{it:opt-feas-cert} If either of the primal or dual programs are infeasible, then we can obtain a Farkas certificate of primal or dual infeasibility in the same running time as in \ref{it:opt-feas}, plus $O(n^3 m^2\log\log(\bar\chi_A+n))$ additional arithmetic operations. 
\end{enumerate}
\end{theorem}

This theorem yields the first LP algorithm achieving $\bar{\chi}_A$ dependence that is
not based of the analysis of the central path. At a high level, we achieve this
by more deeply exploiting the power of LP proximity theorems, which are already
at the core of Tardos's framework. In the rest of this section, we explain some
of the key ideas behind the above theorem and how it compares to Tardos's original
algorithm as well as that of Vavasis and Ye. 

\paragraph{Overview of the approach}
Both Tardos's and our approach use variants of Hoffman's proximity bounds, see Section~\ref{sec:proximity-hoffman}. 
 The fundamental difference is that while Tardos uses an exact solver where the 
perturbed objective and right hand side vectors are fixed in
advance before calling the solver, we decide these perturbations ``on the fly'' as a function of the returned approximate solutions we receive. 

Let us illustrate Tardos's and our approaches  on the dual feasibility LP 
\begin{equation}\label{eq:dual}\tag{$D$}
A^\top y+s=c, s\ge 0\, .
\end{equation}
 The feasibility algorithm in \cite{Tardos85} proceeds as follows. Define $\tilde b_i=\sum_{i=1}^n (\Delta_A+1)^{i-1}a_i$, where $a_i$ is the $i$-th column vector of $A$, and
consider the primal system 
\begin{equation}\label{eq:primal-t}\tag{$\tilde P$}
\min\ \pr{c}{x}\quad\textrm{s.t. }Ax=\tilde b\, ,\ x\ge 0\, .
\end{equation}
Note that by the choice of $\tilde b$, this system is always feasible. If it is unbounded, then we may conclude infeasibility of \eqref{eq:dual}.
\bl{The reason for the particular choice of $\tilde b$ is that whenever the system is bounded, the dual of \eqref{eq:primal-t}  has a unique optimal solution; this can be shown by 
a determinant argument. Consequently, for any optimal solution $x^*$ to  \eqref{eq:primal-t} and $S^*=\supp(x^*)$, the system $a_i^\top y=c_i$, $i\in S^*$ yields a feasible solution to \eqref{eq:dual}.}
The exact LP solver will be applied to a series of rounded problem instances of the form 
\begin{equation}\label{eq:primal-h}\tag{$\hat P$}
\min\  \pr{\tilde c}{x}\quad\textrm{s.t. }Ax=\tilde b\, \, ,\ x\ge 0\, ,\ x_T=0\, ,
\end{equation}
 where $\tilde c\in \Z^n$, $\|\tilde c\|_\infty\le n^2 \Delta_A$, and $T\subseteq [n]$ is a set of indices $i$ where we have already concluded that $x^*_i=0$ in every optimal solution to \eqref{eq:primal-t}. This is initialized as $T=\emptyset$, and every call to the LP solver enables the addition of  at least one new index; thus, we need $O(n)$ oracle calls to solve feasiblity. According to the definition of $\tilde b$, this is an integer vector with \bl{$\|\tilde b\|=\Theta(\sqrt{m}\Delta_A^n)$}. As explained above, we can obtain an exact solution to \eqref{eq:primal-h} by calling \eqref{eq:apx-lp} for accuracy $\varepsilon=1/\left(n^{O(1)}\tilde \Delta_A^{O(n)}\right)$.

 To conclude that $i\in T$ for some $i\in T$, Tardos uses a proximity theorem that is a variant of \Cref{lem:primal_prox}. It implies that if $\|\tilde c-c\|_\infty$ is ``small'', then \eqref{eq:primal-t} has a dual optimal solution that is ``close'' to the dual optimal solution obtained for \eqref{eq:primal-h}.

\medskip

 In contrast, our approach in Section~\ref{sec:feasibility} proceeds as follows. If $c\ge 0$, we simply return $s=c$. Otherwise, the norm of the negative coordinates $
 \|c^-\|_1$ will play a key role. We can
strengthen \eqref{eq:dual} by adding the  constraint 
 \begin{equation}\label{eq:prox-constraint}
 \|s-c\|_\infty \le 16\kappa_A^2 n\|c^-\|_1\, ,
 \end{equation}
  where $\kappa_A$ is the circuit imbalance measure; for integer matrices $\kappa_A\le \Delta_A$. A proximity result (\Cref{cor:feas}) implies that whenever \eqref{eq:dual} is feasible, there is a feasible solution also satisfying \eqref{eq:prox-constraint}.

We can use \eqref{eq:apx-lp} directly to obtain a solution $(\tilde y,\tilde s)$ such that $A^\top \tilde y+\tilde s=c$, $\|\tilde s-c\|_\infty\le 3\kappa_A^2 n\|c^-\|_1$, and $\|x^-\|_\infty \le \varepsilon \|c^-\|_1$ for $\varepsilon=1/O(n^4\kappa_A^4)$. Again, note that in addition to approximate feasiblity, we also require proximity of $s$ to $c$; \bl{we can obtain such a solution with this extra property without an increase in the running time cost.}

From here, we can identify a set $K$ of coordinates such that $\tilde s_i$ is large enough to conclude that there exists a feasible solution $s$ to \eqref{eq:dual} with $\tilde s_i>0$ for $i\in K$; this is done similarly as in Tardos's approach.

 We project out all variables in $K$, meaning that we remove the inequalities $a_i^\top y+s_i=c_i$ for $i\in K$ from the system. We recurse on the smaller subsystem. From the recursive call, we obtain a feasible solution $y'$ to \eqref{eq:dual} in the smaller system that also satisfies \eqref{eq:prox-constraint}. The proximity constraints enables us to easily map back $y'$ to a feasible solution $y$ to \eqref{eq:dual} by a simple `pullback' operation.

\medskip

As noted above, the very existence of an exact LP oracle heavily relies on the integrality assumption of $A$. This integrality is also used to establish the relation between the optimal solutions of \eqref{eq:primal-t} and the solutions of \eqref{eq:dual}, using a determinant argument.
\bl{ In contrast, the proximity arguments as in \Cref{lem:primal_prox} and \Cref{cor:feas} does not rely on integrality; we can use here $\kappa_A$ instead of $\Delta_A$.}

Even for integer matrices and $\kappa_A=\Theta(\Delta_A)$, and using the same solver for \eqref{eq:apx-lp}, our algorithm is faster by a factor $\Omega(n^2/m)$. 
A key ingredient in the running time improvement is to strengthen the system with \eqref{eq:prox-constraint}. This allows us to use $\varepsilon=1/(n^{O(1)}\kappa_A^{O(1)})$; otherwise, we would need to require a higher precision $\varepsilon=1/(n^{O(1)}\kappa_A^{O(n)})$. This yields a factor $n$ improvement over \cite{Tardos85}. 

Another factor $n/m$ improvement can be obtained as follows. In the approach sketched above, if the set of ``large'' coordinates $K$ is nonempty, we get a bound $n$ on the number of recursive calls. Using a slightly more careful recursive setup, we can decrease the rank of the system  at each iteration, improving this bound to $m$.

\medskip

Let us now turn to optimization. Our algorithm will be more similar to the one in \cite{Tardos85}, and for integer matrices with $\kappa_A=\Theta(\Delta_A)$ \bl{ and $m=\Omega(n)$}, the asymptotic running time bounds will be the same. 

We now outline Tardos's approach. Given an optimization problem \eqref{LP_primal_dual}, we first check for both primal and dual feasibility. If these are both feasible, then we go through $\le m$ main loops. In each main loop, we use the same approach as above to solve \eqref{eq:primal-t} with a perturbed $\tilde b\in \Z^m$ with $\|\tilde b\|_\infty\le n^2\Delta_A$. Using $\le n$ oracle calls, we obtain optimal primal and dual solutions $(x,y,s)$ Again, proximity guarantees that if $\tilde b$ is ``close'' to $b$, then we can identify an index $i$ with a ``large'' $x_i>0$ where we can conclude $s^*_i=0$ in every optimal solution. Equivalently, $x_i$
is in the support of some optimal solution, and hence we may delete the
constraint $x_i \geq 0$, and proceed to the next main loop after projecting out the variable $x_i$. We note that the bound $n$ on the inner loops is in reality $n-m$, and this can be improved to $m$ by swapping the primal and dual sides.

In our approach in Section~\ref{sec:optimization}, the goal is to end up in the same place as Tardos at the end of
the main loop, where the difference will be how we get there. As mentioned
above, in Tardos's setting, one already knows beforehand that the final
objective and right hand side for which one will have optimal primal-dual
solutions will be $\tilde{b}$, a rounded version of $b$, and the original
$c$. However, the only important property is that at the end of the loop we end
up with a primal-dual optimal pair for the original objective $c$, and
\emph{some} right hand side $b'$ close enough to the original $b$. In
particular, $b'$ need not be known at the beginning of the algorithm and can
thus be chosen \emph{adaptively} depending on the outcome of the approximate LP
solves. 

For the above purpose, we utilize proximity theorems (see
Section~\ref{sec:proximity-hoffman} for precise statements) to allow us to stitch
together the ``large'' coordinates of approximate dual solutions to achieve
feasibility. At the same time, we perform a similar complementary stitching of
primal approximate solutions, where we judiciously perturb ``small'' coordinates
to $0$, inducing a corresponding change of right hand side, to enforce
complementarity with the dual solution. Here proximity allows us to control how
much the solutions will change in future iterations, which is crucial to not
destroying the structure of the solutions built so far. 

\medskip

We also note that Gr\"otschel, Lov\'asz, and Schrijver \cite[Theorem 6.6.3]{gls} give a different proof for Tardos's result using simultaneous Diophantine approximation (see also \cite{franktardos}). This shows that \eqref{LP_primal_dual} can be solved by creating a single perturbed instance with integer $\tilde b$ and $\tilde c$ bounded in terms of the encoding length of $A$ such that the set of optimal bases coincide in the two systems. The perturbed instance can be solved in poly$(n,m,\log \Delta_A)$; we simply take an optimal basis and compute the corresponding primal and dual optimal solutions for the original $b$ and $c$.
\bl{However, this reduction inherently relies on integrality arguments.}

\iftoggle{focs}{}{
\paragraph{Comparison to layered least squares IPM methods} To setup a comparison, we first
recall that standard log-barrier based IPMs follow the \emph{central path}
$\{(x(\mu),s(\mu),y(\mu)): \mu > 0\}$ defined by the equations $x_i(\mu)
s_i(\mu) = \mu$, $\forall i \in [n]$, together with feasibility
for~\eqref{LP_primal_dual}. $\mu$ represents the normalized duality gap and
$(x_\mu,y_\mu,s_\mu)$ converges to optimal solutions as $\mu \rightarrow 0$.
The number of calls to the approximate LP solver above can be usefully
compared to the number of so-called disjoint \emph{crossover events} on the
central path used in the analysis of the Vavasis--Ye algorithm \cite{Vavasis1996}. A crossover
event occurs for a pair of distinct indices $(i,j)$ between the times $\mu^1 <
\mu^0$, if $x_i(\mu^0) \bar{\chi}_A^n \geq x_j(\mu^0)$ and for all times $\mu' <
\mu^1$, $ x_j(\mu^0) > x_i(\mu^0)$. In words, an $(i,j)$ crossover happens
between time $\mu^0$ and $\mu^1$ if the variables $x_i,x_j$ are ``close'' to
being in the wrong order at time $\mu^0$ and are in the correct order at all
times after $\mu^1$. The Vavasis and Ye LLS step was in fact designed to ensure
that a new cross-over event occurs a ``short time'' after the step,
i.e., sometime before $\mu/\bar{\chi}_A^n$ if the step ends at $\mu$. From here,
it is obvious that the number of distinct crossover events, i.e., on a new pair
of indices, is bounded by $\binom{n}{2}$. 

The approximate LP solves in our algorithm have the effect of inducing similar
crossover type events, though this number is $O(mn)$ instead of $O(n^2)$.
Precisely, after each LP solve, we are able identify two non-empty disjoint
subsets of variables $I,J \subseteq [n]$, such that at least one of the
variables $x_j, j \in J$, will end up being substantially larger than all the
variables $x_i, i \in I$ in the final optimal solution. Lastly, the accuracy
requirement of $\eps = 1/(n\bar{\chi}_A)^{O(1)}$ for each LP solve is in a sense
analogous to moving down the central path by that amount. We note
that~\cite{DHNV20} gave an improved analysis of the Vavasis and Ye algorithm,
showing that on ``average'' one sees $\Omega(1/\log n)$ (slightly different)
crossover events after $(n\bar{\chi}_A)^{O(1)}$ time units, which is slightly
worse than what we achieve here per approximate LP solve.}

\paragraph{Failure will be certified} Our algorithm requires an estimate on the circuit imbalance parameter $\kappa_A$ (see definition in Section~\ref{sec:chi-kappa}). This is a common assumption shared by most previous literature: Tardos's algorithm uses an estimate of $\Delta_A$; Vavasis and Ye require a bound on $\bar\chi_A$. These parameters are hard to compute \cite{Khachiyan:1995:CAE:192503.192510,Tuncel1999}. However, knowing these values are not required, and we can use the following simple guessing procedure, attributed to J. Renegar in \cite{Vavasis1996}. We start with a low guess on $\bar\chi_A$ (or some other parameter), say $M=100$. If the algorithm fails to return the required solution, then we conclude that the estimate was too low, and replace the guess $M$ by $M^2$. Thus, we can still obtain a dependence on $\log(\bar\chi_A+n)$, without knowing the value.

A new aspect of our algorithm is that in case of a failure, we do not simply conclude that our estimate was too low indirectly from the failure of the algorithm, but we also obtain an explicit \emph{certificate}. Namely, an elementary operation is to compute \emph{lifts} mentioned previously: for the subset $W=\ker(A)$, an index set $I\subseteq [n]$, and a vector $y\in \pi_I(W)$, we compute the minimum-norm vector $z\in W$ such that $z_I=y$. Our parameter $\kappa_A$ satisfies $\|z\|_\infty\le \kappa_A\|y\|_1$ (\Cref{lem:kappa-lift}). Whenever our algorithm fails due to underestimating $M<\kappa_A$, this will be certified by an index set $I\subseteq [n]$ and a vector $y\in \pi_I(W)$, and lift $z$ with  $\|z\|_\infty> M\|y\|_1$.

\iftoggle{focs}{}{

\subsection{Organization}

In Section~\ref{sec:prelims}, we recall the subspace formulation of LP and
review important properties of the condition numbers $\bar{\chi}_A$ and its
combinatorial cousin, the circuit imbalance measure $\kappa_A$. In Section~\ref{sec:proximity-hoffman}, we give a
self-contained review of known LP proximity results, based on Hoffman type
bounds, as well as novel variant (Theorem~\ref{thm:dual-fix}), which leverages
the structure of the linear independence matroid on $A$. In
Section~\ref{sec:proximity-alg}, we present a constructive strongly polynomial
time variant of Hoffman's proximity theorem, which will be useful for extracting
Farkas infeasibility certificates from approximate solutions. In
Section~\ref{sec:blackbox}, we review the current state of the art approximate
LP solvers and state our main theorems for extracting the solutions we need from
these solvers in both the feasibility and optimization context, Theorems
\ref{thm: best_LP_algorithms_real_model_feasibility} and \ref{thm:
best_LP_algorithms_real_model} respectively. The proofs of these theorems are
deferred to Section~\ref{sec:approx}, where we also describe the
LP extended system we use~\eqref{new_initialization_LP_subspace}. In Section~\ref{sec:feasibility}, give the
describe our framework for LP feasibility, and in Section~\ref{sec:optimization}
our framework for LP optimization.

}

\section{Preliminaries}
\label{sec:prelims}

\iftoggle{focs}{}{
We let $[n] := \{1, \ldots, n\}$.
Let $\R_{++}$ denote the set of positive reals, and $\R_+$ the set of nonnegative reals. We denote the support of a vector $x \in \R^n$ by $\supp(x) = \{i\in [n]: x_i \neq 0\}$.  We let $\1_n$ denote the $n$-dimensional all-ones vector, or simply $\1$, whenever the dimension is clear from the context. Let $e_i$ denote the $i$-th unit vector.
}

For vectors $v, w \in \R^n$ we denote by $\min\{v,w\}$ the vector $z \in \R^n$ with $z_i = \min\{v_i, w_i\}, i \in [n]$; analogously for $\max\{v,w\}$. Further, we use the notation $v^+ = \max\{v, 0_n\}$ and $v^- = \max\{-v, 0_n\}$ 
\iftoggle{focs}{.}{;note that both $v^+$ and $v^-$ are nonnegative vectors.

We will use $\ell_1,\ell_2$ and $\ell_\infty$ vector norms, denoted as $\|.\|_1,\|.\|_2$, and $\|.\|_\infty$, respectively. By $\|v\|$, we always mean the 2-norm $\|v\|_2$. Further, for a matrix $A\in\R^{m\times n}$, $\|A\|$ will refer to the $\ell_2\to\ell_2$ operator norm, $\|A\|_F=\sqrt{\sum_{i,j}|A_{ij}|^2}$ to the Frobenius-norm, and $\|A\|_{\max}=\max_{i,j}|A_{ij}|$ to the max-norm.
}

For a vector $v \in \R^n$, we denote by $\Diag(v)$ the diagonal matrix whose $i$-th diagonal entry is $v_i$.
For two vectors $x, y \in \R^n$, we let $\pr{x}{y}=x^\top y$ denote their scalar product. We denote by the binary operation $\circ$ the element-wise multiplication $x \circ y = \Diag(x)y$.  We let $\mathbf D$ denote the set of all positive definite $n\times n$
diagonal matrices.

For an index subset $I\subseteq [n]$, we use $\pi_I: \R^n \rightarrow \R^I$ for the coordinate
projection. That is, $\pi_I(x)=x_I$, and for a subset $S\subseteq
\R^n$, $\pi_I(S)=\{x_I:\, x\in S\}$.
We let $\R^n_I = \{x \in \R^n : x_{[n]\setminus I} = 0\}$.

For a subspace $W\subseteq \R^n$, we let $W_I=\pi_I(W\cap \R^n_I)$. It is easy to see that $\pi_I(W)^\perp = (W^\perp)_I$. 
Assume we are given a matrix $A\in \R^{m\times n}$ such that $W=\ker(A)$. Then, $W_I=\ker(A_I)$, and we can obtain a matrix $A'$ from $A$ such that $\pi_I(W)=\ker(A')$ by performing a Gaussian elimination of the variables in $[n]\setminus I$.

For a subspace $W \subseteq \R^n$ and a vector $d \in \R^n$ we define by $d/W$ the orthogonal projection of $d$ onto $W^\perp$, that is $d/W = \pi_{W^\perp}(d)$. In particular, $d/W$ is the minimum-norm vector in $W + d$.
Further, for a subspace $W \subseteq \R^n$, we let
$W_+=W\cap \R^n_+$.

\paragraph{Linear programming in subspace formulation}
Let $A\in\R^{m\times n}$, and $W=\ker(A)\subseteq \R^n$. 
For $c,d\in \R^n$,
we can write  \eqref{LP_primal_dual} in the following equivalent form, where $d\in \R^n$ such that $Ad=b$.
\begin{tagequation}[\LPPrimalDual$(W,d,c)$]\label{LP-subspace-f}
\min \; & \pr{c}{x} \\ 
x &\in W + d \\
x &\geq 0\, ,
\end{aligned}
\quad\quad
\begin{aligned}
\max \; &  \pr{d}{c-s} \\ 
s  &\in W^\perp+c \\
s &\geq 0\, .
\end{tagequation}
Note that $(x,s)$ are optimal primal and dual solutions if and only if they are feasible and $\pr{x}{s}=0$.
Thus, \ref{LP-subspace-f} is equivalent to the 
following feasibility problem:
\begin{equation}\label{LP-subspace}
  \iftoggle{focs}{x \in W + d, \,s \in W^\perp + c, \,\pr{x}{s} = 0, \,(x, s) \ge 0.}
  {
\begin{aligned}
x &\in W + d \\
s  &\in W^\perp+c \\
\pr{x}{s}&=0\\
x, s &\geq 0\, .
\end{aligned}
  }
\end{equation}
\paragraph{Circuits}
For a linear subspace $W \subseteq \R^n$ and a matrix $A$ such that
$W = \ker(A)$, a circuit is an inclusion-wise
minimal dependent set of columns of $A$. 
\iftoggle{focs}{}{Recall that this corresponds to standard notion of circuits in the linear matroid
associated with  $A$. 

This notion only depends on the subspace $W$, and not on the particular representation $A$; an equivalent definition is that 
$C \subseteq [n]$ is a circuit if and only if $W \cap \R^n_C$ is one-dimensional and that no strict subset of $C$
has this property.}
The set of circuits of $W$ is denoted $\circuits_W$.

For a subset $I\subseteq [n]$, we let $\cl(I)$ denote its \emph{closure} in the matroidal sense. 
\iftoggle{focs}{
  We will make the assumption that $|C|>1$ for all $C\in \circuits_W$. 
}{
That is, $\cl(I)=J$ is the unique maximal set containing $J\supseteq I$ such that $\rk(A_J)=\rk(A_I)$. Equivalently,
\[
\cl(I)=I\cup\{j\in [n]\setminus I:\, \exists C\in \circuits_W, j\in C\subseteq I\cup \{j\}\}
\]

We will make the assumption that 
\begin{equation}\label{a:no-loop}
|C|>1\quad \forall C\in \circuits_W\, ,
\end{equation}
that is, $W$ does not contain any loops (trivial circuits). 

If $W=\ker(A)$, then a loop corresponds to zero-columns of $A$. For solving the feasibility problem $x\in W+d, x\ge 0$, we can eliminate all variables $i\in [n]$ that form  a loop, without affecting feasibility. For the optimization problem \ref{LP-subspace-f}, if $i\in [n]$ forms a loop, and the primal problem $x\in W+d, x\ge 0$ is feasible, then $c_i<0$ means that the problem is unbounded, and if $c_i\ge 0$, then we can find an optimal solution with $x_i=0$.
}

\iftoggle{focs}{}{
\paragraph{Sign-consistent circuit decompositions}
We say that the vector $y \in \R^n$ is \emph{sign-consistent} with
$x\in\R^n$ if $x_iy_i \geq 0$ for all $i \in [n]$ and $x_i = 0$ implies $y_i = 0$
for all $i \in [n]$. Given a subspace $W\subseteq \R^n$, a \emph{sign-consistent circuit decomposition} of a vector $z\in W$ is a decomposition
\[
z=\sum_{k=1}^h  g^k,
\]
where
$h\le n$, the vectors $g^1,g^2,\ldots,g^h\in W$ are sign-consistent with
$z$ and $\supp(g^k)=C_k$ for all $k\in [h]$ for circuits
 $C_1,C_2,\ldots,C_h\in
\circuits_W$.
}

\begin{lemma} \label{lem:sign-cons}
For every subspace $W\subseteq \R^n$,  every $z\in W$ admits a sign-consistent circuit decomposition. 
\end{lemma}
\iftoggle{focs}{}{\begin{proof}
Let $F \subseteq W$ be the set of vectors sign-consistent with $z$. $F$ is a polyhedral cone; its
 faces  correspond to inequalities of the form $y_k \geq 0$, $y_k\leq 0$,
or $y_k = 0$.
The rays (edges) of $F$ are of the form $\{\alpha g:\,  \alpha\ge 0\}$  with $\supp(g) \in \circuits_W$. Clearly, $z\in F$, and thus, $z$ can be written as a conic combination of at most $n$ rays by the Minkowski-Weyl theorem. Such a decomposition yields a sign-consistent circuit decomposition.
\end{proof}}

\subsection[The condition numbers chi-bar and kappa]{The condition numbers {$\bar\chi$} and $\kappa$}\label{sec:chi-kappa}
For a matrix $A\in \R^{m\times n}$, the condition number $\bar\chi_A$ is defined as
\begin{equation}
\begin{aligned}
\bar\chi_A&=\sup\left\{\left\|A^\top \left(A D A^\top\right)^{-1}AD\right\|\, : D\in
    {\mathbf D}\right\}
    \iftoggle{focs}{}{
\\
&=\sup\left\{\frac{\norm{A^\top y}}{\norm{p}} :
\text{$y$ minimizes $\norm{D^{1/2}(A^\top y - p)}$ for some $0 \neq p \in \R^n$ and $D \in \mathbf D$}\right\}
}
.
\end{aligned}
\end{equation}
This quantity was first studied by  Dikin \cite{dikin}, Stewart \cite{stewart}, and
Todd \cite{todd-90}, and  has been extensively studied in the context of interior point methods; we refer the reader to 
\cite{ho2002,Monteiro2003,Vavasis1996} for further results
and references.

It is important to note that $\bar\chi_A$ only depends on the subspace
$W=\ker(A)$. Hence, we can also write $\bar\chi_W$ for a subspace
$W\subseteq \R^n$, defined to be equal to $\bar\chi_A$ for some matrix
$A\in \R^{k\times n}$ with $W=\ker(A)$. We will use the notations $\bar\chi_A$ and
$\bar\chi_W$ interchangeably. 

\iftoggle{focs}{}{

\begin{proposition}[{\cite{Todd2001,Vavasis1996}}]\label{prop:chi-bar-norm}
For every matrix $A\in \R^{m\times n}$,
\[
\bar\chi_A = \max\left\{ \|A_B^{-1} A\| : A_B \text{ non-singular $m \times m$-submatrix of } A\right\}\, .
\]
\end{proposition}
}
Let us define the {\em lifting map}
$L_I^W : \pi_{I}(W) \to W$ by
\[
L_I^W(p) = \argmin\left\{\|z\| : z_I = p, z \in W\right\}.
\]
\iftoggle{focs}{}{
Note that $L_I^W$ is the unique linear map from $\pi_{I}(W)$ to $W$ such that
$L_I^W(p)_I = p$ and $L_I^W(p)$ is orthogonal to $W \cap \R^n_{[n] \setminus
I}$. If $|I| = \Theta(n)$ any explicit computation of the lifting map requires $\Omega(n^2)$ operations as $L_I^W \in \R^{n\times |I|}$. The following lemma shows that the lift of a single vector can be obtained more efficiently if $m \ll n$.

\begin{lemma}
\label{lem:compute-lift}
Let $A \in \R^{m \times n}$, $W = \ker(A)$ and $I \subseteq [n]$. 
Then, both of the following can be done in time $O(\min\{m^2n, n^\omega\})$. 
\begin{enumerate}[label=(\roman*)]
  \item \label{it:primal_lift_comp} Computing $L_I^W(p)$, for any $p \in \pi_I(W)$. 
  \item \label{it:dual_lift_comp} Computing $L_I^{W^\perp}(q)$, for any $q \in \pi_I(W^\perp)$. 
\end{enumerate}
\end{lemma}
}
\iftoggle{focs}{}{\begin{proof}
Let us begin with \ref{it:primal_lift_comp}. First, obtain a vector $w \in W$ such that $w_I = p$. This can be done by solving the linear system $A_{[n]\setminus I} x = -A_I p$ and setting $w = (x,p)$ in time $O(\min\{m^2n, n^\omega\})$. The components in $[n]\setminus I$ of $L_I^W(p)$ are now given by the orthogonal projection of $x$ onto $(W_{[n]\setminus I})^\perp$. Note that $W_{[n]\setminus I} = \ker(A_{[n]\setminus I})$. So, for $B := A_{[n]\setminus I}$ we have 
\[[L_I^W(p)]_{[n]\setminus I} = B^\top(BB^\top)^{-1}Bx,\]
where $BB^\top \in \R^{m \times m}$ can be computed in $O(\min\{m^2n, n^\omega\})$ and its inverse in $O(m^\omega) \le O(\min\{m^2n, n^\omega\})$. All other operations are matrix-vector products and can be computed in time $O(mn)$. Therefore, the overall running time of computing $L_I^W(p)$ is $O(\min\{m^2n, n^\omega\})$.

To show\ref{it:dual_lift_comp},  analogously to \ref{it:primal_lift_comp} we first obtain a vector $\hat w \in W^\perp$ such that $\hat w_I = q$. Note that $W^\perp = \im(A^\top)$ and so we find $y \in \R^m$ such that $(A_I)^\top y = q$. This can again be done in time $O(\min\{m^2n, n^\omega\})$. Then, we can set $\hat w = A^\top y \in W^\perp$. The coordinates in $[n]\setminus I$ of $L_I^{W^\perp}(q)$ are now given by by the projection of $\hat w_{[n]\setminus I}$ onto $\pi_{[n]\setminus I}(W)$. The space $\pi_{[n]\setminus I}(W)$ can be represented as the kernel of a matrix $\hat A$, that arises by performing Gaussian elimination and pivoting the entries in $I$ in time $O(nm^{\omega - 1}) \le O(\min\{m^2n, n^\omega\})$ \cite{Bunch1974,Ibarra1982}. Then, we have
\[ 
[L_I^{W^\perp}(q)]_{[n]\setminus I} = I - \hat A^\top (\hat A \hat A^\top)^{-1}\hat A \hat w_{[n]\setminus I}.
\]
As in \ref{it:primal_lift_comp} this can be computed in time $O(\min\{m^2n, n^\omega\})$.
\end{proof}}

A useful characterization of $\bar\chi_W$ can be given in terms of the operator norm of the lifting map.
This is was shown in \cite{DHNV20}, by using results from \cite{stewart} and  \cite{OLeary1990}.

\begin{proposition}[\cite{DHNV20}]\label{prop:subspace-chi}
For a linear subspace $W  \subseteq \R^n$,
\[
\bar{\chi}_W =\max\left\{ \|L_I^W\|\, : {I\subseteq [n]}, I\neq\emptyset\right\}\, .
\]
\end{proposition}

\paragraph{The circuit imbalance measure}
Consider now the circuits $\circuits_W$ of the subspace $W$. For a circuit $C \in \circuits_W$, let $g^C \in W$ be such that $\supp(g^C) = C$. Note that $g^C$ is unique up to multiplication by scalar. 
We define the {\em circuit imbalance measure}
\[
\Le_W=\max\left\{\frac{\max_{i\in C} \left|g^C_i\right|}{\min_{j\in C} \left|g^C_j\right|}:\, C\in \circuits_W\right\}\, .
\]
as the largest ratio between two entries of any minimum support nonzero vector in $W$. This was studied in 
\cite{Vavasis1994,ho2002,DHNV20}. Note that $\Le_W=1$ corresponds to totally unimodular spaces. As shown in \Cref{prop:kappa-chi} below, the condition 
measures $\bar\chi_W$ and $\Le_W$ are closely related: $O(\log(\bar\chi_W+n))=O(\log(\kappa_W+n))$ holds.
However, $\kappa_W$ has several advantageous combinatorial properties. It fits particularly nicely with the proximity results in Section~\ref{sec:proximity-hoffman} and \ref{sec:proximity-alg}. In fact, the argument in the proof of Tardos's main proximity result using Cramer's rule is implicitly bounding circuit imbalances, see discussion of \Cref{lem:primal_prox}. Therefore, we will use $\kappa_W$ instead of $\bar\chi_W$ throughout the paper.

We can give a characterization analogous to \Cref{prop:chi-bar-norm}, using max-norm instead of $\ell_2$-norm. \begin{proposition}\label{prop:kappa-max}
For every matrix $A\in \R^{m\times n}$ with $\rk(A) = m$,
\[
  \iftoggle{focs}{
\Le_A = \max\left\{ \|A_B^{-1} A\|_{\max} : B \text{ basis}\right\}\, .
  }{
    \Le_A = \max\left\{ \|A_B^{-1} A\|_{\max} : A_B \text{ non-singular $m \times m$-submatrix of } A\right\}\, .
  }
\]
\end{proposition} 
\iftoggle{focs}{}{\begin{proof}
Consider the matrix $A'=A_B^{-1}A$ for any non-singular $m\times m$ submatrix $A_B$. Let us renumber the columns such that $B$ corresponds to the first $m$ columns. 
Then, for every $m+1\le j\le n$, the $j$th column of $A'$ represents a circuit where $x_j=1$, and $x_i=-A'_{ij}$ for $i\in [m]$. Hence, $\|A'\|_{\max}$ gives a lower bound on $\Le_A$. 

Conversely, take the circuit solution $g^C$ that gives the maximum in the definition of $\Le_A$; let $g_j^C$ be the minimum absolute value element.
For any circuit $C$, we can select a base $B$ such that $C\setminus \{j\}\subseteq B$. Then, the largest absolute value in the $j$-th column of $A_B^{-1}A$ will be $\Le_A$.
\end{proof}}

\Cref{prop:subspace-chi} asserts that $\bar\chi_W$ is the maximum $\ell_2\to\ell_2$ operator norm of the mappings $L_I^W$ over $I\subseteq [n]$. We now show that  the maximum $\ell_1\to\ell_\infty$ operator norm of these mappings is  $\kappa_W$, even though the lifting map is defined with respect to $\ell_2$ norms.

\begin{proposition}\label{lem:kappa-lift}
For a linear subspace $W  \subseteq \R^n$,
\[
\kappa_W =\max\left\{ \frac{\|L_I^W(p)\|_\infty}{\|p\|_1}\, : {I\subseteq [n]}, I\neq\emptyset, p\in \pi_I(W)\setminus\{0\}\right\}\, .
\]
\end{proposition}

\iftoggle{focs}{}{\begin{proof}
We first show that for any $I\neq\emptyset$, and $p\in \pi_I(W)\setminus\{0\}$, $\|L_I^W(p)\|_\infty\le \kappa_W \|p\|_1$ holds. 

Let $z=L_I^W(p)$, and take a sign-consistent decomposition $z=\sum_{k=1}^h g^k$ as in \Cref{lem:sign-cons}. For each $k\in [h]$, let $C_k=\supp(g^k)$. We claim that all these circuits must intersect $I$. Indeed, assume for a contradiction that one of them, say $C_1$ is disjoint from $I$, and let $z'=\sum_{k=2}^h g^k$. Then, $z'\in W$ and $z'_I=z_I=p$. Thus, $z'$ also lifts $p$ to $W$, but $\|z'\|_2<\|z\|_2$, contradicting the definition of $z=L_I^W(p)$ as the minimum-norm lift of $p$.

By the definition of $\kappa_W$, $\|g^k\|_\infty\le \kappa_W \|g^k_I\|_1$ for each $k\in [h]$. The claim follows since $p=z_I=\sum_{k=1}^h g^k_I$, moreover, sign-consistency guarantees that $\|p\|_1=\sum_{k=1}^h \|g^k_I\|_1$. Therefore,
\[
\|z\|_\infty\le \sum_{k=1}^h \|g^k\|_\infty\le \kappa_W \sum_{k=1}^h \|g^k_I\|_1=\kappa_W\|p\|_1\, .\]
We have thus shown that the maximum value in the statement is at most $\kappa_W$. To show that equality holds, let $C\in\circuits_W$ be the circuit and $g^C\in W$ the corresponding vector and $i,j\in C$ such that $\kappa_W=|g^C_i|/|g^C_j|$.

Let us set $I=([n]\setminus C)\cup \{j\}$, and define $p_k=0$ if $k\in [n]\setminus C$ and $p_j=g^C_j$. Then  $p\in \pi_I(W)$,  and the unique extension to $W$ is $g^C$; thus, $L_I^W(p)=g^C$. We have $\|L_I^W(p)\|_\infty=|g^C_i|$. Noting that $\|p\|_1=|g^C_j|$, it follows that 
$\kappa_W=\|L_I^W(p)\|_\infty/\|p\|_1$.
\end{proof}}

  Using \Cref{lem:kappa-lift}, we can easily relate the quantities $\bar\chi_W$ and $\kappa_W$. The upper bound was already shown in \cite{Vavasis1994}; the slightly weaker lower bound $\sqrt{\bar\chi_W^2-1}/n\le \kappa_W$ was given in \cite{DHNV20}.

  \begin{proposition}\label{prop:kappa-chi}
      For a linear subspace $W \subseteq \R^n$,
      \[\frac{1}{n} \bar\chi_W \le \kappa_W \le \sqrt{\bar\chi_W^2 - 1}.\]
  \end{proposition}
  \iftoggle{focs}{}{\begin{proof}
  
    Note that the inequality $\frac{1}{n} \bar\chi_W \le \kappa_W \le \bar\chi_W$ is a simple consequence from \Cref{prop:subspace-chi}, \Cref{lem:kappa-lift} and bounds between the $\ell_\infty$, $\ell_1$, and $\ell_2$-norms.
     Let us now turn to the second inequality. Note that by assumption \eqref{a:no-loop}, \Cref{prop:subspace-chi} implies $\bar\chi_W \ge \sqrt{2}$. So the inequality already holds if $\kappa_W = 1$. For the rest, assume that $\kappa_W > 1$. We pick an element $p \in \R^I, p \in \pi_I(W)$ for which $\kappa_W$ attains its value. Then for $J = [n] \setminus I$
    \[\kappa_W = \frac{\|L_I^W(p)\|_\infty}{\|p\|_1} = 
    \frac{\|{L_I^W(p)}_J\|_\infty}{\|p\|_1} 
    \le
     \frac{\|{L_I^W(p)}_J\|_2}{\|p\|_2} = \sqrt{\frac{\|L_I^W(p)\|_2^2}{\|p\|_2^2} - 1} \le \sqrt{\bar\chi_W^2 - 1},\]
     where the second equality used $\|{L_I^W(p)}_J\|_\infty > \|{L_I^W(p)}_I\|_\infty = \|p\|_\infty$, as $\kappa_W > 1$.
  \end{proof}}

  \iftoggle{focs}{}{

We state two more useful properties of  $\bar\chi_W$ and $\Le_W$. 
The first property shows that  $\bar\chi_W$ and $\Le_W$ are self-dual, shown respectively in \cite{gonzaga_lara} and in \cite{DHNV20}.
\begin{proposition}[\cite{gonzaga_lara,DHNV20}]\label{prop:kappa-dual}
For any linear subspace $W\subseteq \R^n$, we have $\bar\chi_W=\bar\chi_{W^\perp}$ and 
$\Le_W=\Le_{W^\perp}$.
\end{proposition}

The next proposition asserts that these condition numbers are non-increasing under fixing and projecting  variables.

\begin{proposition}\label{prop:proj-fix}
For any linear subspace $W\subseteq \R^n$ and $J\subseteq[n]$, we have 
\[
\bar\chi_{W_J}\le \bar\chi_{W}\, ,\quad \bar\chi_{\pi_J(W)}\le \bar\chi_{W}\, ,\quad \Le_{W_J}\le\Le_W\, ,\quad \mbox{and}\quad \Le_{\pi_J(W)}\le\Le_W\, .
\]
\end{proposition}
\iftoggle{focs}{}{\begin{proof}
The statements on $W_J$ are immediate from \Cref{prop:subspace-chi} and \Cref{lem:kappa-lift}. The two other statements follow from the statements on $W_J$, \Cref{prop:kappa-dual} and noting that $\pi_J(W)^\perp = (W^\perp)_J$.
\end{proof}}
  }

\paragraph{The estimate $M$ and lifting certificates} The value of $\kappa_W$ and $\bar\chi_A$ may not be known. In fact, these are hard to approximate even within a factor $2^{\mathrm{poly}(n)}$ \cite{Tuncel1999}. Throughout our algorithms, we maintain  a guess $M$ on the value of $2\kappa_W$, initialized as $M=2$. At certain points in the algorithm,  we may find an index set $I\subseteq [n]$ and a vector $p\in \pi_I(W)$  such that
$\|L_I^W(p)\|_\infty>M \|p\|_1$. In this case, we conclude that $M<\kappa_W$ by \Cref{lem:kappa-lift}. Such a pair $(I,p)$ is called a \emph{lifting certificate} of $M>\kappa$.
We can then restart the algorithm with an updated estimate $M'=\max\{2\|L_I^W(p)\|_\infty/ \|p\|_1,M^2\}$. 

\iftoggle{focs}{}{
\begin{remark} \label{remark:lift_carry_over} 
  During the algorithm, we project out variable sets $J \subset [n]$ and work recursively with the space $W':= \pi_{[n]\setminus J}(W)$. A lifting certificate for $W'$ is then a pair $(I,p)$ with $I \subset [n]\setminus J$, $p \in \R^I$, such that $\|L_I^{W'}(p)\|_\infty > M \|p\|_1$. While \Cref{prop:proj-fix} already certifies that the guess $M$ is wrong for $W$, it is unsatisfactory that the obtained certificate holds for a different space. But it is easy to see that the certificate still holds up to a factor of $\sqrt{n}$ also for the original space $W$:
Let $\hat p \in W$ be an arbitrary vector such that $p_{[n]\setminus I} = L_I^{W'}(p)$. Then
\[
  \sqrt{n}\|L_I^W(p)\|_\infty \ge \|L_I^W(p)\| \ge \|[L_I^W(p)]_{[n]\setminus J}\| \ge \|L_I^{W'}(p)\| \ge \|L_I^{W'}(p)\|_\infty > M \|p\|_1.
\]
In particular, the inequality above shows $\|L_I^W(p)\| > M\|p\|_1 \ge M\|p\|$, so $(I,p)$ is a certificate for $W$ in the classical $\ell_2$-norm.
For ease of presentation we disregard this detail in the remainder of the paper.
\end{remark}
}

\paragraph{Optimal rescalings}
For every $D\in {\mathbf D}$, we can consider the condition numbers
$\bar\chi_{WD}=\bar\chi_{AD^{-1}}$ and  $\Le_{WD}=\Le_{AD^{-1}}$. We let
\[
\begin{aligned}
\bar\chi^*_W&=\bar\chi^*_A=\inf\{\bar\chi_{WD}\, : D\in {\mathbf D}\}\, \\
\Le^*_W&=\Le^*_A=\inf\{\Le_{WD}\, : D\in {\mathbf D}\}\, \\
\end{aligned}
\]
denote the best possible values of $\bar\chi$ and $\Le$ that can be attained by
rescaling the coordinates of $W$. A near-optimal rescaling can be found in strongly polynomial time
\iftoggle{focs}{\cite{DHNV20}.}{.
\begin{theorem}[{\cite{DHNV20}}]\label{thm:bar-chi-star}
There is an $O(n^2m^2 + n^3)$ time algorithm that for any matrix $A\in\R^{m\times n}$ and $W=\ker(A)$, computes a value $t$  such that
\[
\begin{aligned}
t &\leq \bar\chi_W \leq t (\bar\chi_W^*)^2\\
\end{aligned}
\]
and a $D\in {\mathbf D}$ such that
\[
\begin{aligned}
\Le_{WD}\le (\Le_W^*)^3\, \mbox{ and } \bar\chi_{WD}\le n(\bar\chi_W^*)^3\, .
\end{aligned}
\]
\end{theorem}
}
As a consequence, after using this preprocessing step, any algorithm that has running time dependence on $\log (\kappa_W+n)$ is turned into an algorithm with dependence on $\log(\kappa_W^*+n)$. We note however that for small values of $\log (\kappa_W+n)$, this preprocessing may turn out to be a bottleneck operation for our feasibility algorithm.

\section{Proximity via Hoffman-bounds}
\label{sec:proximity-hoffman}

Hoffman's seminal work \cite{Hoffman52} has analyzed proximity of LP solutions. Given $P=\{x\in \R^n:\, Ax\le b\}$, $x_0\in \R^n$, and norms $\|.\|_\alpha$ and $\|.\|_\beta$, we are interested in the minimum of 
$\|x-x_0\|_\alpha$ over $x\in P$.
Hoffman showed that this can be bounded as $H_{\alpha,\beta}(A)\|(Ax_0-b)^+\|_\beta$, where the Lipschitz-bound $H_{\alpha,\beta}(A)$ is a constant that only depends on $A$ and the norms. Such bounds have been shown for several different problem forms and norms; we refer the reader to \cite{Pena2020} for results and references. 

We will use a Hoffman-bound for a system of the form $x\in W$, $\ell\le x\le u$.
We show that $H_{\infty,1}=\kappa_W$ for such a system. Related bounds using $\bar\chi_A$ have been shown in \cite{ho2002}; here, we present a self-contained proof.

For vectors $d,c\in \R^n$, let us define the set
\begin{equation}\label{theta-def}
\Lambda(d,c):=\supp(d^-)\cup \supp(c^+)\, .
\end{equation}

\begin{theorem}[Hoffman Proximity Theorem]\label{thm:master-proximity}
Let $W \subseteq \R^n$ be a subspace and $\ell\in (\R\cup\{-\infty\})^n$,
$u \in (\R\cup\{\infty\})^n$ be lower and upper bounds, and assume that $P = \{x \in W: \ell \le x \le u\}$ is non-empty. Then, 
for every $x\in P$ we have
\[
\|\ell^++ u^-\|_1 \le \|x_{\Lambda(u,\ell)}\|_1\, ,
\]
and there exists $x \in P$ such that 
\[
\|x\|_\infty \le \kappa_W \|\ell^++ u^-\|_1\, . \]  
\end{theorem}

\iftoggle{focs}{}{\begin{proof}
Let us start with the first statement. Clearly, $\supp(u^-)\cap \supp(\ell^+)=\emptyset$.
If $u_i<0$, then $|x_i|\ge |u_i|$, and if $\ell_i>0$, then $|x_i|\ge |\ell_i|$. Thus,  $\|\ell^+ + u^-\|_1\le \|x_{\Lambda(u,\ell)}\|_1$ follows for every $x\in P$.

For the second statement, select $x \in P$ such that $\|x\|_1$ minimal and let $x = \sum_{k=1}^h g^k$ be a sign-consistent circuit decomposition of $x$ as in \Cref{lem:sign-cons}. For each $k\in [h]$, we show that $C_k=\supp(g^k)$ must either contain an element $i\in C_k$ with $x_i=u_i<0$, or with $x_i=\ell_i>0$. For a contradiction, assume that one of them, say $C_1$, contains no such element. Then, for some $\varepsilon>0$, $x'=(1-\varepsilon) g^1+ \sum_{k=2}^h g^k\in P$, giving a contradiction, since $\|x'\|_1<\|x\|_1$.

The inequality $\|x\|_\infty \le \kappa_W \|\ell^++ u^-\|_1$ then follows as in the proof of \Cref{lem:kappa-lift}.
\end{proof}}

We can derive  useful corollaries for feasibility and optimization problems.
\begin{corollary}\label{cor:feas}
Let $W\subseteq \R^n$ be a subspace and $d\in \R^n$. If the system
$x\in W+d,\, x\ge 0$ is
 feasible, then the system
 \iftoggle{focs}{
   \[
  x\in W+d,\;
\|x-d\|_\infty\le \kappa_W \|d^-\|_1,\;
x\ge 0 ,
\]
 }{
\begin{align*} 
x&\in W+d\\
\|x-d\|_\infty&\le \kappa_W \|d^-\|_1\\
x&\ge 0,
\end{align*}
 }
is also feasible.
\end{corollary}
\iftoggle{focs}{}{\begin{proof}
Using the variable $z=d-x$, the system can be equivalently written as $z\in W$, $z\le d$. Thus, \Cref{thm:master-proximity} guarantees a solution with $\|z\|_\infty\le \kappa_W\|d^-\|_1$, as required.
\end{proof}}

\begin{corollary}\label{thm:prox-basic}
Let $W\subseteq \R^n$ be a subspace and $c,d\in \R^n$, and assume $c\ge 0$.
If \ref{LP-subspace-f} is feasible, 
then there is an optimal  solution $(x,s)$ such that 
\[
\|x-d\|_\infty\le \Le_W \|d_{\Lambda(d,c)}\|_1\, .
\]
\end{corollary}
\iftoggle{focs}{}{\begin{proof}
Let $(x^*, s^*)$ be in \ref{LP-subspace-f}. That is, $x^*$ minimizes $\pr{c}{x}$ over $x\in W+d, x\ge 0$. Consider the feasibility system $x\in W+d, x\ge 0$, and $x_i\le x_i^* \ \forall i\in\supp(c)$. Note that $x^*$ is a feasible solution. In fact, the inequality $x_i\le x_i^*$ must be tight for all $i\in\supp(c)$, as otherwise $\pr{c}{x}<\pr{c}{x^*}$, contradicting the optimality of $x^*$.

For $z=d-x$, we get the system $z\in W$, $\ell\le z\le d$, where $\ell_i=d_i-x_i^*$ if $i\in \supp(c)$, and $\ell_i=-\infty$ if $c_i=0$. Note that $\ell^+_i\le d^+_i$, for $i\in \supp(c)=\supp(c^+)$, and therefore,
\[
\|d_{\Lambda(d,c)}\|_1=\|d^+_{\supp(c)}\|_1+\|d^-\|_1\ge \|\ell^+ +d^-\|_1\, .
\]
Thus, the claim follows by \Cref{thm:master-proximity}.
\end{proof}}

The next lemma will be used to conclude that a primal variable $s^*_i=0$ in every solution $(x^*,s^*)$ to \ref{LP-subspace-f}. For integer matrices, a similar statement was given by Cook et al. \cite[Theorem 5]{Cook1986}, see also \cite[Theorem 10.5]{SchrijverLPIP} with a bound in terms of the maximum subdeterminant $\Delta_A$. A variant of this statement is used by Tardos \cite[Lemma 1.1]{Tardos85} as the main underlying proximity statement of her algorithm. Ho and Tun\c{c}el \cite[Theorem 6.3]{ho2002} generalized this bound to arbitrary matrices, using the condition number $\bar \chi_A$. This implies our statement with $n\kappa_W$ instead of $\kappa_W+1$. We note that the arguments in \cite{Cook1986,Tardos85} are based on Cramer's rule. In essence, this is used to bound the circuit imbalances in terms of $\Delta_A$. Hence, our formulation with $\kappa_W$ can be seen as a natural extension.

\begin{lemma}\label{lem:primal_prox}
Let $W\subseteq \R^n$ be a subspace and $c, d, \tilde d\in \R^n$.
Let $(\tilde x, s)$ be an optimal solution to
 \hyperref[LP-subspace-f]{\LPPrimalDual$(W,\tilde d,c)$}. Then there exists an optimal solution $(x^*, s^*)$ to {\LPPrimalDual$(W,d,c)$} such that
\[
    \|x^* - \tilde x\|_\infty \le (\kappa_W + 1)\|d-\tilde d\|_1 \, .
 \]
\end{lemma}
\iftoggle{focs}{}{\begin{proof}
Let $x=\tilde x-\tilde d+d$. Note that $W+x=W+d$, and also $W^\perp+s=W^\perp+c$. Thus, the systems
  \hyperref[LP-subspace-f]{\LPPrimalDual$(W,d,c)$} and   \hyperref[LP-subspace-f]{\LPPrimalDual$(W,x,s)$} define the same problem.

We apply \Cref{thm:prox-basic}  to $W,x,s$. This guarantees the existence of an optimal $(x^*,s^*)$
 to  \hyperref[LP-subspace-f]{\LPPrimalDual$(W,x,s)$}  such that $\|x^*-x\|_\infty\le \kappa_W\|x_{\Lambda(x,s)}\|_1$.
 Recall that $\Lambda(x,s)=\supp(x^-)\cup\supp(s^+)$, and thus, 
 $\|x_{\Lambda(x,s)}\|_1 = \|x^-\|_1+\|x_{\supp(s^+)}^+\|_1$.
 
 Since $\tilde x\ge 0$, we get that $\|x^-\|_1\le \|x_{\supp(x^-)}-\tilde x_{\supp(x^-)}\|_1$. Second, by the optimality of $(\tilde x, s)$, we have $\tilde x_{\supp(s^+)}=0$, and thus $x_{\supp(s^+)} = x_{\supp(s^+)}-\tilde x_{\supp(s^+)}$. Noting that $x-\tilde x=d-\tilde d$, these together imply that
 \[
 \|x^* - \tilde x\|_\infty \le \|x^* - x\|_\infty + \|x - \tilde x\|_\infty \le (\kappa_W + 1)\|d-\tilde d\|_1\, . \qedhere
 \]
\end{proof}}
We can immediately use this theorem to derive a conclusion on the support of the optimal dual solutions to  \hyperref[LP-subspace-f]{\LPPrimalDual$(W,d,c)$}, using  the optimal solution to  \hyperref[LP-subspace-f]{\LPPrimalDual$(W,\tilde d,c)$}. 
\begin{corollary}\label{cor:dual-fix-weak}
    Let $W\subseteq \R^n$ be a subspace and $c, d, \tilde d\in \R^n$.
     Let $(\tilde x, s)$ be an optimal solution to
      \hyperref[LP-subspace-f]{\LPPrimalDual$(W,\tilde d,c)$}. Let 
    \[
    R:=\{i\in [n]:\,  \tilde x_i> (\kappa_W + 1)\|d-\tilde d\|_1\}\, .
    \]
     Then for every dual optimal solution $s^*$ to \ref{LP-subspace-f}, we have $s^*_R=0$.
\end{corollary}
\iftoggle{focs}{}{\begin{proof}
    By \Cref{lem:primal_prox} there exists an optimal solution $(x', s')$ to  \hyperref[LP-subspace-f]{\LPPrimalDual$(W,d,c)$} such that $\|x' - \tilde x\|_\infty \le (\kappa_W + 1)\|d - \tilde d\|_1$.
    Consequently, $x'_R>0$, implying $s^*_R=0$ for every dual optimal $s^*$ by complementary slackness.
\end{proof}}
We now formulate a strengthening of this corollary. We show that besides setting dual variables in $R$ to 0, we are also able to set certain primal variables to 0. This will be the key to decrease the number of recursive calls from $n$ to $m$.

More precisely, we show the following. Assume $x'$ in the previous proof contains a `large' coordinate set $I_L$, significantly larger than the threshold for $R$ in \Cref{cor:dual-fix-weak}. Assume that the closure $\cl(I_L)$ contains some indices from $[n]\setminus R$. Then, we can transform $x'$ in the proof to another optimal solution $x''$ where all these coordinates are set to 0. This can be achieved by changing the coordinates in $I_L$ only, and their high value in $x''$ guarantees that they remain positive.
\begin{theorem}\label{thm:dual-fix}
Let $W\subseteq \R^n$ be a subspace and $c, d, \tilde d\in \R^n$.
     Let $(\tilde x, s)$ be an optimal solution to
      \hyperref[LP-subspace-f]{\LPPrimalDual$(W,\tilde d,c)$}, and let $\tau\ge (\kappa_W + 1)\|d-\tilde d\|_1$ and $T\ge (2n\kappa_W+1)\tau$. Let us define the following partition of $[n]$ into large, medium, and small indices.
      \iftoggle{focs}{
        \begin{align*}
          I_L&=\{i\in [n]:\,  \tilde x_i>T\}, \\ 
     I_M&=\{i\in [n]:\,  T\ge \tilde x_i> \tau\}, \\
     I_S&=\{i\in [n]:\,  \tau\ge \tilde x_i\}\, . 
        \end{align*}
      }{
     \[
     I_L=\{i\in [n]:\,  \tilde x_i>T\}\, , \quad 
     I_M=\{i\in [n]:\,  T\ge \tilde x_i> \tau\}\, , 
     \quad I_S=\{i\in [n]:\,  \tau\ge \tilde x_i\}\, . 
     \]
      }
    We further partition $I_S$ as
     \[
     I_S^0=I_S\cap \cl(I_L)\, ,\quad I_S^+=I_S\setminus \cl(I_L)\, .
     \]
     Then, there exists a primal  optimal solution $x''$  to \ref{LP-subspace-f} such that $x''_{I_L\cup I_M}>0$, and $x''_{I_S^0}=0$.
\end{theorem}
\iftoggle{focs}{}{\begin{proof}
Note that $I_L\cup I_M\subseteq R$ in \Cref{cor:dual-fix-weak}, and consider the same optimal $x'$ with $\|x'-\tilde x\|_\infty \le \tau$ as guaranteed by \Cref{lem:primal_prox}. Thus, $x'_i>2n\kappa_W\tau$ for $i\in I_L$, $x'_i>0$ for $i\in I_M$, and $2\tau\ge x'_i\ge 0$ for $i\in I_S$.

For every $i\in I_S^0$, there exists a circuit $C\in\circuits_W$ with $i\in C\subseteq I_L\cup\{i\}$. Consequently, there exists a vector $z\in W$ such that $z_{I_S^0}=x'_{I_S^0}$ and $z_{I_M\cup I_S^+}=0$ (note that we can find such a $z$ with arbitrary values on ${I_S^0}$). 
We now take the lift 
\[
w=L_{I_M\cup I_S}^W\left(z_{I_M\cup I_S}\right)\, .
\]
 Since $|x'_i|\le 2\tau$ for all $i\in I_{S^0}$, we get $\|z_{I_M\cup I_S}\|_1\le 2\tau n$. Therefore, $\|w\|_\infty\le 2n\kappa_W\tau$ by \Cref{lem:kappa-lift}.
 Consequently, for $x'':=x-w$ we have $x'' \ge 0$, $x''_{I_{S^0}}=0$ and $x''_{I_L\cup I_M}>0$. Note that  $x''$ is also an optimal solution, since $\supp(x'')\subseteq \supp(x')$. The claim follows.
\end{proof}}

\section{Constructive proximity algorithms}\label{sec:proximity-alg}
In this section, we give an algorithmic implementation of the Hoffman Proximity Theorem (\Cref{thm:master-proximity}), assuming that a feasible solution is already given: we obtain another feasible solution in strongly polynomial time that also satisfies the proximity bounds. This is given in \Cref{lem:prox-find}.

This can be derived from a  `Carath\'eodory-type' algorithm that we present in a more general form, as it may be of independent interest.
We present the algorithm in two stages, with a  basic subroutine in 
\Cref{lem:slow-lift}, and the main algorithm described in \Cref{lem:prox-find-new}. 

\begin{lemma}\label{lem:slow-lift}
Let 
$A\in \R^{m\times n}$ with  $\rk(A)=m$. Let $W=\ker(A)$, $y\in W \setminus \{0\}$,  and $J\subseteq [n]$ such that $y_J\neq 0$.
Then, in $O(m(n-m)^2 + nm^{\omega - 1})$ time, we can find a vector
\[
z=y-\sum_{j=1}^t \alpha_j y^{(j)}\, ,
\]
such that
\begin{itemize}
\item   $z\in W$  such that $z$ is sign-consistent with $y$ and  $z_J=0$ ($z=0$ is possible);
\item	$t\le n$, $\alpha_1,\ldots,\alpha_t\ge 0$, $\sum_{j=1}^t \alpha_j=1$, and
\item for $j\in [t]$, $y^{(j)}\in W$ and $y^{(j)}$ has inclusion-wise minimal support subject to $y^{(j)}_J=y_J$.
\end{itemize}
\end{lemma}
\iftoggle{focs}{}{\begin{proof}
If $y_J=0$, then we can simply output $z=y$. Otherwise, we set $\hat y=y$, and gradually modify $\hat y$. 

At the beginning iteration $\tau$ of the algorithm, we maintain $\hat y\in W$ and  $\hat \alpha\in [0,1]$, such that $\hat y$ is sign-consistent with $y,$
$\hat y=y-\sum_{j=1}^{\tau-1} \alpha_j y^{(j)}$, $\alpha_j\ge 0$, $\sum_{j=1}^{\tau-1}\alpha_j=1-\hat\alpha$, and $y^{(j)}\in W$ satisfy the required property.
We have $\hat\alpha=1$ at the beginning and $\hat \alpha=0$ at termination. At termination, we return $z=\hat y$, noting that $\hat y_J=0$. 

Further, we will maintain an index set $T\subseteq [n]\setminus J$, initialized as $T=\emptyset$, such that $\hat y_T=0$. Every iteration will add at least one new index into $T$.
If the method fails at any point, then we obtain a support minimal $y'\in W$ as in the second alternative.

Throughout, we maintain a maximal set $B\subseteq [n]\setminus (J\cup T)$ of linearly independent columns; $|B|=k\le m$. 
This can be obtained initially using Gaussian elimination; we update $B$ each time $T$ is extended. 
Thus, we transform $A$ with  row operations  to the following form.
\[
\begin{blockarray}{ccc}
B & [n]\setminus (B\cup J) & J  \\
\begin{block}{(c|c|c)}
        I_k &   \star      &   \star \\
              \cline{1-2} 
          \multicolumn{2}{c|}{$0$} & \multicolumn{1}{c}{$\star$}  \\
\end{block}
\end{blockarray} 
\]
 Note that row operations do not affect $\ker(A)$; thus, we still have $W=\ker(A)$. 

Consider now any iteration $\tau\ge 1$, with the matrix given in the above form, and the current vector $\hat y\in W$ with $\hat y_J=\hat\alpha y_J$. Let us define the vector $y'$ as 
 \[
 y'_i=\begin{cases} 
 y_i,\, &\textrm{if }i\in J\, ,\\
-A_{i,J} y_J,\, &\textrm{if }i\in B\, ,\\
0,\, &\textrm{otherwise}\, .
 \end{cases}
 \]
 Here, $A_{i,J}$ is the $i$-th row of the matrix restricted to the columns in $J$.
 First, note that the construction guarantees $y'\in W$, with $y'_J=y_J$. Further, $y'$ has minimal support subject to $y'_J=y_J$, because the columns of $A$ corresponding to $B$ are pairwise supported on disjoint rows.

Let us update $\hat y$ to
 $\hat y'=\hat y-\alpha y'$ for the smallest value $\alpha\ge 0$ such that
$\hat y'_i=0$ for some $i\in B$, or $\alpha=\hat \alpha$; we set the new value of $\hat \alpha$ to $\hat \alpha'=\hat\alpha-\alpha$. Clearly, $\hat y'\in W$ in either case, and $\hat y'_J=\hat\alpha' y_J$. 
We add $y^{(\tau)}=\hat y'$ to the combination with $\alpha_\tau=\alpha$.

If $\alpha=\hat\alpha$, then we terminate  at the end of this iteration. 
If $\alpha<\hat\alpha$, then we add all columns where $\hat y'_i=0$ to the set $T$; thus, $T$ is extended by at least one element.
We then update $B$ by pivot operations, removing all indices from $B\cap T$, and replacing them by columns in $[n]\setminus (J\cup T)$ if possible; the size of $B$ may also decrease.

Since $T$ is extended in each step, the number of iterations is bounded by $n$. The initial Gaussian elimination to find a basis can be performed in $O(nm^{\omega - 1})$ time \cite{Bunch1974,Ibarra1982}. There are at most $n-m$ additional pivot operations, each of these taking  $O(m(n-m))$ time, giving a total running time bound $O(m(n-m)^2 + nm^{\omega - 1})$.
\end{proof}}

We will use the following proximity variant later:
\begin{corollary}\label{lem:slow-lift-prox}
Let 
$A\in \R^{m\times n}$ with  $\rk(A)=m$. Let $W=\ker(A)$, $y\in W$, $M\ge 1$, and $J\subseteq [n]$.
In $O(m(n-m)^2 + nm^{\omega - 1})$
time, we can either find a vector $z\in W$,  such that $z$ is sign-consistent with $y$, $\|z-y\|_\infty\le M\|y_J\|_1$, and $z_J=0$, or 
find a lifting certificate of $M<\kappa_W$.
\end{corollary}
\iftoggle{focs}{}{\begin{proof}
If $y_J=0$, we simply output $z=y$. Otherwise, we consider the output $z\in W$ and $y^{(j)}$, $\alpha_j$, $j\in [t]$ as in \Cref{lem:slow-lift}. We have
\[
\|z-y\|_\infty=\left\|\sum_{j=1}^t \alpha_j y^{(j)}\right\|_\infty\le \sum_{j=1}^t \alpha_j \left\|y^{(j)}\right\|_\infty.
\]
If $\max_j \|y^{(j)}\|_\infty\le  M\|y_J\|_1$, then $z$ satisfies the requirements.

Assume now  that for some $y'=y^{(j)}$, we have $\|y'\|_\infty> M\|y_J\|_1$. In this case, we claim that the set $I=J\cup ([n]\setminus \supp(y'))$ and $p=y'_I$ form a lifting certificate of $M<\kappa_W$.

Clearly, $y'_I\in \pi_I(W)$, and by the support minimality of $y'$, there is a unique lift of $y'_I$ to $W$, namely, $L_I^W(p)=y'$. We have
$\|p\|_1=\|y'_J\|_1=\|y_J\|_1$. Hence,
$\|L_I^W(p)\|_\infty/\|p\|_1>M$ follows. 
\end{proof}}

The next statement formulates the outcomes of the general algorithm, using  \Cref{lem:slow-lift} as a subroutine.

\begin{lemma}\label{lem:prox-find-new}
Let $A\in \R^{m\times n}$, $\rk(A)=m$, $W=\ker(A)$,  $\ell,u\in \R^n$,  and let $x\in W$, $\ell\le x\le u$.
Then, in $O(nm(n-m)^2 + nm^{\omega - 1})$ time, we can  find vectors
$y,y^{(1)},\ldots, y^{(t)}\in W$, $t\le n$, a set $J \subseteq [n]$ and coefficients $\alpha_1,\ldots,\alpha_t\ge 0$, $\sum_{j=1}^t \alpha_j=1$ such that
\begin{enumerate}[label=(\roman*)]
\item \label{it: feas_and_eq_bound} $\ell\le y\le u$ and $y_J = (\ell^+ - u^-)_J$, and
\item  $y=\sum_{j=1}^t \alpha_j y^{(j)}$, and each $y^{(j)}\in W$ has minimal support subject to $y^{(j)}_J=(\ell^+ - u^-)_J$.
\end{enumerate}
\end{lemma}

\iftoggle{focs}{}{\begin{proof}
We repeatedly use the algorithm in \Cref{lem:slow-lift} to transform $x$ to such a vector $y$. Let 
\[
P=\{i\in [n]:\, \ell_i>0\ \ \textrm{or}\ \ u_i<0\},\ \ N=[n]\setminus P\, .
\]
If $P = \emptyset$ we return $y=0$ and $J = \emptyset$. Otherwise,
we initialize $y=x$, and 
\begin{equation}\label{eq:J-form}
J=\{i\in P:\ y_i=\ell_i>0\ \ \textrm{or}\ \ y_i=u_i<0\}\cup\{i\in N:\ y_i=0\}\, .
\end{equation}
The vector $y$ and the set $J$ will be iteratively updated, maintaining $y \in W$, $\ell\le y\le u$ and \eqref{eq:J-form}. Thus, property \ref{it: feas_and_eq_bound} is maintained throughout, since $y_J = (\ell^+ - u^-)_J$. Note that $J$ must include all indices $i\in [n]$ with $\ell_i=u_i$.

In every iteration, we call the  subroutine in \Cref{lem:slow-lift} for $y$ and $J$.
If $z=0$, then we terminate with the convex combination returned.
 Otherwise, we obtain $z\in W\setminus \{0\}$ that is sign-consistent with $y$, and $z_J=0$. 
In the very first step, $J=\emptyset$ is possible, in which case the output will be $z=y$.

We now update $y$ to the new value $y'=y-\alpha z$,
where $\alpha \ge 0$ is chosen as the smallest value where either $y'_i=\ell_i>0$, or $y'_i=u_i<0$, or $y'_i=0$ for some $i\notin J$. This $\alpha$ is finite, since there exists a coordinate $z_i\neq 0$, and $z$ is sign-consistent with $y$. 

This finishes the description of the algorithm. It is clear that the set $J$ will be strictly increased in such an iteration, and hence, the number of calls to \Cref{lem:slow-lift}  is bounded by $n$, giving the total running time bound $O(nm(n-m)^2 + n^2m^{\omega - 1})$. The second term corresponds to computing a basis of $A$, which only has to be done a single time. So the runtime reduces to $O(nm(n-m)^2 + nm^{\omega - 1})$. 
\end{proof}}

The above algorithm directly corresponds to a constructive version of \Cref{thm:master-proximity}. The proof is analogous to the proof of \Cref{lem:prox-find-new}.
\begin{corollary}\label{lem:prox-find}
Let $A\in \R^{m\times n}$, $\rk(A)=m$, $W=\ker(A)$,  $\ell,u\in \R^n$, $M\ge 1$, and let $x\in W$, $\ell\le x\le u$.
Then, in $O(nm(n-m)^2 + nm^{\omega - 1})$ time, we can either find a vector $y\in W$,  such that 
$\ell\le y\le u$, and $\|y\|_\infty\le M\|\ell^++u^-\|_1$, or 
find a lifting certificate of $M<\kappa_W$.
\end{corollary}

\section{Black-box linear programming solvers}\label{sec:blackbox}

Our feasibility and optimization algorithms in Sections~\ref{sec:feasibility} and \ref{sec:optimization}
use oracles that return approximate LP solutions. These can be implemented by using any weakly-polynomial algorithm that returns approximately optimal approximately feasible solutions as in \eqref{eq:apx-lp}. We will use the following result that summarizes recent developments on interior point methods. Whereas the papers only formulate the main statements on primal solutions, they all use primal-dual interior-point methods, and also find dual solutions with similar properties. We present the results in such a primal-dual form.
\begin{theorem}[\cite{LS19, vdb20, br2020solving}] \label{thm: best_LP_algorithms_bit_model}
Consider 
\eqref{LP_primal_dual} for $A \in \R^{m \times n}$ with $\rk(A)=m$. Assume both the primal and dual programs are feasible, let  $\delta \in [0,1]$ and $d \in \R_+^n$ be such a feasible primal solution i.e. $Ad = b$. 
Let $R_P$ and $R_D$ be the diameters of the primal and dual solution sets in $\ell_2$ norm, i.e., $\|x\|_2 \le R_P$ for all primal feasible solutions $x$, and 
$\|s\|_2\le R_D$ for all 
dual feasible solutions $(y,s)$. 
Then, we can find a vector $(x,y,s)\in \R^{n+m+n}$ with $x,s\ge 0$ such that
\begin{enumerate}[(i)]
\item
$\pr{c}{x} \le \pr{b}{y}+ \delta \cdot (\|c\|_2 \cdot R_P+\|d\|_2\cdot R_D)$,
\item $\|Ax - b\|_2 \le \delta \cdot (\|A\|_F \cdot R_P + \|b\|_2)$, and
\item $\|A^\top y+c - s\|_2 \le \delta \cdot (\|A\|_F \cdot R_D + \|c\|_2)$.
\end{enumerate}
in the following running time bounds:
\begin{enumerate}[(1)]
\item  In $O((mn + m^3) \log^{O(1)}(n) \log (n/{\delta}))$ expected runtime \cite{br2020solving}.
\item In deterministic $O(n^\omega \log^2(n) \log (n/{\delta}))$ running time, assuming $\omega \ge 13/6$
\cite{vdb20}.
\item In $O\big((\nnz(A) + m^2)\sqrt{m}\log^{O(1)} (n) \log (n/{\delta}))$ expected running time \cite{LS19}, where $\nnz(A)$ denotes the number of nonzero entries in $A$. 
\end{enumerate}
\end{theorem}
We use the notation $\CPU(A)$ to denote the `cost per unit' in these results. Namely, a $\delta$-approximate solution can be obtained in time $O(\CPU(A)\log (n/{\delta}))$, where
\begin{equation}\label{eq:CPU-def}
    \CPU(A) \le \log^{O(1)}(n) \min \{mn + m^3, n^\omega, \sqrt{m}(\nnz(A) + m^2)\}\, .
\end{equation}
We note that the third bound will not be directly applicable, since we will use the oracle in various subspaces, where the number of nonzero entries may increase.

\iftoggle{focs}{}{

We now state the main forms of feasibility and optimization oracles we use.
In Section~\ref{sec:approx}, we derive these results from \Cref{thm: best_LP_algorithms_bit_model}, by running the algorithms on an extended system.
The oracles used in Sections~\ref{sec:feasibility} and \ref{sec:optimization} can be implemented from \Cref{thm: best_LP_algorithms_real_model_feasibility} and \Cref{thm: best_LP_algorithms_real_model}. 

The main difference is in the allowable violation of primal and dual feasibility. 
 \Cref{thm: best_LP_algorithms_bit_model} returns $x,s\ge 0$, but they might not be in the right subspace: proximity bounds are given on $\|Ax - b\|_2 $ and $\|A^\top y+c - s\|_2$. In contrast, our theorems require $Ax=b$ and $A^\top y+c-s$, but allow small negative components. For a suitable represented matrix, as assumed in the statements, we can easily convert the first type of violation into the second. 

\iftoggle{focs}{
    \begin{theorem}
       \label{thm: best_LP_algorithms_real_model_feasibility}
Let 
$A\in \R^{m\times n}$, $\rk(A)=m$, $W=\ker(A)$, and  $d\in \R^n$, and assume the matrix $A_B^{-1}A$ is provided for some $m\times m$ non-singular submatrix $A_B$. Let $0< \varepsilon\le 1/2$, and let $M\ge 2$ be an estimate of $\kappa_W$. Consider a linear feasibility problem $x\in W+d, x\ge 0$.
There exists an algorithm that returns either of the following outcomes:
\begin{enumerate}[leftmargin=1.cm,label = (F\arabic*)]
\item \label{case:near-feasible} A near-feasible solution $x\in W+d$, 
$\|x^-\|_1 \le \varepsilon \|d/W\|_1$, $\|x\|_\infty\le 2M\|d\|_1$.
 \item \label{case:feasible-farkas} A Farkas certificate of  infeasibility: $s\in (W^\perp)_+$, $\pr{d}{s}<0$.
\item \label{case:certificate-feasible} A lifting certificate of $M<\kappa_W$.
\end{enumerate}
The running time is 
$O(\CPU(A)\log((M+n)/\varepsilon))$ for outcome \ref{case:near-feasible}, and 
$O(\CPU(A)\log((M+n)/\varepsilon)+(n-m)m^2 + n^\omega)$ in outcomes \ref{case:feasible-farkas} and \ref{case:certificate-feasible}.     \end{theorem}
}{
\begin{restatable}{theorem}{theoremFeasibilityProx}
\label{thm: best_LP_algorithms_real_model_feasibility}
Let 
$A\in \R^{m\times n}$, $\rk(A)=m$, $W=\ker(A)$, and  $d\in \R^n$, and assume the matrix $A_B^{-1}A$ is provided for some $m\times m$ non-singular submatrix $A_B$. Let $0< \varepsilon\le 1/2$, and let $M\ge 2$ be an estimate of $\kappa_W$. Consider a linear feasibility problem $x\in W+d, x\ge 0$.
There exists an algorithm that returns either of the following outcomes:
\begin{enumerate}[leftmargin=1.cm,label = (F\arabic*)]
\item \label{case:near-feasible} A near-feasible solution $x\in W+d$, 
$\|x^-\|_1 \le \varepsilon \|d/W\|_1$, $\|x\|_\infty\le 2M\|d\|_1$.
 \item \label{case:feasible-farkas} A Farkas certificate of  infeasibility: $s\in (W^\perp)_+$, $\pr{d}{s}<0$.
\item \label{case:certificate-feasible} A lifting certificate of $M<\kappa_W$.
\end{enumerate}
The running time is 
$O(\CPU(A)\log((M+n)/\varepsilon))$ for outcome \ref{case:near-feasible}, and 
$O(\CPU(A)\log((M+n)/\varepsilon)+(n-m)m^2 + n^\omega)$ in outcomes \ref{case:feasible-farkas} and \ref{case:certificate-feasible}. \end{restatable}
}

\iftoggle{focs}{
    \begin{theorem}
        \label{thm: best_LP_algorithms_real_model}
Let 
$A\in \R^{m\times n}$, $\rk(A)=m$, $W=\ker(A)$, and  $d\in \R^n$, and assume the matrix $A_B^{-1}A$ is provided for some $m\times m$ non-singular submatrix $A_B$. Let $0< \varepsilon\le 1/2$, and let $M\ge 2$ be an estimate of $\kappa_W$.
Consider a linear program of the form  \ref{LP-subspace-f}. 
There exists an algorithm that returns either of the following outcomes:
\begin{enumerate}[leftmargin=1.cm,label = (M\arabic*)]
\item \label{case:near-feasible-optimal} A pair of primal and dual near-feasible and near-optimal solutions $x \in W + d, s \in W^\perp +c$, that is,
\begin{equation}
    \label{eq: primal_dual_near_optimal}
    \begin{aligned}
        \|x^-\|_1 &\le \varepsilon \|d/W\|_1, \\ 
        \|s^-\|_1  &\le \varepsilon \|c/W^\perp\|_1,  \\
        \|x \circ s\|_1 &\le 5\varepsilon M \|d/W\|_1\|c/W^\perp\|_1, \\
        \|x\|_\infty & \le 2M\|d/W\|_1, \; \text{and}  \\
        \|s\|_\infty & \le 2M\|c/W^\perp\|_1 
    \end{aligned}
\end{equation}
Further, if $\kappa_W\le M$, and both primal and dual sides of \ref{LP-subspace-f} are feasible, and the optimum value is $\textrm{OPT}$, then 
\begin{equation}
\label{eq:lem_statement_on_opt}
\pr{c}{x}-5\varepsilon M \|d/W\|_1\|c/W^\perp\|_1\le \textrm{OPT}\le \pr{d}{c-s}+5\varepsilon M \|d/W\|_1\|c/W^\perp\|_1\, .
\end{equation}

\item \label{case:certificate-farkas-primal} A Farkas certificate of primal infeasibility: $s \in (W^\perp)_+$, $\pr{d}{s} < 0$.
\item \label{case:certificate-farkas-dual} A Farkas certificate of dual infeasibility:
$x \in W_+$, $\pr{c}{x} < 0$.
\item \label{case:certificate-guess-wrong} A lifting certificate of  $M<\kappa_W=\kappa_{W^\perp}$.
\end{enumerate}
The running time is 
$O(\CPU(A)\log((M+n)/\varepsilon))$ for outcome \ref{case:near-feasible-optimal}. For \ref{case:certificate-farkas-primal}, \ref{case:certificate-farkas-dual} and \ref{case:certificate-guess-wrong} the runtime is   
$O(\CPU(A)\log((M+n)/\varepsilon)+n^{2}m + n^\omega)$ if $d = 0$ or $c = 0$ and $O(\CPU(A)\log((M+n)/\varepsilon)+n^3 m)$ if $c \neq 0$ and $d \neq 0$.      \end{theorem}
}{
\begin{restatable}{theorem}{theoremOptimalityProx}
\label{thm: best_LP_algorithms_real_model}
Let 
$A\in \R^{m\times n}$, $\rk(A)=m$, $W=\ker(A)$, and  $d\in \R^n$, and assume the matrix $A_B^{-1}A$ is provided for some $m\times m$ non-singular submatrix $A_B$. Let $0< \varepsilon\le 1/2$, and let $M\ge 2$ be an estimate of $\kappa_W$.
Consider a linear program of the form  \ref{LP-subspace-f}. 
There exists an algorithm that returns either of the following outcomes:
\begin{enumerate}[leftmargin=1.cm,label = (M\arabic*)]
\item \label{case:near-feasible-optimal} A pair of primal and dual near-feasible and near-optimal solutions $x \in W + d, s \in W^\perp +c$, that is,
\begin{equation}
    \label{eq: primal_dual_near_optimal}
    \begin{aligned}
        \|x^-\|_1 &\le \varepsilon \|d/W\|_1, \\ 
        \|s^-\|_1  &\le \varepsilon \|c/W^\perp\|_1,  \\
        \|x \circ s\|_1 &\le 5\varepsilon M \|d/W\|_1\|c/W^\perp\|_1, \\
        \|x\|_\infty & \le 2M\|d/W\|_1, \; \text{and}  \\
        \|s\|_\infty & \le 2M\|c/W^\perp\|_1 
    \end{aligned}
\end{equation}
Further, if $\kappa_W\le M$, and both primal and dual sides of \ref{LP-subspace-f} are feasible, and the optimum value is $\textrm{OPT}$, then 
\begin{equation}
\label{eq:lem_statement_on_opt}
\pr{c}{x}-5\varepsilon M \|d/W\|_1\|c/W^\perp\|_1\le \textrm{OPT}\le \pr{d}{c-s}+5\varepsilon M \|d/W\|_1\|c/W^\perp\|_1\, .
\end{equation}

\item \label{case:certificate-farkas-primal} A Farkas certificate of primal infeasibility: $s \in (W^\perp)_+$, $\pr{d}{s} < 0$.
\item \label{case:certificate-farkas-dual} A Farkas certificate of dual infeasibility:
$x \in W_+$, $\pr{c}{x} < 0$.
\item \label{case:certificate-guess-wrong} A lifting certificate of  $M<\kappa_W=\kappa_{W^\perp}$.
\end{enumerate}
The running time is 
$O(\CPU(A)\log((M+n)/\varepsilon))$ for outcome \ref{case:near-feasible-optimal}. For \ref{case:certificate-farkas-primal}, \ref{case:certificate-farkas-dual} and \ref{case:certificate-guess-wrong} the runtime is   
$O(\CPU(A)\log((M+n)/\varepsilon)+n^{2}m + n^\omega)$ if $d = 0$ or $c = 0$ and $O(\CPU(A)\log((M+n)/\varepsilon)+n^3 m)$ if $c \neq 0$ and $d \neq 0$.  \end{restatable}
}
}

\section{The feasibility algorithm}\label{sec:feasibility}
Given a matrix $A\in \R^{m\times n}$ and $d\in \R^n$, we let $W=\ker(A)$. In this section, we consider the feasibility problem $x\in W+d, x\ge 0$. 

A key insight is to work with a stronger system, including a proximity constraint.
According to \Cref{cor:feas}, whenever the  problem $x\in W+d,\, x\ge 0$ is feasible and $\kappa_W\le M$, then the following system is also feasible. In fact, this would be true even with the stronger bound $M$ instead of $16M^2n$; we use this weaker bound to leave sufficient slack for the recursive argument. Note that if $d\ge 0$, then the only feasible solution is $x=d$.

\begin{tagequation}[\FeasLP$(W,d,M)$] 
\label{primal-bound}
x&\in W+d \\
\|x-d\|_\infty&\le 16M^2n\|d^-\|_1\\
x &\geq 0.
\end{tagequation}
We use a black-box approach assuming  an oracle that returns an approximately feasible solution. 
We will assume that an oracle \hyperref[oracle: prox_feas_solver]{\ProxFeasSolver$(W,d,M,\varepsilon)$} is given as in Oracle \ref{oracle: prox_feas_solver}. Outcome \ref{it:prox} gives an approximately feasible solution with a bound on the negative components and a somewhat stronger proximity guarantee as in \eqref{primal-bound}.
Outcome \ref{it:infeas} gives a Farkas certificate of infeasibility, whereas outcome \ref{it:chi} gives a lifting certificate of $M<\kappa_W$. 

\iftoggle{focs}{
\begin{oracle}
  \caption{\ProxFeasSolver$(W,d,M,\varepsilon)$}
  \label{oracle: prox_feas_solver}
  \SetKwInOut{Input}{Input}
  \SetKwInOut{Output}{Output}
  \SetKw{And}{\textbf{and}}
  
  \Input{A subspace $W\subseteq \R^n$, given as $W=\ker(A)$ for $A\in \R^{m\times n}$,
  a vector $d\in \R^n$,
 $M,\varepsilon>0$.} 
  \Output{One of the following three outcomes
  \begin{enumerate}[label=(\roman*), ref=(\roman*)]
  \item \label{it:prox} A solution $x$ to the system
  \begin{tagequation}[\ProxFeas$(W,d,M, \varepsilon)$] 
x &\in W + d \\
\|x-d\|_\infty &\leq 3M^2n \|d^-\|_1 \\
\|x^-\|_\infty &\leq \varepsilon \|d^-\|_1\ 
\label{Proximal_LP} 
\end{tagequation}
  \item \label{it:infeas} A vector $y\in W^\perp$, $y\ge 0$, $\pr{d}{y}<0$,
  \item \label{it:chi} A subset $I\subseteq [n]$ and a vector $p\in \pi_I(W)$ such that
    $\|L_I^W(p)\|_\infty>M \|p\|_1$.
  \end{enumerate}
} \end{oracle}
}{
\begin{restatable}{oracle}{oracleOne}
 \caption{\ProxFeasSolver$(W,d,M,\varepsilon)$}
  \label{oracle: prox_feas_solver}
  \SetKwInOut{Input}{Input}
  \SetKwInOut{Output}{Output}
  \SetKw{And}{\textbf{and}}
  
  \Input{A subspace $W\subseteq \R^n$, given as $W=\ker(A)$ for $A\in \R^{m\times n}$,
  a vector $d\in \R^n$,
 $M,\varepsilon>0$.} 
  \Output{One of the following three outcomes
  \begin{enumerate}[label=(\roman*), ref=(\roman*)]
  \item \label{it:prox} A solution $x$ to the system
  \begin{tagequation}[\ProxFeas$(W,d,M, \varepsilon)$] 
x &\in W + d \\
\|x-d\|_\infty &\leq 3M^2n \|d^-\|_1 \\
\|x^-\|_\infty &\leq \varepsilon \|d^-\|_1\ 
\label{Proximal_LP} 
\end{tagequation}
  \item \label{it:infeas} A vector $y\in W^\perp$, $y\ge 0$, $\pr{d}{y}<0$,
  \item \label{it:chi} A subset $I\subseteq [n]$ and a vector $p\in \pi_I(W)$ such that
    $\|L_I^W(p)\|_\infty>M \|p\|_1$.
  \end{enumerate}
} \end{restatable}
}

This oracle can be derived from \Cref{thm: best_LP_algorithms_real_model_feasibility}, by finding an approximately feasible solution to a modification of the system $x\in W+d, x\ge 0$. The derivation is given in \Cref{subsec: oracle_2_implementation}; the running time is stated as follows. Recall the definition of $\CPU(A)$ from \eqref{eq:CPU-def}.

\iftoggle{focs}{
  \begin{lemma}
    \label{thm:approx-implement}
There exists an $O(\CPU(A)\cdot \log(M+n) + nm^{\omega - 1})$ time algorithm, that either 
returns a solution to \ref{Proximal_LP}, or concludes that \ref{it:infeas} or \ref{it:chi} should be the outcome of
\hyperref[oracle: prox_feas_solver]{\ProxFeasSolver$(W,d,M,\varepsilon)$}. In the latter case, these outcomes can be obtained in additional time $O(nm^2 + n^\omega)$.    \end{lemma}
}{
\begin{restatable}{lemma}{approxImplement}
\label{thm:approx-implement}
There exists an $O(\CPU(A)\cdot \log(M+n) + nm^{\omega - 1})$ time algorithm, that either 
returns a solution to \ref{Proximal_LP}, or concludes that \ref{it:infeas} or \ref{it:chi} should be the outcome of
\hyperref[oracle: prox_feas_solver]{\ProxFeasSolver$(W,d,M,\varepsilon)$}. In the latter case, these outcomes can be obtained in additional time $O(nm^2 + n^\omega)$.  \end{restatable}
}

The next lemma will be the key technical tool of the algorithm. 
It allows to solve \ref{primal-bound} by combining an approximate
solution to \hyperref[Proximal_LP]{\ProxFeas$(W,d',M,\varepsilon)$} for some $d'\in W+d$ with an exact solution to \ref{primal-bound} obtained
recursively from a smaller system.

 Recall that for a set $K\subseteq [n]$, $\cl(K)$ denotes the closure of $K$ i.e., the unique largest set $J\subseteq [n]$ such that $K \subset J$ and $\rk(A_J)=\rk(A_K)$. 

We select a set $K$ of indices $i$ where $x_i$ is very large in the approximate solution $x$; for such indices, proximity guarantees that there must be a feasible solution $x^*\in W+d$, $x^*\ge 0$ with $x^*_i>0$. We project out all these indices, along with all other indices $J=\cl(K)\setminus K$ in their closure, and recurse on the remaining index set $I$.
 We note that the purpose of the set $J$ is to maintain property \eqref{a:no-loop} for the recursive call, i.e. to avoid creating loops from the recursive instances.

 The choice of the proximity bounds allow us to `stich together' the solution obtained on $\pi_I(W)$ from the recursive call with the approximate solution $x$ to a feasible solution to the original system. Roughly speaking, the amount of change required to cancel out all negative coordinates in $x_I$ is small enough so that $x$ remains positive on $K$.

An important feature in the scheme is  the choice of the vector $d'$ for the approximate system. This will be either $d'=d$ or $d'=d/W$; hence $W+d'=W+d$. However, this choice makes a difference due to the proximity bounds: the system \ref{primal-bound} features $\|d^-\|_1$ as well as a bound on $\|x-d\|_\infty$.

In particular, if  $\|d^-\|_1$ is  `too big', then we may end up with an empty index set $K$ and cannot recurse. In this case, we swap to $d'=d/W$; otherwise, we 
keep $d'=d$. We note that always swapping to $d'=d/W$ does not work either: \ref{primal-bound} features the bound $\|x-d\|_\infty$, and using $\|x-d/W\|_\infty$ in the approximate system may move us too far from $d$.
Fortunately, the bad cases for these two choices turn out to be complementary.
\iftoggle{focs}{}{
We note that the distinguished role of $d/W$ is due to the bound $\|x\|_2\ge \|d/W\|_2$ for any $x\in W+d$.  We formulate the following simple consequence:
\begin{lemma} \label{cl:x-d-norm}
For any subspace $W\subseteq \R^n$, and vectors $x,d\in\R^n$, $x\in W+d$, we have
\[\|x\|_\infty \ge  \frac{\|x\|_2}{\sqrt{n} }\ge \frac{\|d/W\|_1}{n}\, .\]
\end{lemma}
}
\iftoggle{focs}{}{\begin{proof}
By definition of $d/W$, we have $\|x\|_2\ge \|d/W\|_2$ for any $x\in W+d$. The claim follows combining this with the norm inequalities
$\|x\|_\infty\ge \|x\|_2/\sqrt{n}$ and $\|d/W\|_2\ge \|d/W\|_1/\sqrt{n}$.
\end{proof}}

\begin{lemma}
\label{lem:feas_lift}
Let $M,\varepsilon>0$ such that $\varepsilon \le 1/(16M^4n^4)$, let $d \in
\R^n$ and define 
\[
d'=\begin{cases}
d&\textrm{if }\|d^-\|_1<\max\left\{M\|d/W\|_1, \frac{\|d\|_\infty}{4M^2n}\right\},\\
d/W&\textrm{otherwise}.
\end{cases}
\]
Let $x$ be a feasible solution
to \hyperref[Proximal_LP]{\ProxFeas$(W,d',M,\varepsilon)$}, and let
\begin{equation}\label{eq:I}
  \begin{aligned}
K &= \set{i \in [n]\colon x_i \ge 16n^2M^3\|x^-\|_1}\, , \\
J &= \cl(K)\setminus K\, ,\\
I &= [n]\setminus \cl(K).
\end{aligned} 
\end{equation}
Then $K \neq \emptyset$. Further, if \hyperref[primal-bound]{\FeasLP$(\pi_I(W),x_I,M)$} is feasible, then 
let 
$w\in \R^I$ be any feasible solution, and
let
\[
x'=x+L_{I \cup J}^W\big((w-x_I, x_J^-)\big).
\] Then,
either $x'$ is feasible to \hyperref[primal-bound]{\FeasLP$(W,d,M)$}, or 
$\|L_{I \cup J}^W((w-x_I, x_J^-))\|_\infty >M\|(w-x_I, x_J^-)\|_1$, that is, $M<\kappa_W$.

If $y\in \pi_I(W)^\perp$, $y\ge 0$, $\pr{x_I}{y}<0$ is an infeasibility certificate to $x'\in \pi_I(W)+x_I$, $x_I\ge 0$, then  $y'=(y,0_{J\cup K})$ is an infeasibility certificate to $x\in W+d$, $x\ge 0$.
\end{lemma}
\iftoggle{focs}{}{\begin{proof}
We first show that $x'$ is feasible to \hyperref[primal-bound]{\FeasLP$(W,d,M)$}. Then, we verify that $K\neq \emptyset$.

\paragraph{Feasibility of $x'$:}
Let us first show that $x'$ is well-defined; this needs that $(w-x_I, x_J^-)\in \pi_{I \cup J}(W)$. By definition, $w \in \pi_I(W) + x_I$ means that $w=\hat w_I$ for some $\hat w\in W + x$. Clearly, $w- x_I =(\hat w-x)_I\in \pi_I(W)$. By the definition of $J$, $(0_I,z)\in \pi_{I\cup J}(W)$ for any $z\in \R^J$.
Thus, $(w-x_I,z')\in \pi_{I \cup J}(W)$ for any $z'\in \R^J$.

The containment $x'\in W+d$ is immediate, since $L_{I \cup J}^W((w-x_I, x_J^-))\in W$, and
$W+x=W+ d$. For the rest of the
proof, assume that
$\|L_{I \cup J}^W((w-x_I, x_J^-))\|_\infty\le M\|(w-x_I, x_J^-)\|_1$. 

Let us verify
$x'\ge0$. By definition, if $i\in I$, then the corresponding coordinate of
$L_{I \cup J}^W((w-x_I, x_J^-))$ equals $w_i-x_i$, and thus $x'_i=w_i\ge 0$.
Analogously, for $j \in J$, the corresponding coordinate of
$L_{I \cup J}^W((w-x_I, x_J^-))$ equals $x_j^-$ and so $x_j' = x_j + x_j^- \ge 0$. 
 For $k \in
K$, we have
$x'_k\ge x_k-\|L_{I \cup J}^W((w-x_I, x_J^-))\|_\infty$.
By definition of $K$, $x_k \ge 16n^2M^3 \|x_{I \cup J}^-\|_1$. 
 Then, $x'_k\ge 0$ follows as
\begin{equation}\label{eq:L-bound}
\|L_{I \cup J}^W\big((w-x_I, x_J^-)\big)\|_\infty\le M\|(w-x_I, x_J^-)\|_1\le nM\|(w-x_I, x_J^-)\|_\infty\le
16n^2M^3 \|x_{I \cup J}^-\|_1
\end{equation}
The first inequality is by the assumption we made; the second inequality estimates the 1-norm by the $\infty$-norm; the third inequality is since $w$ is feasible to \hyperref[primal-bound]{\FeasLP$(\pi_I(W),x_I,M)$}.

To complete the proof that $x'$ is feasible to \hyperref[primal-bound]{\FeasLP$(W,d,M)$}, it remains to verify the proximity bound $\|x'- d\|_\infty \le 16M^2n\|d^-\|_1$. First, we need an auxiliary claim.
\begin{claim} \label{cl: primal-prox}
 $\|x-d\|_\infty \le 8M^2n\|d^-\|_1$ and $\|x^-\|_1 \le n\varepsilon \|d^-\|_1$.
\end{claim}
\begin{claimproof}
If $d' = d$, then from the feasibility of $x$ to \hyperref[Proximal_LP]{\ProxFeas$(W, d', M)$}, we have $\|x- d\|_\infty\le 3M^2n\|d^-\|_1$ and $\|x^-\|_1 \le n \|x^-\|_\infty \le n\varepsilon \|d^-\|_1$. 
If $d' \neq d$, then $d' = d/W$, $\|d^-\|_1 > M\|d/W\|_1$ and $\|d^-\|_1 > \|d\|_\infty/(4M^2n)$ and so
\begin{equation*}
  \begin{aligned}
  \|x - d\|_\infty &\le \|x - d/W\|_\infty  + \|d/W\|_\infty + \|d\|_\infty \le 3M^2n \|(d/W)^-\|_1 + \|d/W\|_\infty + 4M^2n\|d^-\|_1 \\
  &\le (3M^2n + 1)\|d/W\|_1 + 4M^2n\|d^-\|_1 \le (3M^2n + 1)M^{-1} \|d^-\|_1 + 4M^2n \|d^-\|_1 \\
  &\le 8M^2n \|d^-\|_1\, , \quad \text{and} \\
  \|x^-\|_1 &\le n\|x^-\|_\infty \le n \varepsilon \|(d')^-\|_1 = n\varepsilon \|(d/W)^-\|_1 \le n\varepsilon M^{-1}\|d^-\|_1
  \end{aligned} 
\end{equation*}
proving the claim.
\end{claimproof}
Using \Cref{cl: primal-prox} as well as the bound in \eqref{eq:L-bound}, we see that 
\[
\|x'-d\|_\infty\le \|x- d\|_\infty+\|L_{I \cup J}^W\big((w-x_I, x_J^-)\big)\|_\infty\le
8M^2n\|d^-\|_1+ 16n^2M^3\|x_{I \cup J}^-\|_1\le 16M^2n\|d^-\|_1.
\]
\paragraph{Recursion on smaller subset:} We show that $K \neq \emptyset$. 
If $d = d'$ then either $\|d^-\|_1 < M \|d/W\|_1$ or $\|d^-\|_1 < \|d\|_\infty/(4M^2n)$. If $\|d^-\|_1 < M \|d/W\|_1$ then, using \Cref{cl:x-d-norm},
\[
\|x\|_\infty \ge \frac{\|d/W\|_1}{n} > \frac{\|d^-\|_1}{Mn} \ge \frac{\|x^-\|_\infty}{\varepsilon Mn} \ge \frac{\|x^-\|_1}{\varepsilon Mn^2} \ge 16n^2M^3 \|x^-\|_1 \, , 
\]
and if $\|d^-\|_1 < \|d\|_\infty/(4M^2n)$, then $\|x - d\|_\infty \le 3M^2n \|d^-\|_1$ by the call to \hyperref[Proximal_LP]{\ProxFeas$(W, d', M)$} and so 
\[
  \|x\|_\infty \ge \|d\|_\infty - \|x - d\|_\infty > 4M^2n \|d^-\|_1 - 3M^2n \|d^-\|_1 = M^2 n \|d^-\|_1 \ge \frac{M^2n}{\varepsilon n} \|x^-\|_1 \ge 16n^2M^3 \|x^-\|_1.
\]
The remaining case is $d' = d/W$. Again by \Cref{cl:x-d-norm},
\[
  \|x\|_\infty \ge \frac{\|d/W\|_1}{n} \ge \frac{\|x^-\|_\infty}{\varepsilon n} \ge \frac{\|x^-\|_1}{\varepsilon n^2} \ge 16n^2M^3 \|x^-\|_1.
\] 
\paragraph{Infeasibility certificate:} Consider now the case when we have an infeasibility cerificate
$y\in \pi_I(W)^\perp$, $y\ge 0$, $\pr{x_I}{y}<0$ to the system $x'\in \pi_I(W)+x_I$, $x_I\ge 0$.
Recall that $\pi_I(W)^\perp=(W^\perp)_I$, that is, $y'=(y,0_{J\cup K})\in W^\perp$. Clearly, $y'\ge0$, and $\pr{d}{y'}=\pr{x}{y'}=\pr{x_I}{y}<0$. In the first equality, we used that $x\in W+d$ and $y'\in W^\perp$.
\end{proof}}

\begin{algorithm}[htb!]
  \iftoggle{focs}{\SetInd{0.2em}{0.4em}}{}
  \caption{\FeasibilityAlgorithm}
  \label{alg:feasibility}
  \SetKwInOut{Input}{Input}
  \SetKwInOut{Output}{Output}
  \SetKw{And}{\textbf{and}}
  
  \Input{$A\in\R^{m\times n}$,
  $W=\ker(A) \subseteq \R^n$, $\rk(A)=m$ \iftoggle{focs}{}{satisfying \eqref{a:no-loop}}, $d\in \R^n$, $M\ge 2$.} 
  \Output{One of the following:
  {\em (i)} a solution $x$ to \ref{primal-bound}; {\em (ii)}
 a Farkas certificate of infeasibility, or {\em (iii)}
 a lifting certificate that $M < \kappa_W$.}

 \If{$\|d^-\|_1 \ge \max\left\{M\|d/W\|_1, {\|d\|_\infty}/{4M^2n}\right\}$}{
  $d \gets d / W$\;
}
  \lIf{$d \ge 0$}{
    \Return{solution $d$}
  }

  $\varepsilon \gets 1/(2Mn)^4$ \;
\Switch{$\ProxFeasSolver(W, d, M, \varepsilon)$}{
     \uCase{\ref{it:prox}}{
        $x \gets \ProxFeasSolver(W, d, M, \varepsilon)$ \;
        $K \gets \{i \in [n] : x_i \ge 16n^2M^3\|x^-\|_1\}$ \;
        $J \gets \cl(K)\setminus K$ \; 
$I \gets [n] \setminus (J \cup K)$ \;
        \Switch{output of $\FeasibilityAlgorithm(\pi_I(W), d_I, M)$}{
                 \uCase{\ref{it:prox}}{
            $x' \gets \FeasibilityAlgorithm(\pi_I(W), d_I, M)$\; \label{line: recursive_call} 
            $w \gets L_{I \cup J}^W(x' - x_I,x_J^-)$ \label{line: perform_lift}\;
            \lIf{
              $\|w\|_\infty \le M \|(x' - x_I,x_J^-)\|_1$}{
              \Return{$x + w$} \label{line: return_lift}
            }
            \lElse {
              \Return{lifting certificate of $M<\kappa_W$}.
            }
          }
   \uCase{\ref{it:infeas}}{ 
          $y\gets \FeasibilityAlgorithm(\pi_I(W), d_I, M)$\; 
          $y'\gets (y,0_{J\cup K})$\; 
            \Return{Farkas certificate $y'$ to \ref{primal-bound}}
          }
   \lCase{\ref{it:chi} }{ 
            \Return{lifting certificate of $M<\kappa_W$}
          }
  }
  }
     \lCase{\ref{it:infeas}}{ 
            \Return{Farkas certificate}
          }
   \lCase{\ref{it:chi} }{ 
            \Return{lift certificate of $M<\kappa_W$}
      }
      }
\end{algorithm}
The overall feasibility algorithm is given in \Cref{alg:feasibility}, a recursive implementation of 
\Cref{lem:feas_lift}.
 The output can be {\em(i)} a feasible solution to \ref{primal-bound}; {\em(ii)} a Farkas certificate of infeasibility, or {\em (iii)} a lifting certificate of $M<\kappa_W$. The latter will always be of the form of 
an index set $I\subseteq [n]$ and a vector $p\in
\pi_I(W)$  such that
$\|L_I^W(p)\|_\infty>M \|p\|_1$. In this case, we can restart the entire algorithm, after updating $M$ to  $\max\{\|L_I^W(p)\|_\infty/ \|p\|_1,M^2\}$.

The algorithm calls Oracle~\ref{oracle: prox_feas_solver}. For outputs \ref{it:infeas} and \ref{it:chi}, we return the Farkas certificate or the lifting certificate for $M<\kappa_W$. For  output \ref{it:prox}, we construct the sets $I$, $J$, and $K$ and recurse on $\pi_I(W)$, as in \Cref{lem:feas_lift}.

We are now ready to state the central theorem of this section, which in particular proves \Cref{thm:tardos-v2-feas}.

\begin{theorem}\label{thm:feas}
\Cref{alg:feasibility} is correct. If $M > \kappa_W$ and the system $x\in W+d, x\ge 0$ is feasible, then the algorithm returns a solution in time $O(m \CPU(A) \log(M+n)+ mn^\omega)$. 
If the system is infeasible or $M<\Le_W$, then a Farkas certificate or a lifting certificate can be obtained in additional running time $O(nm^2 + n^\omega)$.
\end{theorem}
\iftoggle{focs}{}{\begin{proof}
The bulk of of the correctness proof follows by \Cref{lem:feas_lift}. 
We prove correctness by induction on the number of variables. For $n = 1$ the statement is trivial. Assume we have already proved correctness for $<n$ variables; we now prove it for $n$ variables.

The recursive step calls the algorithm for $\pi_I(W)$; this has fewer than $n$ variables as \Cref{lem:feas_lift} shows that $K \neq \emptyset$.
 By induction hypothesis, this recursive call returns a solution to \ref{primal-bound} which can be lifted to a feasible solution in line \ref{line: return_lift} of  \Cref{alg:feasibility} by \Cref{lem:feas_lift}, unless we conclude in time $O(nm^2 + n^\omega)$ with either a Farkas certificate or $M < \kappa_W$ by \Cref{thm:approx-implement}. The second part of the same lemma asserts that  a Farkas certificate to the smaller system can be lifted to a Farkas certificate of the original system.

Also note that by \Cref{remark:lift_carry_over} we are able to transfer a lifting certificate for a subspace into one for $W$.

Let us show that the required property \eqref{a:no-loop} holds in every recursive call. Indeed, assume that some $i\in I$ forms a loop in $\pi_I(W)$. Then, there is a circuit $C\in \circuits_W$ such that $C\cap I=\{i\}$. Consequently, $\rk(K\cup J\cup\{i\})=\rk(K\cup J)$, contradicting the choice of $J\cup K=\cl(K)$.
(Note that the linear matroid on $\pi_I(W)$ corresponds to the contraction of the set $K\cup J$ in the linear matroid of $W$.)

We next show that in every recursive call, $\dim(W^\perp) > \dim(\pi_I(W)^\perp)$.
Consequently, the total number of recursive calls to \Cref{alg:feasibility}, and thus, to \hyperref[oracle: prox_feas_solver]{\ProxFeasSolver$(W,d,M,\varepsilon)$} can be bounded by $m = n - \dim(W) = \dim(W^\perp)$. 

For a contradiction, assume that $\dim(W^\perp) = \dim(\pi_I(W)^\perp)$.
Recall that $\pi_I(W)^\perp = (W^\perp)_I $. Thus, $W^\perp = \set{0}^{J \cup K} \times  (W^\perp)_I$ and therefore $W = \R^{K \cup J} \times \pi_I(W)$, which by $K \neq \emptyset$ contradicts the assumption \eqref{a:no-loop} on $W$.

The running time bounds are dominated by the running time of the oracle calls as stated in \Cref{thm:approx-implement} as well as the additional linear algebra. 
These calculations are to identify the index sets $J$, the projections $\pi_I(W)$ and computing the lifting vector as in \ref{line: return_lift}. In particular, we need to maintain throughout $W=\ker(A)$ along with the form $A^{-1}_B A$, as required in \Cref{thm: best_LP_algorithms_real_model}.

We can maintain such a representation as follows. Let $A$ denote the original matrix. Initially, we perform a Gaussian elimination to obtain a basis $B\subseteq [n]$ and the form $A'=A^{-1}_B A$ in time $O(m^2 n)$. At any point of the algorithm, let $H\subseteq [n]$ denote the current index subset; we let $W^H\subseteq \R^H$ denote the current subspace. We also maintain the basis $B\subseteq H$, and the form $A'=A^{-1}_B A_H$; thus, $W^H=\ker(A'_H)$. At every iteration, we need to partition $H=I\cup J\cup K$, and recurse to $\pi_I(W^H)$. Identifying $K$ is straightforward, but finding $J$ requires pivot operations. First, we use row operations to exchange a maximal number of columns in $K$ into $B$. Now, $\cl(K)$ is given by $K\cap B$ along with the columns generated by them in the current form; we set $J=\cl(K)\setminus K$ and let $I=H\setminus \cl(K)$.

At this point, the submatrix $A'_{(I \cap B) \times I}$ represents $\pi_I(W^H)$. Here, the rows $I \cap B$ correspond in abuse of notation to the columns of $B$ that are in $I$ of the identity submatrix of $A'$. 
Hence, we can charge to every removed column at most two basis exchanges, giving a total of $O(n^2 m)$ for these operations. 

It is left to compute the lift from line \ref{line: perform_lift}, which recursively has to be performed $m$ times and takes $O(\min\{m^2n, n^\omega\})$ by \Cref{lem:compute-lift}. 

\end{proof}}

\section{The optimization algorithm}\label{sec:optimization}
In this section, we show how \ref{LP-subspace-f} can be solved using an approximate LP solver.
As in the feasibility algorithm, we let $M$ denote our current upper estimate on $2\kappa_W$.
We  present an algorithm that comprises an Inner and an Outer Loop. The calls to the approximate LP solver will happen inside the Inner Loop.

The outer loop gives an algorithmic implementation of \Cref{thm:dual-fix}. 
The subroutine \InnerLoop$(W,d,c,M)$ returns a solution $(\tilde d,\tilde x,s)$, where $\tilde d$ is a `perturbed' version of $d$, and $(\tilde x, s)$ are optimal solutions to \hyperref[LP-subspace-f]{\LPPrimalDual$(W,\tilde d,c)$}.
We get $(\tilde d,\tilde x,\tilde s)$  as solutions to the following system. 
\begin{tagequation}[\InnerSysF$(W,d,c,M)$]
\label{eq:Inner}
\|\tilde d-d\|_1&\le \frac{\|\tilde x\|_\infty}{4n^2M^2}\\
\tilde x&\in W+\tilde d \\
\tilde s&\in W^\perp + c\\
\pr{\tilde x}{\tilde s}&=0\\
\tilde x,\tilde s&\geq 0.
\end{tagequation}
The subroutine will be described in Section~\ref{sec:inner}; we now state the running time.
\begin{theorem}\label{thm:inner}
Assume we are given a matrix $A\in\R^{m\times n}$, vectors $c\in\R^n$ and $d\in \R^n_+$; let $W=\ker(A)$, and $M$ be an estimate on $\kappa_W$.
There exists an $O(n\CPU(A) \log(M + n) + n^{\omega + 1})$ time algorithm (\Cref{alg:inner_loop}) that returns a solution $(\tilde d,\tilde x,\tilde s)$ to \ref{eq:Inner}, or decides that the system is either primally infeasible or $M < \kappa_W$. To obtain a Farkas certificate of primal infeasibility or a certificate that $M< \kappa_W$ is obtained in additional time $O(n^3m)$.
\end{theorem}
 
The overall algorithm described in Section~\ref{sec:overall-1} repeatedly calls \InnerLoop{} to set primal and dual variables to 0 according to \Cref{thm:dual-fix}, and recurses to lower dimensional subspaces. The final optimal solutions are obtained via calling the feasibility algorithm on both primal and dual side. 
The drawback of this variant is that, in case the algorithm fails, we do not obtain lifting certificates of $M<\kappa_W$---a guarantee we can achieve for feasibility in Section~\ref{sec:feasibility}.
A post-processing, as described in \Cref{lem: certificate_for_optimization} is able to create a certificate of $M < \Le_W$, at an additional computational expense.

\subsection{The Outer Loop}\label{sec:overall-1}
Consider an instance of \ref{LP-subspace-f} and an estimate $M$ on $\kappa_W$.
We first use the feasibility algorithm and check if both systems $x\in W+d$, $x\ge0$ and $s\in W^\perp+c$, $s\ge 0$ are feasible. For the remainder of this section, let us assume both these systems are feasible, and consequently, \ref{LP-subspace-f} is also feasible. Moreover, 
 we can write an equivalent system with nonnegative $d\ge 0$. (We could also impose $c\ge 0$, but this will not be used).

 The overall algorithm is presented in Algorithm~\ref{alg:simple_tardos}. We let $(W^0,d^0,c^0)$ denote the original input, where $W^0$ is a subspace of $\R^n$, and $d^0,c^0\in \R^n$, $d^0\ge 0$. We will maintain an index set $I\subseteq [n]$, initialized as $I=[n]$. 
 We gradually move every index into the set $B$ or $N$. 
We apply \Cref{thm:dual-fix} with thresholds $\tau=\|x\|_\infty/(3n^2 M)$ and $T=\|x\|_\infty/n$. The bound $\|\tilde d-d\|_1\le {\|\tilde x\|_\infty}/{4n^2M^2}$ in \ref{eq:Inner} guarantees that these are suitable choices.
 Assuming that $M\ge \kappa_W$, we obtain $s^*_B=0$ for every dual optimal solution $s^*$, and that there exists a primal optimal solution $x^*$ with $x^*_N=0$. Note the asymmetry between the two sides; the weaker guarantee on $N$ already suffices for correctness.

 At every iteration, we have an index set $I\subseteq [n]$ of `undecided indices' and a subspace $W\subseteq \R^I$. We consider the partition $I=I_L\cup I_M\cup I_S^+\cup I_S^0$ according to  \Cref{thm:dual-fix}, and add $I_L\cup I_M$ to $B$ and $I_S^0$ to $N$. The optimal solution $x''$ guaranteed in the theorem implies that the optimum value is the same on $W$ and $W'=W\cap \R^I_{I_L\cup I_S\cup I_S^+}$, i.e. the subspace where all entries in $I_S^0$ are forced to 0. We then update $I=I_S^+$ as the remaining set of undecided indices and recurse on the subspace $\pi_I(W')$.

 The algorithm terminates when $d$ is contained in $W$. The remaining indices $I$ are split up between $B$ and $N$ based on whether they are in the support of the optimal dual solution of the perturbed system. Finally, we obtain the primal and dual solutions by solving feasibility problems on the subsets $B$ and $N$. If both are feasible, they form a complementary pair of primal and dual solutions, and hence they are optimal. In case of a failure, we conclude that the underlying assumption $M\ge \kappa_W$ was wrong.

 \begin{algorithm}[htb]
  \iftoggle{focs}{\SetInd{0.2em}{0.4em}}{}
  \caption{Optimization Algorithm}
\label{alg:simple_tardos}
  \SetKwInOut{Input}{Input}
  \SetKwInOut{Output}{Output}
  \SetKw{And}{\textbf{and}}
  
  \Input{$W^0\subseteq \R^n$\, ; $c^0\in \R^n$, $d^0\in \R_+^n$ such that \LPPrimalDual$(W^0,d^0,c^0)$
  is feasible, and $M\ge 2$.
  }
  \Output{Solution $(x^*,s^*)$ to \LPPrimalDual$(W^0,d^0,c^0)$ or the answer $M<\kappa_{W^0}$.}

  $W\gets W^0$\, ;\, $d\gets d^0$\, ;\, $c\gets c^0/W^\perp$ \;
  $I\gets [n]$\, ; 
  $B \gets \emptyset$; $N \gets \emptyset$ \;
  \While(\label{line: while_loop}){$I\neq \emptyset$ and $d \notin W$}{
    \lIf{\InnerLoop$(W,d,c,M)$ returns  $M < \kappa_W$}{
      \Return $M < \kappa_W$.
    }
    $(\tilde d, \tilde x, \tilde s)\gets$ \InnerLoop$(W,d,c,M)$ \;
    $I_L\gets \{i\in I:\,  \tilde x_i>\|\tilde x\|_\infty/n\}$ \;
     $I_M\gets \{i\in I:\,  \|\tilde x\|_\infty/n\ge \tilde x_i> \|\tilde x\|_\infty/(3n^2 M)\}$ \; 
    $I_S\gets\{i\in I:\,  \|\tilde x\|_\infty/(3n^2 M)\ge \tilde x_i\}$ \;
    $I_S^0\gets I_S\cap \cl(I_L)$\, ; $I_S^+\gets I_S\setminus \cl(I_L)$ \; 
   $B \gets B \cup I_L\cup I_M$ 
   \iftoggle{focs}{;}{\;}
    $N \gets N \cup I_S^0$ \;
    $W'\gets W\cap \R^I_{I_L\cup I_M\cup I_S^+}$ \;
    $I \gets I_S^+$ \label{line: I_update}\;
    $W\gets \pi_{I}(W')$\, ;\, $d\gets d_I$\, ;\, $c\gets \tilde s_I$ ;
  }
  $N \gets N \cup \big(I \cap \supp(\tilde s)\big)$ \label{line: final_N_update}
  \iftoggle{focs}{;}{\;}
  $B \gets B \cup \big(I \setminus \supp(\tilde s)\big)$ \label{line: final_B_upate}\;
  \If{\FeasibilityAlgorithm$(W^0_B,d^0)$ and \FeasibilityAlgorithm$((W^0)^\perp_N,c^0)$ are feasible}{
    $x^*\gets \FeasibilityAlgorithm(W^0_B,d^0)$ \label{alg_tardos: primal_feas}\; 
    $s^* \gets \FeasibilityAlgorithm((W^0)^\perp_N,c^0)$ \label{alg_tardos: dual_feas}\;
    \Return{$(x^*, s^*)$}
    } 
  \lElse{
    \Return{$M<\Le_{W^0}$}
  }
\end{algorithm}

\begin{theorem}\label{thm:optimization-main}
Assuming that \hyperref[LP-subspace-f]{\LPPrimalDual$(W^0,d^0,c^0)$} is both primal and dual feasible,
Algorithm~\ref{alg:simple_tardos} either finds an optimal solution to  \hyperref[LP-subspace-f]{\LPPrimalDual$(W^0,d^0,c^0)$} or correctly concludes that $M<\kappa_{W^0}$ by at most $m$ calls to the subroutine \hyperref[sec:inner]{\InnerLoop}. The total runtime is  $O(mn\CPU(A)\log(M + n)+ mn^{\omega +1})$to find a primal-dual optimal pair. Obtaining a lifting certificate requires additional time $O(n^3m^2)$.
\end{theorem}
\iftoggle{focs}{}{\begin{proof}
Let us assume that $M\ge \Le_{W^0}$. The correctness of constructing the sets $B$ and $N$ follows from \Cref{thm:dual-fix}. Clearly, $I_L\neq \emptyset$.

The bound of at most $m$ recursive calls can be established with the same argumentation as in the proof of \Cref{thm:feas} for the feasibility algorithm, by showing that $\dim(W^\perp)$ decreases by at least one in every step.
Assume $\dim(\pi_{I_S^+}(W')^\perp) = \dim(W^\perp)$. Note that $\pi_{I_S^+}(W')^\perp = \pi_{I_S^+}(W)^\perp = (W^\perp)_{I_S^+}$ and and so $\dim(W^\perp) = \dim((W^\perp)_{I_S^+})$. But then $W = \R^{I\setminus I_S^+} \times \pi_{I_S^+}(W)$, a contradiction as $I \setminus I_S^+ \neq \emptyset$ and assumption \eqref{a:no-loop} on $W$ holds. Also, assumption \eqref{a:no-loop} holds for $\pi_{I_S^+}(W')$ as $\cl(K) \cap I_S^+ = \emptyset$.

The running time bound then follows from \Cref{thm:inner}, noting that the computational cost is dominated by the calls to \InnerLoop{} and recomputing the representation of the subspace $W$. This can be done similarly as in the proof of \Cref{thm:feas}. The claim on the lifting certificates will be proved in \Cref{lem: certificate_for_optimization}.
\end{proof}}

By proving \Cref{thm:optimization-main}, we also proved the special case \Cref{thm:tardos-v2-opt}.

\subsection{The Inner Loop}\label{sec:inner}

For the Inner Loop and $d \ge 0$ we formulate the stronger version \ref{eq:Inner_Prox} of \hyperref[eq:Inner]{$\InnerSysF(W,d,c,M)$}, which maintains dual proximity and therefore---in similar vein as the feasibility algorithm---only needs an oracle with precision $(Mn)^{-O(1)}$ for recursive calls.

\begin{tagequation}[\textsl{F-Primal-Prox}$(W,d,c,M,\varepsilon)$]
  \label{eq:Inner_Prox}
x&\in W+d \\
  s&\in W^\perp + c\\
  \|s - c\|_\infty &\le 16M^2 n\|c_{\Lambda(c,d)}\|_1 \\
  \|x_{\Lambda(x,s)}\|_\infty &\le 2\varepsilon n \|x\|_\infty\\
s&\ge 0
\end{tagequation}
Given an output $(x,s)$ to this system with $\varepsilon\le 1/(32 M^4n^4)$, we obtain a solution to  \ref{eq:Inner} as
\begin{equation}\label{eq:def-tilde}
\tilde x_i=
\begin{cases} 0 & \textrm{if }i\in \Lambda(x,s)\, ,\\
x_i & \textrm{otherwise},
\end{cases}\, , \quad
\tilde d=d-x+\tilde x\,  ,\quad \tilde s=s\, .
\end{equation}
We will assume that the following oracle (Oracle \ref{oracle: prox_opt_solver}) is available, that returns a solution $(\tilde c,x, s)$ to the system \ref{eq:ProxOpt}, a primal or dual infeasibility certificate, or a lifting certificate.
For the input, we require the nonnegativity $c \ge 0$. This will be perturbed to $\tilde c$,  and $(x,s)$ will be near-optimal and near-feasible primal and dual solutions with respect to the perturbed system, satisfying a primal proximity constraint.

\iftoggle{focs}{
  \begin{oracle}
  \caption{\ProxOpt$(W,d,c,M,\varepsilon)$}
\label{oracle: prox_opt_solver}
\SetKwInOut{Input}{Input}
\SetKwInOut{Output}{Output}
\Input{$W \subseteq \R^n$, $c\in \R^n_+$, $d \in \R^n$, $M \ge 2, \varepsilon > 0$.}
\Output{One of the following:
\begin{enumerate}[label=(\roman*)]
\item \label{it:proxopt_succ} A solution to the system
\begin{tagequation}[\textsl{Prox-Opt}$(W,d,c,M,\varepsilon)$]
  \label{eq:ProxOpt}
  x &\in W +d\\
   s &\in W^\perp + \tilde c \\
\|x_{\Lambda(x,s)}\|_\infty & \leq       \varepsilon \|d_{\Lambda(d,c)}\|_1\\
\|x-d\|_\infty&\le 3M^2n\|d_{\Lambda(d,c)}\|_1 \\
       \|c-\tilde c\|_\infty &\le \frac{\varepsilon}{n} \|c/W^\perp \|_1\\
       c-\tilde c &\ge 0 \\
    s & \ge 0 \\
\end{tagequation}
\item \label{it:infeas_primal} A vector $y\in W^\perp$, $y\ge 0$, $\pr{d}{y}<0$,
\item \label{it:infeas_dual} A vector $x\in W$, $x\ge 0$, $\pr{c}{x}<0$,
\item \label{it:oracle_2_cert} A lifting certificate of $M<\kappa_W=\kappa_{W^\perp}$.
\end{enumerate}
}   \end{oracle}
}{
\begin{restatable}{oracle}{oracleTwo}
\caption{\ProxOpt$(W,d,c,M,\varepsilon)$}
\label{oracle: prox_opt_solver}
\SetKwInOut{Input}{Input}
\SetKwInOut{Output}{Output}
\Input{$W \subseteq \R^n$, $c\in \R^n_+$, $d \in \R^n$, $M \ge 2, \varepsilon > 0$.}
\Output{One of the following:
\begin{enumerate}[label=(\roman*)]
\item \label{it:proxopt_succ} A solution to the system
\begin{tagequation}[\textsl{Prox-Opt}$(W,d,c,M,\varepsilon)$]
  \label{eq:ProxOpt}
  x &\in W +d\\
   s &\in W^\perp + \tilde c \\
\|x_{\Lambda(x,s)}\|_\infty & \leq       \varepsilon \|d_{\Lambda(d,c)}\|_1\\
\|x-d\|_\infty&\le 3M^2n\|d_{\Lambda(d,c)}\|_1 \\
       \|c-\tilde c\|_\infty &\le \frac{\varepsilon}{n} \|c/W^\perp \|_1\\
       c-\tilde c &\ge 0 \\
    s & \ge 0 \\
\end{tagequation}
\item \label{it:infeas_primal} A vector $y\in W^\perp$, $y\ge 0$, $\pr{d}{y}<0$,
\item \label{it:infeas_dual} A vector $x\in W$, $x\ge 0$, $\pr{c}{x}<0$,
\item \label{it:oracle_2_cert} A lifting certificate of $M<\kappa_W=\kappa_{W^\perp}$.
\end{enumerate}
} \end{restatable}
}

\iftoggle{focs}{}{
The implementation of Oracle \ref{oracle: prox_opt_solver} is given in \Cref{subsec: oracle_2_implementation}, using \Cref{thm: best_LP_algorithms_real_model}.
}

\iftoggle{focs}{
\begin{lemma}
  \label{lem:oracle_2_implementation}
  Assume we are given a matrix $A\in\R^{m\times n}$, vectors $c\in\R^n_+$ and $d\in \R^n$; let $W=\ker(A)$, and $M$ be an estimate on $\kappa_W$. Further, let $0<\varepsilon<1$. There exists an $O(\CPU(A) \cdot \log (M + n) + nm^{\omega - 1})$ time algorithm, that either returns a solution to \ref{eq:ProxOpt} or concludes that \ref{it:infeas_primal}, \ref{it:infeas_dual} or \ref{it:oracle_2_cert} should be the outcome of \hyperref[oracle: prox_opt_solver]{\ProxOpt$(W,d,c,M,\varepsilon)$}. These latter outcomes require an additional computational time $O(n^3m)$.  \end{lemma}
}{
\begin{restatable}{lemma}{oracleTwoImplementation} \label{lem:oracle_2_implementation}
  Assume we are given a matrix $A\in\R^{m\times n}$, vectors $c\in\R^n_+$ and $d\in \R^n$; let $W=\ker(A)$, and $M$ be an estimate on $\kappa_W$. Further, let $0<\varepsilon<1$. There exists an $O(\CPU(A) \cdot \log (M + n) + nm^{\omega - 1})$ time algorithm, that either returns a solution to \ref{eq:ProxOpt} or concludes that \ref{it:infeas_primal}, \ref{it:infeas_dual} or \ref{it:oracle_2_cert} should be the outcome of \hyperref[oracle: prox_opt_solver]{\ProxOpt$(W,d,c,M,\varepsilon)$}. These latter outcomes require an additional computational time $O(n^3m)$.  \end{restatable}
}

\InnerLoop$(W,d,c,M,\varepsilon)$ recursively calls itself and Oracle \ref{oracle: prox_opt_solver}, while maintaining dual proximity. The Oracle will be called for the dual system. 

\iftoggle{focs}{}{
 For convenience we now state it with primal and dual side flipped:

 \begin{tagequation}[\textsl{Prox-Opt}$(W^\perp,c,d,M,\varepsilon)$]
  \label{eq:ProxOptDual}
  s &\in W^\perp +c\\
   x &\in W + \tilde d \\
\|s_{\Lambda(s,x)}\|_\infty & \leq \varepsilon\|c_{\Lambda(c,d)}\|_1 \\
    \|s-c\|_\infty&\le 3M^2n\|c_{\Lambda(c,d)}\|_1 \\
      \|d-\tilde d\|_\infty &\le \frac{\varepsilon}{n} \|d/W\|_1\\
      d-\tilde d &\ge 0 \\
    x & \ge 0 \\
\end{tagequation}
}

\begin{lemma}
 \label{lem: opt_recursive_call} 
 Let $M$ be an estimate on $\kappa_W$, $0<\varepsilon<1/(32M^4 n^4)$. Let us define
\[
c'=\begin{cases}
c&\textrm{if }\|c_{\Lambda(c,d)}\|_1<\max\left\{M\|c/W^\perp\|_1, \frac{\|c\|_\infty}{4M^2n}\right\},\\
c/W^\perp&\textrm{otherwise}.
\end{cases}
\]
Let $(x,s)$ be a feasible solution to \hyperref[eq:ProxOpt]{\textsl{Prox-Opt}$(W^\perp, c', d, M, \varepsilon)$} and let 
 \begin{equation*}
  I = \{i \in [n]: s_i \le 16n^3M^3\|s_{\Lambda(s,x)}\|_1\}\, , \quad J = [n]\setminus I\, .
 \end{equation*}
Then, the following hold:
 \begin{enumerate}[label=(\roman*)]
 \item If
 $I=\emptyset$, then we must have $s\ge 0$ and $d\in W$, and $(0,s)$ is feasible to \ref{eq:Inner_Prox}.
 \item 
 If $I\neq \emptyset$, then  let $(w,z)$ be a solution to \hyperref[eq:Inner_Prox]{\textsl{F-Primal-Prox$(W_I,x_I,s_I,M,\varepsilon)$}}, and define
 \begin{equation*}
\tilde x = (w, x_J)+d-\tilde d,   \quad \textrm{and}\quad 
    \tilde s = s + L_I^{W^\perp}(z - s_I)\, .
 \end{equation*}
 Then  either $(\tilde x, \tilde s)$ is feasible to \ref{eq:Inner_Prox} or 
we obtain a lifting certificate of
$M < \kappa_{W^\perp} = \kappa_W$.
  \item $J\neq \emptyset$.
\end{enumerate}
\end{lemma}
\iftoggle{focs}{}{\begin{proof}
Let us start with two  simple claims.
\begin{claim}\label{cl:x-J-bound}
 $x_J= 0$. 
 \end{claim}
 \begin{claimproof}
The vector $x$ is required to be nonnegative. Assume that for some $j\in J$, $x_j>0$, and thus $j\in \Lambda(s,x)$. For such an index, $s_j > 16n^3M^3\|s_{\Lambda(s,x)}\|_1\ge 16n^3M^3 s_j$, a clear contradiction.
\end{claimproof}

\begin{claim}\label{cl:s-c-bound}
 $\|s-c\|_\infty\le  8M^2n\|c_{\Lambda(c,d)}\|_1$.
 \end{claim}
 \begin{claimproof}
If $c'=c$, then the stronger bound with coefficient $3M^2n$ is included in \hyperref[eq:ProxOpt]{\textsl{Prox-Opt}$(W^\perp, c', d, M, \varepsilon)$}.
Assume $c'=c/W^\perp$. Then,
\[
\begin{aligned}
\|s-c\|_\infty&\le \|s-c/W^\perp \|_\infty+\|c/W^\perp \|_\infty+\|c\|_\infty\\
&\le 3M^2n\|c/W^\perp\|_1+\|c/W^\perp\|_1+ \|c\|_\infty\\
&\le \frac{3M^2n+1}{M} \|c_{\Lambda(c,d)}\|_1+4M^2n\|c_{\Lambda(c,d)}\|_1\\
&\le 8M^2n\|c_{\Lambda(c,d)}\|_1\, .
\end{aligned}
\]
The second inequality used the bound in $\|s-c'\|_\infty$ in the proximal solver, upper bounding the left hand side as $3M^2n\|c'\|_1$. The third inequality uses the bounds on $\|c/W^\perp\|$ and $\|c\|_\infty$ in the case  when $c'=c/W^\perp$
is chosen.
\end{claimproof}

\paragraph{Part (i)} Assume $I=\emptyset$ and $J=[n]$. Then, for every $j\in [n]$, $s_j>0$ by definition, and $x_j=0$ by \Cref{cl:x-J-bound}.  Since $x=0\in W+\tilde d$, we have $\tilde d\in W$, and therefore $d-\tilde d\in W+d$. By \Cref{cl:x-d-norm}, $\|d-\tilde d\|_\infty\ge \|d/W\|_1/n$. On the other hand, we have the upper bound $\|d-\tilde d\|_\infty\le \varepsilon \|d/W\|_1/n$. This is only possible if $d/W=0$, that is, $d\in W$. Using also \Cref{cl:s-c-bound}, we conclude that  $(0,s)$ is feasible to \ref{eq:Inner_Prox}.

\paragraph{Part (ii)}
The vector $\tilde s$ is well-defined, since $z-s_I\in (W_I)^\perp=\pi_I(W^\perp)$.
  We have $\tilde s \in W^\perp + s = W^\perp +c'=W^\perp+c$ by definition. 
  Further,  $(w,0_J)\in W+(x_I,0_J)$, and therefore, $(w,x_J)\in W+x=W+\tilde d$.
  Then, $(w,x_J)-\tilde d+d\in W+d$.

Let us assume that 
 $\|L_I^{W^\perp}(z - s_I)\|_\infty \le M\|z - s_I\|_1$; otherwise, we can return a lifting certificate of $M<\kappa_{W^\perp}$. Thus,
\begin{equation}\label{eq:lift-I}
  \|L_I^{W^\perp}(z - s_I)\|_\infty \le M \|z- s_I\|_1 \le 16M^3 n^2\|s_{\Lambda(s_I,x_I)}\|_1 \le 16M^3 n^3\|s_{\Lambda(s,x)}\|_\infty\, ,
\end{equation}
where the second inequality comes from the bound on $\|z-s_I\|_1$ in  \hyperref[eq:Inner_Prox]{\textsl{F-Primal-Prox$(W_I,x'_I,s_I,M,\varepsilon)$}}.
We are ready to show $\tilde s \ge 0$. For $j\in I$, $\tilde s_j= z_j \ge 0$, and for $j \in J$, we have  
\begin{equation}
  \tilde s_j \ge s_j - \|L_I^{W^\perp}(z-s_I)\|_\infty \ge 16M^3n^3 \|s_{\Lambda(s,x)}\|_1 -16M^3n^3\|s_{\Lambda(s,x)}\|_\infty \ge 0\, .
\end{equation}
We now turn to dual proximity:
\begin{align*}
\|\tilde s - c\|_\infty & \le \|\tilde s - s\|_\infty + \|s - c\|_\infty \le \|L_I^{W^\perp}(z - s_I)\|_\infty + 8 M^2n \|c_{\Lambda(c,d)}\|_1 \\
&\le 16M^3n^3 \|s_{\Lambda(s,x)}\|_\infty + 8 M^2n  \|c_{\Lambda(c,d)}\|_1 \\
& \le 16  M^3n^3 \varepsilon \|c_{\Lambda(c,d)}\|_1 + 8M^2n  \|c_{\Lambda(c,d)}\|_1 \le 16M^2n \|c_{\Lambda(c,d)}\|_1\, ,
\end{align*}
using \eqref{eq:lift-I}, \Cref{cl:s-c-bound} and that $\varepsilon<1/(2Mn^2)$.

\medskip

It is left to show $\|\tilde x_{\Lambda(\tilde x,\tilde s)}\|_\infty\le 2\varepsilon n\|\tilde x\|_\infty$. Recall that $\tilde x=x+d-\tilde d$, and that $d-\tilde d\ge0$. Further,
 $\|d-\tilde d\|_\infty\le \varepsilon\|d/W\|_1/n\le \varepsilon\|\tilde x\|_\infty$ using \Cref{cl:x-d-norm}.

Let us fix an index $j\in [n]$ such that $\tilde x_j<0$ or $\tilde s_j>0$; our goal is to show $|\tilde x_j|\le 2\varepsilon n\|\tilde x\|_\infty$.

Assume first $j\in J$. Then, $\tilde x_j=x_j+d_j-\tilde d_j=d_j-\tilde d_j$ by \Cref{cl:x-J-bound}. Hence,  we obtain the stronger bound $|\tilde x_j|\le \|d-\tilde d\|_\infty\le \varepsilon\|\tilde x\|_\infty$.
For the rest, let us assume $j\in I$. Then, $\tilde x_j=w_j+d_j-\tilde d_j$. 
\begin{claim}
 $|w_j|\le  2\varepsilon(n-1) \|w\|_\infty$.
 \end{claim}
 \begin{claimproof}
 From the recursive call, we know that $\|w_{\Lambda(w,z)}\|_\infty\le 2\varepsilon(n-1) \|w\|_\infty$, using that $|I|\le n-1$. The claim follows by showing that $j\in \Lambda(w,z)$.
  Indeed, $\tilde s_I=z$, and thus $\tilde s_j>0$ means $j\in \supp(z^+)$. If $\tilde x_j<0$, then  $0>\tilde x_j=w_j+d_j-\tilde d_j\ge w_j$ by the nonnegativity of $d-\tilde d$. Hence, $j\in \supp(w^-)$.
\end{claimproof}

  Since $\tilde x_I-w=d_I-\tilde d_I$, we obtain the bound $\|\tilde x_I-w\|_\infty\le \varepsilon\|\tilde x\|_\infty$, implying $\|w\|_\infty\le (1+\varepsilon)\|\tilde x\|_\infty$. We can thus bound
  \[
\begin{aligned}
|\tilde x_i|&\le |w_i|+\varepsilon \|\tilde x\|_\infty\le 2\varepsilon(n-1)\|w\|_\infty+ \varepsilon \|\tilde x\|_\infty\\
&\le \varepsilon\left((2n-2)(1+\varepsilon)+1\right)\|\tilde x\|_\infty\le 2n\varepsilon\|\tilde x\|_\infty\, ,
\end{aligned}
\]
using that $\varepsilon<1/(2n)$.

\paragraph{Part (iii)}
 By the definition of $I$, we have
\begin{equation}\label{eq:s-I}
\|s_I\|_\infty < 16M^3n^2\|s_{\Lambda(s,x)}\|_1 \le 16M^3n^3 \|s_{\Lambda(s,x)}\|_\infty\le 16\varepsilon M^3n^3 \|c_{\Lambda(c,d)}\|_1\, .
\end{equation}
First, let us assume $c'=c$. This can happen in two cases: if $\|c_{\Lambda(c,d)}\|_1<M\|c/W^\perp\|_1$ or if $\|c_{\Lambda(c,d)}\|_1<
{\|c\|_\infty}/{(4M^2n)}$. In the first case, 
\[
\|s_I\|_\infty< 16\varepsilon M^3 n^3 \|c_{\Lambda(c,d)}\|_1\
\le \frac{\|c_{\Lambda(c,d)}\|_1}{2Mn} \le \frac{\|c/W^\perp\|_1}{2n}\le \frac{\|s\|_2}{2\sqrt{n}}\, , 
\]
showing that $I\neq [n]$. The second inequality uses the choice of $\varepsilon < 1/(32M^4n^4)$, and the last inequality uses the minimum-norm property \Cref{cl:x-d-norm}. In the second case, we can bound
\begin{equation}
\begin{aligned}
\|s\|_\infty &\ge \|c\|_\infty-\|c-s\|_\infty\ge 4M^2n \|c_{\Lambda(c,d)}\|_1- 3M^2n\|c_{\Lambda(c,d)}\|_1 = M^2n \|c_{\Lambda(c,d)}\|_1 \\
& \ge 16M^6n^5\|s_{\Lambda(s,x)}\|_\infty > \|s_I\|_\infty \, ,
\end{aligned}
\end{equation}
again implying $I\neq [n]$.

Let us now turn to the case $c'=c/W^\perp$. From \eqref{eq:s-I}, we get
\[
\|s_I\|_\infty< 16\varepsilon M^3 n^3 \|c'_{\Lambda(c',d)}\|_1\le 
16\varepsilon M^3 n^3 \varepsilon \|c/W^\perp\|_1\le \frac{\|c/W^\perp \|_2}{2M\sqrt{n}}\le \frac{\|s\|_2}{2M\sqrt{n}}\, ,
\]
showing again $I\neq [n]$.
\end{proof}}

Algorithm~\ref{alg:inner_loop} implements the recursive calls in accordance with \Cref{lem: opt_recursive_call}.
This is the algorithm asserted in \Cref{thm:inner}.

\begin{algorithm}[htb]
  \iftoggle{focs}{\SetInd{0.2em}{0.4em}}{}
  \caption{\textsc{InnerLoop}}
  \label{alg:inner_loop}
  \SetKwInOut{Input}{Input}
  \SetKwInOut{Output}{Output}
  \SetKw{And}{\textbf{and}}
  
  \Input{$W \subseteq \R^n$, $c \in \R^n$, $d\in \R_+^n$  such that \LPPrimalDual$(W,d,c)$ is feasible, $M \ge 1$.} 
  \Output{Either a solution $(\tilde d, \tilde x, \tilde s)$ to \ref{eq:Inner_Prox} or lifting certificate of $M < \kappa_W$.}

  $\varepsilon \gets 1/{(32M^4 n^4)}$ \;
  \lIf{$\|c_{\Lambda(c,d)}\|_1 \ge \max\left\{M\|c/W^\perp\|_1, \frac{\|c\|_\infty}{4M^2n}\right\}$}{
    $c \gets c / W^\perp$
  }
  \If{\hyperref[oracle: prox_opt_solver]{\ProxOpt$(W^\perp,c,d,M,\varepsilon)$} in \ref{it:proxopt_succ}}{
      $(\tilde d, s, x) \gets $ \hyperref[oracle: prox_opt_solver]{\ProxOpt$(W^\perp,c,d,M,\varepsilon)$} \;
      $I = \{i \in [n]: s_i < 16M^3 n^3\|s_{\Lambda(s,x)}\|_1\}, \; J = [n]\setminus I$\;
      \lIf{$I=\emptyset$}{\Return $(x,s)$}
      \Else{
      \If{\textsc{InnerLoop}$(W_I, x_I, s_I, M,\varepsilon)$ returns a solution}{
          $(w ,z) \gets \textsc{InnerLoop}(W_I, x_I, s_I, M,\varepsilon)$ \;
          \If{$\|L_I^{W^\perp}(z - s_I)\|_\infty \le M\|z - s_I\|_1$ \label{line: inner_loop_lift}}{
            $\tilde x \gets (w, x_J)+d-\tilde d$ \;
            $\tilde s \gets s + L_I^{W^\perp}(z - s_I)$ \;
            \Return $(\tilde x, \tilde s)$ \; 
          } \lElse{
            \Return Lifting certificate of $M < \Le_W$
          }    
      }
      \lElse{
        \Return Lifting certificate of $M < \Le_W$
      }}
  }
  \lElse {
    \Return Lifting certificate of $M < \Le_W$
  }
\end{algorithm}

\iftoggle{focs}{}{\begin{proof}[Proof of \Cref{thm:inner}]
The output of Algorithm \ref{alg:inner_loop}  can be transformed to a solution to the system \ref{eq:Inner} using \eqref{eq:def-tilde}.

We show that the algorithm returns the claimed output and admits the running time bound.
Correctness of the recursive calls follows from \Cref{lem: opt_recursive_call}, and part (iii) also shows that  $I\neq [n]$, i.e., the problem reduces in every call.

 The runtime consists of $n$ calls to the \hyperref[oracle: prox_opt_solver]{\textsl{Prox-Opt-Solver}} which takes $O(\CPU(A) \log(\bar\chi_A + n) + nm^{\omega-1})$ each by \Cref{lem:oracle_2_implementation}. 
 We have to maintain the subspaces $W_I$ throughout as required in form $W_I = \ker(A')$ for some matrix $A'$ in form $A' = A_B^{-1}A$ by \Cref{thm: best_LP_algorithms_real_model} for some basis $B \subset [n]$.
Even though \hyperref[oracle: prox_opt_solver]{\textsl{Prox-Opt-Solver}} is called for the dual subspace $W^\perp_I$, a representation of either subspace suffices due to the symmetry in \Cref{thm: best_LP_algorithms_real_model}.

 Let $A$ denote the original matrix. Initially, we perform Gaussian elimination to obtain a basis $B \subseteq [n]$ and the form $A' = A_B^{-1}A$ in time $O(m^2n)$. At any point in the algorithm let $H \subseteq [n]$ denote the current index subset; we let $W \subseteq \R^H$ denote the current subspace. We also maintain the basis $B \subseteq H$, and the form $A' = A_B^{-1}A_H$; thus $W = \ker(A_H')$. At every iteration we recurse on $I \subset H$ and subspace $W_I$. The matrix representation gets updated by deleting the columns $H \setminus I$. For $I \cap B$, pivot operations are required to expand the set $B \cap I$ to a new full rank basis $B$. This requires $|H \setminus I|$ basis exchanges. Hence the total runtime is $O(n^2m)$ for these operations.
 Recursively we invocate line \ref{line: inner_loop_lift} $n$ times. By \Cref{lem:compute-lift} this takes in total $n \min\{m^2n, n^\omega\}$ time.
\end{proof}}

\subsection{Certificate for the wrong guess}

If the feasibility algorithm is applied to a feasible instance, it always provides a certificate of $M < \Le_W$ if no feasible solution is found. In contrast, \Cref{alg:simple_tardos} only provides the verdict $M < \Le_W$ without providing the corresponding certificate. The {\sc InnerLoop} is able to provide a certificate, but if infeasibility is detected in the calls to the \FeasibilityAlgorithm{} in lines \Cref{alg_tardos: primal_feas} and \Cref{alg_tardos: dual_feas} then this means that the partition $B \cup N$ is wrong. Note that failed calls to \FeasibilityAlgorithm{} could also result in $M < \Le_W$ for which the feasibility solvers provide a certificate. 
We are nonetheless able to provide a certificate of $M < \Le_W$ via repeated application of \Cref{lem:prox-find}.

\begin{lemma}\label{lem: certificate_for_optimization}
In line~\ref{alg_tardos: dual_feas} of
\Cref{alg:simple_tardos}, the dual feasibility solver always succeeds to find a dual solution or finds a certificate of $M < \kappa_W$. If the primal feasibility solver in line~\ref{alg_tardos: primal_feas} fails to find a primal feasible solution, then we can find a lifting 
certificate of $M < \Le_W$ in $O(n^3 m^2)$ time.
\end{lemma}

\iftoggle{focs}{}{\begin{proof}
First, let us show that if the algorithm reaches line~\ref{alg_tardos: dual_feas}, then dual feasibility is guaranteed. By induction, we show that there exists a dual feasible solution in the original subspace $W^0$, supported on $I\cup N$, and that for the current subspace $W\subseteq \R^I$, every dual feasible solution extends to such a solution in $W^0$.

 Initially, $[n]=I$, and we assume that the initial system is both primal and dual feasible. 
Every call to {\sc InnerLoop} delivers a solution $(\tilde d, \tilde x, \tilde s)$ to system \ref{eq:Inner}, supported on the current index set $I$. Here, $\tilde s$ is a dual feasible solution, and $\pr{\tilde x}{\tilde s}=0$. 

The algorithm only moves indices to $B$ where $\tilde x_i>0$, and consequently, $\tilde s_i=0$. 
By induction, $\tilde s$ extends to a dual feasible solution in the original subspace supported on $I\cup N$.
Moreover, the new value of $c\in \R^I$ is set as $\tilde s_I$, verifying the second inductive hypothesis.

\medskip

In contrast, primal feasibility may fail during the algorithm in case of $M<\Le_W$. We now show how to `backtrack' and identify a lifting certificate of $M < \Le_W$. Let $k\le m$ be the number of iterations of the \texttt{while} loop in line \ref{line: while_loop} of \Cref{alg:simple_tardos}.
Let $B_0 = N_0 = \emptyset$, and for $i\in [k]$, let $B_i$ and $N_i$ be the sets $B$ and $N$ at the end of the $i$-th iteration. Similarly, use notation $W^{i}$, $\tilde x^{i}, I^i, I_L^i, I_M^i, I_S^i$ for the spaces, vectors and variables at the end of the $i$-th iteration.

Define $z^{k} = 0_{I_k}$. We iteratively perform following pull-back to obtain $z^i \ge 0$ from $z^{i+1} \ge 0$ for $i = k, \ldots, 0$.
Consider the system $y \ge - \tilde x^i, y \in W^i$. Note that $z^{i} - \tilde x^i$ is a feasible solution. Therefore  \Cref{lem:prox-find} finds in time $O(mn(n-m)^2 + n^2m^{\omega -1}) = O(mn^3)$ a feasible solution $y^i$ to this system such that $\|y\|_\infty \le M \|(\tilde x^i)^-\|_1$ or a certificate $M < \Le_W$. 
If the result is such a certificate, we terminate. Otherwise define 
\begin{equation} \label{eq: define_z_recurse}
    z^i = \tilde x^i + L_{I^i \cup I_S^{0,i} \cup I_M^i}^{W^{i-1}}(y^i, - \tilde x^i_{I_S^{0,i}}, 0_{I_M^i}).
\end{equation}
The lift exists as $y^i \in W^i$ and $I_S^{0,i} \subseteq \cl(I_L^{i})$.
We have $z_{I^i}^i = x_{I^i}^i + y_{I^i}^i \ge 0$, $z^i_{I_S^{0,i}} = 0$, and $z^i_{I_M^i} = \tilde x_{I_M^i}^i \ge 0$.
Furthermore, if the lift does not provide a certificate of $M < \Le_W$, then for $\ell \in I_L^i$ we have
\begin{equation}
 z_{\ell}^i > \frac{\|\tilde x\|_\infty}{2n} 
 - M\|(y^i,\tilde x_{I_S^{0,i}}^i)\|_1 
 \ge \frac{\|\tilde x\|_\infty}{2n}- Mn \max\left\{M \|(\tilde x^i)^-\|_\infty, \frac{\|\tilde x\|_\infty}{4n^2M}\right\} > 0 \, ,
\end{equation}
 where the last inequality follows from the precision $\varepsilon$ in {\sc InnerLoop}. We just have shown that $z^i \ge 0$. Iterating the process would lead to $z^0 \in W + d, z^0 \ge 0$ and $\supp(z^0) \subseteq B$, a contradiction. Therefore either an application of \Cref{lem:prox-find} or a computation as in \eqref{eq: define_z_recurse} must have resulted in a certificate of $M < \Le_W$. 
 The lifts in \eqref{eq: define_z_recurse} can be performed in time $m \min\{m^2n, n^\omega\}$ by \Cref{lem:compute-lift}. Therefore, the runtime is dominated by the repeated application of \Cref{lem:prox-find}, which takes $O(m^2n^3)$.
\end{proof}}

\section{Implementation of oracles and subroutines}\label{sec:approx}

\subsection{A symmetric initialization system}\label{sec:initialization}

Throughout, we let $A\in \mathbb{R}^{m\times n}$, $W=\ker(A)$, and let $b\in\R^m$, $c,d\in \R^n$ such that $Ad=b$.
In order to apply interior point methods to \eqref{LP_primal_dual}, one needs to work with an extended system where a near-central initial solution can be easily obtained---in particular, both primal and dual sides must be strictly feasible.
A common approach is to use the self-dual homogenous initialization \cite{YTM94}; however, this may signficantly increase the condition numbers $\bar\chi_A$ and $\kappa_A$. An alternative initialization that approximately preserves $\bar\chi_A$ and $\kappa_A$ was proposed in \cite{Vavasis1996}, and also used in \cite{DHNV20}. The drawback of this formulation is that the primal and dual side are not symmetric: we would have to prove several properties on the primal and dual side separately. However, as primal and dual are \emph{nearly} symmetric, major parts of the proofs of these lemmas would be identical. 

We now propose a symmetric modification of the system in \cite{Vavasis1996} that also preserves the condition numbers approximately. Throughout, $M$ is an estimate of the $\kappa_W$. Given an instance $(W,c,d)$, we derive two other parameters from $M$.
\begin{equation}\label{eq:M-P-D}
M_P=2\|c\|_1 M,\quad M_D=2\|d\|_1 M\, .
\end{equation} 
 The system is then defined as follows: 
\begin{equation}
    \label{initialization_LP} \tag{\textsl{Init-LP}}
    \begin{aligned} 
    \min \; \pr{c}{x - \ubar x} + &M_P\pr{\1}{\ubar x} \qquad & \max \; \pr{y}{b} - &M_D\pr{\1}{\bar s} \\
    Ax - A\ubar x& = b & A^\top y + s - \bar s &= c \\
    x - \frac{1}{2}\ubar x + \bar x &= M_D\1 &  -A^\top y + \frac{1}{2} \bar s + \ubar s &= M_P\1 - c\\
    x,\bar x,\ubar x& \geq 0 & s,\bar s, \ubar s &\geq 0. \\
    \end{aligned}
\end{equation}
Using the identity
\begin{equation}
   \pr{y}{b} = \pr{y}{Ad} = \pr{A^\top y}{d} = - \pr{s - \bar s}{d} + \pr{c}{d} 
\end{equation}
 an equivalent formulation of primal and dual of \eqref{initialization_LP} is  
\begin{equation}
    \label{new_initialization_LP_subspace} \tag{\textsl{Init-LP-sub}}
    \begin{aligned} 
    \min \; \pr{c}{x - \ubar x} + &M_P \pr{\ubar x}{\1} \qquad & \max \; \pr{d}{c} - &\pr{d}{s - \bar s} - M_D \pr{\1}{\bar s} \\
    x - \ubar x &\in W + d & s - \bar s &\in W^\perp + c \\
    x - \frac{1}{2}\ubar x + \bar x &= M_D\1 & s -\frac{1}{2} \bar s + \ubar s &= M_P\1  \\
    x,\bar x,\ubar x& \geq 0 & s,\bar s, \ubar s &\geq 0. \\
    \end{aligned}
\end{equation}
which displays the desired symmetry via $x \sim s$, $\ubar x \sim \bar s$ and $\bar x \sim \ubar s$.
In the following we will show that the system can be initalized centrally and that the condition number $\bar\chi$ does not increase by too much. Let us begin with the latter.

We denote by $\hat A$ the primal constraint matrix of \eqref{initialization_LP}, that is 
\begin{equation}\label{eq:hat-A}
   \hat A = 
   \begin{pmatrix}
      A & - A & 0 \\
      I & -\frac{1}{2} I & I
   \end{pmatrix}.
\end{equation}
\begin{lemma}
Let $\hat W = \ker(\hat A)$ for $\hat A$ as in \eqref{eq:hat-A}. Then,
$\Le_{\hat W} \le 4\Le_W$.
\end{lemma}
\iftoggle{focs}{}{\begin{proof}
Let $g=(x,y,z)\in \R^{3n}$ denote a minimum support vector in $\hat W$.
First, assume there is an index $i\in [n]$ such that $x_i,y_i,z_i$ all have nonzero values. Then, let $g'=(e_i,e_i,-\frac{1}{2}e_i)$, where $e_i$ is the $i$-th unit vector. We have $g'\in \hat W$, and $\supp(g')\subseteq \supp(g)$. By the minimality of $g$, this implies $g=\alpha g'$ for $\alpha\neq 0$; the ratio between the largest and the smallest absolute value entries is $2$.

For the rest of the proof, we can assume that there is no index $i$ such that $x_i,y_i,z_i$ are all nonzero. By construction, $w=x-y\in W$. Let $T_{xy}\subseteq [n]$ denote the set of indices where $x_iy_i\neq 0$, $T_x\subseteq [n]$ the set with $x_i\neq 0$, $y_i=0$, and $T_y\subseteq [n]$ the set with $y_i\neq 0$, $x_i=0$.
By our assumption, if $i\in T_{xy}$ then $z_i=0$, and therefore $x_i=\frac{1}2 y_i$. If $i\in T_x$ then $z_i=-x_i$, and if $i\in T_y$ then $z_i=\frac{1}2 y_i$. 

We claim that $w$ is a minimum support vector in $W$. Indeed, if there is a smaller support vector $w'\in W$ with $\supp(w')\subsetneq \supp(w)$, then we can map it to $g'=(x',y',z')\in \hat W$ as follows. For each $i\in T_{xy}$, we set $x'_i=2w'_i$, $y'_i=w'_i$. For each $i\in T_x$, we let $x'_i=w'_i$, $z'_i=-w'_i$, and for each $i\in T_y$, we let $y'_i=-w'_i$, $z'_i=-\frac{1}{2}w'_i$. We set all other coordinates of $g'$ to 0. It is easy to verify that $g'\in \hat W$ and $\supp(g')\subsetneq \supp(g)$, giving a contradiction. 

Hence, the largest ratio between the absolute value of elements of $w$ is $\le \kappa_W$. The same construction as described above can be used to map the entries of $w$ to the entries of $g$. This implies a bound $\le 4\kappa_W$ on the ratios between the elements of $g$, since each of $(x_i,y_i,z_i)$ will be one of $w_i,-w_i,2w_i$ and $- \frac{1}{2} w_i$.
\end{proof}}

\paragraph{Initial solutions}
While we use the black box results in \Cref{thm: best_LP_algorithms_bit_model}, it is worth noting that for interior point methods, the system can be easily initialized near the central path with the following solutions.
\begin{equation}\label{eq:init-soln}
   \begin{aligned}
      (x, \ubar x, \bar x ) &= \frac23 M_D(\1_{n}, \1_{n}, \1_{n}) + (d,0_n,-d)\\
      (y, s, \ubar s, \bar s) &= \frac23M_P (0_m, \1_n, \1_n, \1_n) + (0_m,0_n,-\frac12c,-c) \\
   \end{aligned} 
\end{equation}
The duality gap between these solutions is $\approx \frac{4}{3}nM_P M_D$.

\paragraph{Box constraints}
\Cref{thm: best_LP_algorithms_bit_model} requires upper bounds $R_P$ and $R_D$ on the diameters of the primal and dual feasible sets. The formulation  \eqref{initialization_LP} may include unbounded directions. However, we can impose trivial upper bounds 
\begin{equation}\label{eq:box-bounds}
(x,\ubar x,\bar x)\le 2nM_D\1_{3n}\, \quad \mbox{and}\quad (s,\ubar s,\bar s)\le 2nM_P\1_{3n}\, .
\end{equation}
 Any solution violating these bounds would have worse gap than the trivial solution \eqref{eq:init-soln} (see Appendix \ref{app:feas_upper_bounds}); therefore, adding such bounds will not affect optimality.

\subsection[Proofs of the approximate solution theorems]{Proofs of \Cref{thm: best_LP_algorithms_real_model_feasibility} and \Cref{thm: best_LP_algorithms_real_model}}\label{sec:proofs-approx}
In the proofs, we can assume  $d=d/W$ and $c=c/W^\perp$; otherwise, we can replace $d$ and $c$ by their projections to the respective subspaces.
These vectors $d$ and $c$ can be computed as projections in $\min{O(m^\omega, mn^2)}$ time.

\paragraph{The approximate LP setup}
We  use \Cref{thm: best_LP_algorithms_bit_model} to the extended system \eqref{initialization_LP}, with additional box constraints \eqref{eq:box-bounds} (represented via slack variables). Our aim in \Cref{thm: best_LP_algorithms_real_model} is to obtain primal and dual solutions 
$(x,\ubar x,\bar x)$ and $(s,\ubar s,\bar s)$ with duality gap $\le \varepsilon M\|c\|_1 \|d\|_1/2$. In \Cref{thm: best_LP_algorithms_real_model_feasibility} we use it for the special case $c=0$; we now explain the more general case.

We note that, instead of using the final output of \Cref{thm: best_LP_algorithms_bit_model}, one could obtain more direct algorithms by using the interior-point methods of \cite{LS19, vdb20, br2020solving} directly on our extended system. We now give a black-box argument using \Cref{thm: best_LP_algorithms_bit_model}, to demonstrate the compatibility of our results with any approximate LP solver, not just interior-point methods.

Both \Cref{thm: best_LP_algorithms_real_model_feasibility} and \ref{thm: best_LP_algorithms_real_model} require in the input $A_B^{-1} A$ for a nonsingular submatrix $A_B$.
For simplicity of notation, we assume $A$ is already in such a form, with $B=[m]$, that is, $A=(I_m|T)$ for some $T\in\R^{m\times(n-m)}$. (Such a form could be obtained via Gaussian elimination, but that might take more time than a single call of the approximate algorithm itself.) Also, we replace the terms $M_P$ and $M_D$ with $\hat M_P = M_P - \varepsilon \|c\|_1 $ and $\hat M_D = M_D - \varepsilon \|d\|_1$ respectively. The motivation behind this shift will become clear at the end of the paragraph.

The extended system will be of the form $\min \pr{\hat c}{\hat x}$ s.t. $\hat A \hat x=\hat b$, $\hat x\ge 0$ with $\hat A$ as in \eqref{eq:hat-A}; we ignore the box constraints for the simplicity of presentation.  

Let us now bound the parameters.  We have the diameter bounds $R_P=O(n^{3/2} M_D)=O(n^{5/2} M\|d\|_1)$ and $R_D=O(n^{3/2} M_P)=O(n^{5/2} M\|c\|_1)$.
Using the special form of $A$, by \Cref{prop:kappa-max} we have $\|A\|_F\le n \|A\|_{\max}\le n\kappa_W$.
For the matrix of the extended system \eqref{eq:hat-A}, we also have $\|\hat A\|_F\le O(n\kappa_W)=O(nM)$, if our guess $M$ on $\kappa_W$ was correct.
Otherwise, we obtain a lifting certificate of $M<\kappa_W$ that can be derived from one of the circuits in the form of $A$.
 
We next bound $\|\hat b\|$ and $\|\hat c\|$. We have $\|b\|_\infty\le \|A'\|_{\max} \|d\|_1\le n\kappa_W\|d\|_1$. If our guess $M>\kappa_W$ is correct, then this is upper bounded  as $nM\|d\|_1$; otherwise, we obtain a lifting certificate of $M<\kappa_W$. The vector $\hat b$ also contains additional entries of $\hat M_D=(M - \varepsilon)\|d\|_1$; thus, we get $\|\hat b\|_2\le O(n^{3/2} M\|d\|_1)$. Similarly, we obtain $\|\hat c\|_2 \le O(n^{3/2} M\|c\|_1)$. 

The output of \Cref{thm: best_LP_algorithms_bit_model} maintains nonnegativity but allows for violation of the subspaces constraints $\hat A \hat x=\hat b$ and $\hat A^\top\hat y+\hat s=\hat c$.
In  \Cref{thm: best_LP_algorithms_real_model_feasibility} and \Cref{thm: best_LP_algorithms_real_model}, we require the points in the right subspaces, but allow for violation of nonnegativity. We now show how fix the violations in the subspaces to obtain a proximal feasible solution to \ref{initialization_LP}.

We will use the algorithm in \Cref{thm: best_LP_algorithms_bit_model} with $\delta=\gamma\varepsilon n^{-4.5} M^{-2}$ for a suitable constant $\gamma >0$, so that the optimality gap is bounded by at most $\gamma\varepsilon M\|c\|_1 \|d\|_1$,  
$\|\hat A\hat x-\hat b\|_2$ is bounded by $\gamma \varepsilon \|d\|_1/n$, and 
$\|\hat A^\top\hat y+\hat s-\hat c\|$ is bounded by $\gamma \varepsilon \|c\|_1/n$. We let $\hat x=(x',\ubar x',\bar x')$, $\hat y$, and $\hat s=(s',\ubar s',\ubar x')$ denote the primal and dual solutions. 

Recalling the form $A=(I_m|T)$, we can subtract $|A_i(x-\ubar x)-b_i| \le \gamma\varepsilon\|d\|_1/n$ from either $x'_i$ or $\ubar x'_i$ to obtain a feasible solution to the equality system. By subsequently shifting up $x'$ and $\ubar x'$ by $2\gamma \varepsilon \|d\|_1 /n \1$ each, we maintain feasibility of the equality system and achieve nonnegativity of $x'$ and $\ubar x'$. Additionally, the equality $x - \frac12 \ubar x' + \bar x' = M_D \1$ is satisfied. Thus the solution is feasible to the extended primal system \eqref{initialization_LP}. 

The total adjustment to the original solution amounts to $2\gamma\varepsilon\|d\|_1/n$ per coordinate at most. Therefore, the objective value changed by at most 

\[3n \cdot \|\hat c\|_\infty \cdot 2\gamma \varepsilon \|d\|_1/n = 6 \gamma\varepsilon \hat M_P\|d\|_1 < 6 \gamma \varepsilon M \|c\|_1 \|d\|_1.\] 

Adjusting the dual can be done equivalently, noting the symmetry in \eqref{new_initialization_LP_subspace}. Note that this requires a matrix $B \in \R^{(n-m) \times n}$ such that $\im(B) = W^\perp$, which might be computed by Gram-Schmidt orthogonalisation. This could take longer than a single call of the approximate algorithm itself, therefore we assume that this matrix is given. 

All together, we obtain feasible primal and dual solution $(x,\ubar x,\bar x)$, $y$, $(s,\ubar s,\bar s)$ to \eqref{initialization_LP} with objective values within $\varepsilon M\|c\|_1 \|d\|_1$ if $\gamma \le 1/13$ is chosen. The running time can be bounded as
\[
O(\CPU(A)\log(n/\delta))=O(\CPU(A)\log((M+n)/\varepsilon))\, .
\]

\iftoggle{focs}{}
{
\theoremFeasibilityProx*
}

\iftoggle{focs}{}{\begin{proof}[Proof of \Cref{thm: best_LP_algorithms_real_model_feasibility}]
If $d=0$, then we return $x=0$. Otherwise,
we use the approximate solver as described above, for systems of the form \eqref{initialization_LP} with $c=0$, $M_P=1$, and $M_D=2M\|d\|_1$.  We thus obtain solutions 
  $(x, \ubar x, \bar x )$ and   $(y,s, \ubar s, \bar s)$ with duality gap $\le \varepsilon  \|d\|_1$. Assume the dual objective value is $\le 0$. In this case, the primal objective value is $\pr{\1}{\ubar x}\le \varepsilon  \|d\|_1$. Setting $\hat x=x-\ubar x$, we get the required solution
\begin{equation}\label{eq:primal-approx}
  \hat x\in W+d,\quad \|\hat x^-\|_1\le {\varepsilon}\|d\|_1,\quad \|\hat x\|_\infty\le 2M\|d\|_1\, .
   \end{equation}
  Next, assume the dual objective value is $-\pr{d}{s - \bar s} - M_D \pr{\1}{\bar s} >0$ for vector $s,\bar s\ge 0, s - \bar s \in W^\perp$. Then $\hat s=s-\bar s$ satisfies
  \begin{equation}\label{eq:approx-Farkas}
\pr{d}{\hat s}<-2M\|d\|_1\|\hat s^-\|_1\, .
  \end{equation}
We now apply \Cref{lem:slow-lift-prox} to $W^\perp$, the vector $\hat s$, and $J=\supp(\hat s^-)$ to find a vector $z\in W^\perp$ that is sign-consistent with $\hat s$, $z_J=0$, and $\|z-\hat s\|_\infty\le M\|\hat s^-\|_1$ in $O((n-m)m^2 + n(n-m)^{\omega - 1}) = O((n-m)m^2 + n^\omega)$ time. If no such vector is found, then we obtain a lifting certificate of $M<\kappa_{W^\perp} = \kappa_W$ as in \ref{case:certificate-feasible}.

By the sign-consistency, $z\ge 0$, since $\hat s_i\ge 0$ for all $i\in [n]\setminus J$. We show that $\pr{d}{z}<0$, and consequently, $z$ is a Farkas certificate as in \ref{case:feasible-farkas}. This follows since
\[
\pr{d}{z}= 
\allowbreak
\pr{d}{\hat s}+\pr{d}{z-\hat s}
\le 
\allowbreak 
\pr{d}{\hat s}+\|d\|_1 \|z-\hat s\|_\infty\, 
< 
\allowbreak 
 -2M\|d\|_1\|\hat s^-\|_1 + M\|d\|_1\|\hat s^-\|_1 \le 0. \qedhere
 \]
\end{proof}}

\iftoggle{focs}{}
{
\theoremOptimalityProx*
}

\iftoggle{focs}{}{\begin{proof}[Proof of \Cref{thm: best_LP_algorithms_real_model}]

First, apply \Cref{thm: best_LP_algorithms_real_model_feasibility} 
to data $(W, d, \varepsilon/4, M)$ and $(W^\perp, c, \varepsilon/4, M)$. If either returns in \ref{case:certificate-feasible} we return in \ref{case:certificate-guess-wrong}.
If the first call returns in \ref{case:feasible-farkas}, we return in \ref{case:certificate-farkas-primal}, if the second call returns in \ref{case:feasible-farkas}, we return in \ref{case:certificate-farkas-dual}. Note that the claimed runtime bounds are observed. From now on assume that both calls returned in \ref{case:near-feasible}. Then we obtain vectors $\hat x$ and $\hat s$ such that

\begin{equation}\label{eq:hat-x-hat-s-output}
  \begin{aligned}
  &\hat x\in W+d\, ,& &\|\hat x^-\|\le \frac{\varepsilon}4 \|d/W\|_1\, , & \quad &\|\hat x\|_\infty\le 2M\|d\|_1\, ,\\ 
  &\hat s\in W^\perp+c\, ,& &\|\hat s^-\|\le \frac{\varepsilon}4 \|c/W^\perp\|_1\, , &\quad  &\|\hat s\|_\infty\le 2M\|c\|_1\, .
  \end{aligned}
\end{equation}

Proceed by calling an approximate interior point solver as described above, with $M_P$ and $M_D$ defined as in \eqref{eq:M-P-D}.
 In time 
$O(\CPU(A)\log((M+n)/\varepsilon))$ we can obtain solutions 
  $(x, \ubar x, \bar x )$ and   $(y,s, \ubar s, \bar s)$ with normalized duality gap $\le \varepsilon M \|c\|_1\|d\|_1/2$.

  We let $\tilde x = x - \ubar x$ and $\tilde s = s - \bar s$. These fulfill $\tilde x \in W + d$ and $\tilde s \in W^\perp + c$. We distinguish two cases.

  \paragraph{Case I: $\|\tilde x^-\|_1\le \varepsilon\|d\|_1$ and $\|\tilde s^-\|_1\le \varepsilon\|c\|_1$.} In this case, we claim that $(\tilde x,\tilde s)$ is a near-feasible near-optimal pair of solutions as required in \ref{case:near-feasible-optimal}. 
First, we have upper bounds $\|\tilde x\|_\infty\le 2M\|d\|_1$, $\|\tilde s\|_\infty\le 2M\|c\|_1$ from
\begin{equation}\label{eq:tilde-bounds}
   \begin{aligned}
      \tilde x &= x - \ubar x \le x - \frac12 \ubar x + \bar x = 2M\|d\|_1 \1 \\ 
      \tilde s &= s - \bar s \le s - \frac12 \bar s + \ubar s = 2M \|c\|_1 \1\, .
   \end{aligned} 
\end{equation}
We now verify the bound on $\|x\circ s\|$. With the upper bounds on $\tilde x$ and $\tilde s$
\[
\begin{aligned}
\|\tilde x \circ \tilde s\|_1 &\le  \|x \circ s\|_1 +\|\tilde x^- \circ \tilde s^-\|_1+\|\tilde x^- \circ \tilde s\|_1 + \|\tilde x \circ \tilde s^- \|_1 \\
 &\le \left(\frac{\varepsilon}{2} M + \varepsilon^2+4\varepsilon M\right)  \|c\|_1\|d\|_1\le 5\varepsilon M \|c\|_1\|d\|_1\, ,
\end{aligned}
\]
using the assumption $\varepsilon \le 1/2$.

Finally, consider \eqref{eq:lem_statement_on_opt}. If $\kappa_W\ge M$, and  $x\in W+d$, $x\ge 0$ is feasible and bounded, then there exists an optimal solution $x^*$ with $\|x^*\|_\infty\le M_D$. We can map this to a solution of the extended system as $(x^*,0,M_D\1_n-x^*)$, with the same objective value $\textrm{OPT}=\pr{c}{x^*}$. 
By weak duality on the extended system, this is lower bounded by the dual objective value of $(\tilde s^+, \tilde s^-, M_P\1_n-\tilde s^++\frac{1}{2}\tilde s^-)$, that is $\pr{d}{c-\tilde s}-M_D \|\tilde s^-\|_1$. The bound $\textrm{OPT}\le \pr{d}{c-\tilde s}+5\varepsilon M 
\|c\|_1\|d\|_1$ follows by this assumption and the bound on the duality gap. The primal bound follows analogously.

 \paragraph{Case II: $\|\tilde x^-\|_1> \varepsilon\|d\|_1$ or $\|\tilde s^-\|_1> \varepsilon\|c\|_1$.} In this case, we use the approximately feasible solutions $\hat x$ and $\hat s$ to find a lifting certificate of $M>\kappa_W$ or $M>\kappa_{W^\perp}$.

Assume that $\|\tilde x^-\|_1> \varepsilon\|d\|_1$. Let us define $\ell=\tilde x-2M\|d\|_1\1$, and $u=\tilde x+\hat x^-$. Consider the system
\[
v\in W\, , \quad \ell\le v\le u\, .
\]
The vector $\tilde x-\hat x$ is a feasible solution. We now apply \Cref{lem:prox-find} to find a lifting certificate of $M>\kappa_W$ or another vector $v$ feasible to the above system such that 
\[
\|v\|_\infty\le M\|\ell^++u^-\|_1=M \|u^-\|_1\le M\|\tilde x^-\|_1\, ,
\]
 since $\ell\le 0$ by \eqref{eq:tilde-bounds}, and $|u_i|\le \tilde x^-_i$ whenever $u_i<0$. We let $z=\tilde x-v$. Now,
\[
z\in W+d,\, \quad -\hat x^-\le z\le 2M\|d\|_1\1,\, \quad \|z-\tilde x\|_\infty\le M\|\tilde x^-\|_1\, .
\]
In particular, $\|z^-\|_1\le \|\hat x^-\|_1\le \frac{\varepsilon}4 \|d\|_1$.
We can map $z$ to a primal feasible solution $(z^+,z^-,M_D\1-z^++\frac{1}{2}z^-)$
of \eqref{initialization_LP}. The objective value of this solution is
\[
\begin{aligned}
\pr{c}{z}+2M\|c\|_1\|z^-\|_1&\le \pr{c}{\tilde x}+\pr{c}{z-\tilde x}+ \frac{\varepsilon}2 M \|c\|_1 \|d\|_1\\
&\le \pr{c}{\tilde x}+\|c\|_1 \|z-\tilde x\|_\infty+\frac{\varepsilon}2 M \|c\|_1 \|d\|_1\\
&\le \pr{c}{\tilde x}+M\|c\|_1 \|\tilde x^-\|_1+\frac{\varepsilon}2 M \|c\|_1 \|d\|_1\\
&< \pr{c}{\tilde x}+2M\|c\|_1 \|\tilde x^-\|_1 - \frac{\varepsilon}2 M \|c\|_1 \|d\|_1 \\
&\le \pr{c}{x - \ubar x} + 2M\|c\|_1\|\ubar x\|_1 - \frac{\varepsilon}2 M \|c\|_1 \|d\|_1\, , 
\end{aligned}
\]
where the penultimate inequality used the assumption $\|\tilde x^-\|_1>\varepsilon \|d\|_1$.
This is a contradiction, since $\pr{c}{x - \ubar x}+2M\|c\|_1 \|\ubar x\|_1$ is the objective value of the solution $(x,\ubar x,\bar x)$ that was chosen to be within
$\frac{\varepsilon}2 M \|c\|_1 \|d\|_1$ from the optimum value.

Hence, such a feasible solution does not exist, and therefore, the algorithm in \Cref{lem:prox-find} must have terminated with a lifting certificate of $M<\kappa_W$ as in \ref{case:feasible-farkas}. The analogous argument applies for the case $\|\tilde s^-\|_1> \varepsilon\|c\|_1$.
The time complexity bound follows: the additional running time term $O(n^3m + n^{\omega + 1})$ comes from \Cref{lem:prox-find}, and is only needed in Case II, that leads to \ref{case:certificate-guess-wrong}. When applied to the primal we get $O(nm(n-m)^2 + nm^{\omega - 1})$ as in \Cref{lem:prox-find}. When applied to the dual, the rank is $n-m$ and so the runtime is $O(nm^2(n-m) + n(n-m)^{\omega - 1})$. Put together, it results in runtime $O(n^3m + n^{\omega}) = O(n^3m)$.
\end{proof}}

\subsection{Implementation of the oracles}
\label{subsec: oracle_2_implementation}

We proceed to show  \Cref{thm:approx-implement} and \Cref{lem:oracle_2_implementation} on the implementations of Oracle \ref{oracle: prox_feas_solver} and Oracle \ref{oracle: prox_opt_solver}, respectively. We restate both statements and oracles for convenience.

\iftoggle{focs}{}{
\approxImplement*
}

\iftoggle{focs}{}{
\oracleTwoImplementation*
}

\setcounter{oraclecf}{0}

\iftoggle{focs}{}
{
\oracleOne*
\oracleTwo*
}

\iftoggle{focs}{}{
We note that 
Oracle \ref{oracle: prox_feas_solver} is a special case of Oracle \ref{oracle: prox_opt_solver}; thus, \Cref{thm:approx-implement} follows as a special case of \Cref{lem:oracle_2_implementation}. We now give the proof of the latter statement.
\begin{proof}[Proof of  \Cref{lem:oracle_2_implementation}]

Let us define
\[
\tau=\|d_{\Lambda(d,c)}\|_1, \, I: = \set{i \in [n] : d_i \le 2M\tau}, \quad 
\textrm{and} \quad \varepsilon' = \frac{\varepsilon^2}{28M^3 n^3}\, .
\]

We run the algorithm in \Cref{thm: best_LP_algorithms_real_model} on the system \hyperref[LP-subspace-f]{\textsl{Primal-Dual}$(\pi_I(W), d_I, c_I)$} and with $\varepsilon'$. Computing the corresponding constraing matrix can be done in $O(nm^{\omega - 1})$. The reason that we do not consider the elements in $[n]\setminus I$ is that we need to bound the minimum-norm vector in the affine space in which we run the solver. Note that $\|d_I/\pi_I(W)\|_1\le 2n^{3/2}M\tau$ and that $\supp(c) \subseteq I$ as $c \ge 0$.
Let us discuss the potential outcomes of the algorithm. We can only be in case \ref{case:certificate-farkas-primal}, if the original system $x \in W + d$ was already infeasible. In particular, a Farkas certificate extends to this original system. 
\ref{case:certificate-farkas-dual} can not happen as $c_{[n] \setminus I} = 0$ by definition of $I$ and so $c_I$ is a feasible dual solution to \hyperref[LP-subspace-f]{\textsl{Primal-Dual}$(\pi_I(W), d_I, c_I)$}. 
If in case
 \ref{case:certificate-guess-wrong}, we output the corresponding lifting certificate of $M< \kappa_W$ and terminate.

Assume now we obtained approximately feasible and approximately optimal primal and dual solutions $(\hat x, \hat s)$ as in \ref{case:near-feasible-optimal} fulfilling  

\begin{equation}
  \label{eq: primal_dual_near_optimal_implementation}
  \begin{aligned}
      \|\hat x^-\|_1 &\le \varepsilon' \|d_I/\pi_I(W)\|_1 \le \varepsilon' \cdot 2n^{3/2}M\tau , \\ 
      \|\hat s^-\|_1  &\le \varepsilon' \|c/W^\perp\|_1,  \\
      \|\hat x \circ \hat s\|_1 &\le 10\varepsilon' M^2 n^{3/2} \tau \|c/W^\perp\|_1, \\
      \|\hat x\|_\infty & \le 4nM^2 \tau, \; \text{and}  \\
      \|\hat s\|_\infty & \le 2M\|c/W^\perp\|_1.
  \end{aligned}
\end{equation}
Further, extend $\hat x$ arbitrarily to an element $x' \in W + d$ with $x'_I = \hat x$ and set $s' = (0_{[n]\setminus I}, \hat s)$. 

We proceed by converting this solution into another one that is proximal to $d$ without increasing the objective value respective to $c$. To this end consider, as in the proof of  \Cref{thm:prox-basic}, a system $z\in W$,
$\ell\le z\le u$ with 
\[ 
\ell = - d - (x')^-
\quad \text{and} \quad
u_i = 
\left\{
\begin{array}{ll}
  x'_i - d_i & \text{if } i \in \supp(c), \\
 \infty &\text{else.}
\end{array}
\right.
\]
We have that 
\[\|\ell^+ +u^-\|_1\le \|d_{\Lambda(d,c)}\|_1 + \|(x')^-\|_1 \le \tau + \varepsilon' \cdot 2Mn^{3/2} \tau \le 2\tau,\] 
As $x'-d$ is a feasible solution we can apply \Cref{lem:prox-find}. 
 Assume it can find a feasible solution $z$ with $\|z\|_\infty \le M\|\ell^+ + u^-\|_1 \le 2M\tau$, as otherwise we terminate with a certificate of $M < \kappa_W$. Defining $\tilde x= d+z$, we note that $\tilde x_{[n]\setminus I} \ge 0$ and $\pr{c}{\tilde x} = \pr{c}{d + z}\le \pr{c}{x'}$ hold as $c\ge 0$ and $d_i + z_i \le x'_i$ for $i \in \supp(c)$.

On the dual side we turn the proximal feasible vector $s'$ into a feasible vector, using the feasible vector $c \ge 0$, again by applying \Cref{lem:prox-find}. Consider the system 
\begin{equation}
w \in W^\perp \cap \R_I^n, w \ge -s', 
\end{equation}
of which $c - s'$ is a feasible solution. Therefore we are able to find a solution $w$ with $\|w\|_\infty \le M\|(s')^-\|_1 = M \|\hat s^-\|_1 \le \varepsilon' M \|c/W^\perp\|_1$ or terminate with a certificate $M < \kappa_W$. In the former case, let $\bar s = s' + w \ge 0$ and define $\tilde s$ as 

\[
\tilde s_i = \begin{cases} 0 &\textrm{if }  \bar s_i \le \frac{\varepsilon}{n} \|c/W^\perp\|_1 \, ,\\
 \bar s_i &\textrm{otherwise}.
\end{cases}
\]

Note that $\bar s \ge \tilde s \ge 0$ and so setting $\tilde c = c - \bar s + \tilde s$ gives 
$c - \tilde c \ge 0$, $\tilde s \in W^\perp + \tilde c$, and
\[
\|\tilde c - c\|_\infty = \|\tilde s - \bar s\|_\infty \le 
\frac{\varepsilon}{n} \|c/W^\perp\|_1.
\]
Further we have 
\begin{equation}
  \label{case: oracle_2_bounds_new}
  \begin{aligned}
  \|\tilde x^-\|_1 \le \|\hat x^-\|_1 \le \varepsilon' \cdot 2Mn^{3/2}\tau, \text{ and } 
   \|\tilde s\|_\infty \le \|\bar s\|_\infty \le \|\hat s\|_\infty + \|w\|_\infty \le 3M \|c/W^\perp\|_1. 
  \end{aligned} 
\end{equation}
To bound the duality gap beween $\tilde x$ and $\tilde s$, we note that $\pr{\tilde x}{\bar s} \le \pr{x'}{\bar s}$ as $\pr{\tilde x}{c} \le \pr{x'}{c}$ and $c - \bar s \in W^\perp$. Also recall that $\tilde s - \bar s \le 0$. Therefore,
\begin{equation}
\begin{aligned}
\pr{\tilde x}{\tilde s} 
&= \pr{\tilde x}{\bar s} + \pr{\tilde x}{\tilde s - \bar s} \\
&\le \pr{x'}{\bar s} + \|\tilde x^-\|_1 \|\tilde s - \bar s\|_\infty \\
&= \pr{x_I'}{s_I'} + \pr{x_I'}{\bar s_I - s_I'} + \|\tilde x^-\|_1 \|\tilde s - \bar s\|_\infty  \\
&= \pr{\hat x}{\hat s} + \pr{\hat x}{w_I} + \|\tilde x^-\|_1\|\tilde s- \bar s\|_\infty  \\
& \le \|\hat x \circ \hat s\|_1 + \|\hat x\|_1 \|w\|_\infty + \|\tilde x^-\|_1\|\tilde s- \bar s\|_\infty \\
&\le 10\varepsilon' M^2 n^{3/2} \tau \|c/W^\perp\|_1 + 4n^2M^2 \tau \cdot \varepsilon'M\|c/W^\perp\|_1 + 2\varepsilon'Mn^{3/2}\tau \cdot \frac{\varepsilon}{n} \|c/W^\perp\|_1
\\
& \le 16 \varepsilon' M^3 n^2 \tau \|c/W^\perp\|_1
\end{aligned}
\end{equation}
and by 
\[\pr{\tilde x^-}{\tilde s} \le \|\tilde x^-\|_1 \|\tilde s\|_\infty \le 6\varepsilon' M^2n^{3/2}\tau \|c/W^\perp\|_1,\]
we also get
\begin{equation}
  \label{eq:duality_gap_tilde}
  \|\tilde x \circ \tilde s \|_1 = \pr{\tilde x}{\tilde s} + 2 \pr{\tilde x^-}{\tilde s} \le 28 \varepsilon' M^3 n^2 \tau \|c/W^\perp\|_1, 
\end{equation}
where the equality used $\tilde s \ge 0$.
Let us verify the requirements of \ref{eq:ProxOpt}. The bounds on $c - \tilde c$ were shown already.
Using the bound on the norm of $z$ we get the primal proximity as 
\[\|\tilde x - d\|_\infty = \|z\|_\infty \le 
2 M \tau \le 3M^2n\tau\, .\] 

Let us now turn to the bound on $\|\tilde x_{\Lambda(\tilde x,\tilde s)}\|_\infty$.
First, let $j\in \supp(\tilde x^-)$. As shown above we have $j \in I$ and further $\tilde x_j \ge x_j'$ 
 and therefore
\[|\tilde x_j|\le |(x'_j)^-|\le \|(x_I')^-\|_1 = \|\hat x^-\|_1\le  \varepsilon' \cdot 2Mn^{3/2}
\tau< \varepsilon\tau\, .\]

Next, let $j\in \supp(\tilde s^+)$. Using \eqref{eq:duality_gap_tilde} we obtain
\begin{equation}
  |\tilde x_j| \le \frac{\|\tilde x \circ \tilde s\|_1}{\tilde s_j} \le \frac{28\varepsilon' M^3 n^2 \tau \|c/W^\perp\|_1}{\frac{\varepsilon}{n}\|c/W^\perp\|_1} \le \varepsilon \tau,
\end{equation}
by choice of $\varepsilon'$.
We conclude that all inequalities of \ref{eq:ProxOpt} are satisfied.
\end{proof}}

\bibliographystyle{alpha}
\bibliography{bib}

\end{document}